\documentclass[12pt,a4paper]{amsart}   

\usepackage[T1]{fontenc}
\usepackage[cp1252]{inputenc}            
\usepackage[dvips]{graphicx}
\usepackage{amssymb}

\newtheorem{lause}{Theorem}               
\newtheorem{maar}[lause]{Definition}    

{\par\vspace{2.5mm}\noindent\textit{Proof.\ }}
{\qed\par\vspace{2.5mm}}

\numberwithin{equation}{section}
\numberwithin{lause}{section}

\allowdisplaybreaks
\begin{document}
\thispagestyle{empty}


\title{Periodic orbits of quadratic polynomials of periods six and seven}

\author{Pekka Kosunen}

\date{\today}

\subjclass[2010]{Primary 37F10; Secondary 30D05}

\keywords{Periodic orbit, iteration, quadratic polynomial, eigenvalue, critical point}

\maketitle

\noindent
\address{Department of Physics and Mathematics, University of Eastern Finland, P. O. Box 111, FI-80101 Joensuu, Finland;}
\email{pekka.kosunen@uef.fi}\\
\thanks{Supported in part by the Academy of Finland grants no. 286877.}

\begin{abstract}
The dynamics of quadratic polynomials is commonly studied by using the family of maps $f_{c}(x)=x^{2}+c$, where $c\in \mathbb{C}$. In this paper we form equations of periodic orbits of periods six and seven on a new $(u,v)$-plane and consider also the corresponding equations on the $(x,y)$-plane. The new $(u,v)$-model produces equations for the periods six and seven with significantly lower degree than the ones obtained by the previous models which enables us to find their solutions.
\end{abstract}

\section{Introduction}
The dynamics of quadratic polynomials is often studied by using the family of maps  $f_{c}(x)=x^{2}+c$, where $c\in \mathbb{C}$ and $x_{i+1}=f_{c}(x_{i})=x_{i}^{2}+c$ (see, for example, \cite{4}, \cite{13} and \cite{25}). In this work a central role is played by the periodic orbits or cycles of $f_{c}$ for which $f_{c}^{n}(x_{0})=x_{0}$ for some $n\in\mathbb{N}$. Now the number $n$ is the period of the orbit and $x_{0},x_{1},x_{2},\ldots x_{n-1}$ are periodic points of period $n$.

In this paper, we obtain the equations for the periodic orbits $6$ and $7$ of the family $f_{c}$ by iterating the function ($\ref{1}$) below and forming the corresponding iterating system on the $(u,v)$-plane by using a model introduced in \cite{5} (see Section $2$). The function $G$ in ($\ref{1}$) is a two-dimensional quadratic polynomial map which is defined in the complex 2-space $\mathbb{C}^{2}$ and its iteration reveals the dynamics of $f_{c}$. We will also compare the $(u,v)$-plane model to the $(x,y)$-plane model, which was introduced by Erkama in \cite{1}. The purpose of developing this new model in \cite{5} was that using it we obtain equations of periodic orbits of lower degree, which are easier to handle. The solutions of the periods $6$ and $7$ obtained in the current paper is a good example of the usefulness of this method.

 When $x$ is the periodic point of the period $n$ in the $(x,c)$-plane, then $f_{c}^{n}(x)=x$. Correspondingly, when $(u,v)$ is the periodic point of the period $n$ in the $(u,v)$-plane, then $G^{n}(u,v)=(u,v)$. The one significant difference between these previous two models is that the $(u,v)$-plane model describes information about both the parameter and dynamical spaces for quadratic polynomial dynamics in a single dynamical system on $\mathbb{C}^{2}$. Also Erkama's model (\cite{1}) can be used to study the dynamics of $f_{c}$ very effectively. In the same way, when $(x,y)$ is the periodic point of the period $n$ in the $(x,y)$-plane, then $F^{n}(x,y)=(x,y)$, where $F(x,y)=(y,y^{2}+y-x^{2})$.

\begin{figure}[h!]
\begin{center}
\includegraphics[width=0.57\textwidth,angle=360]{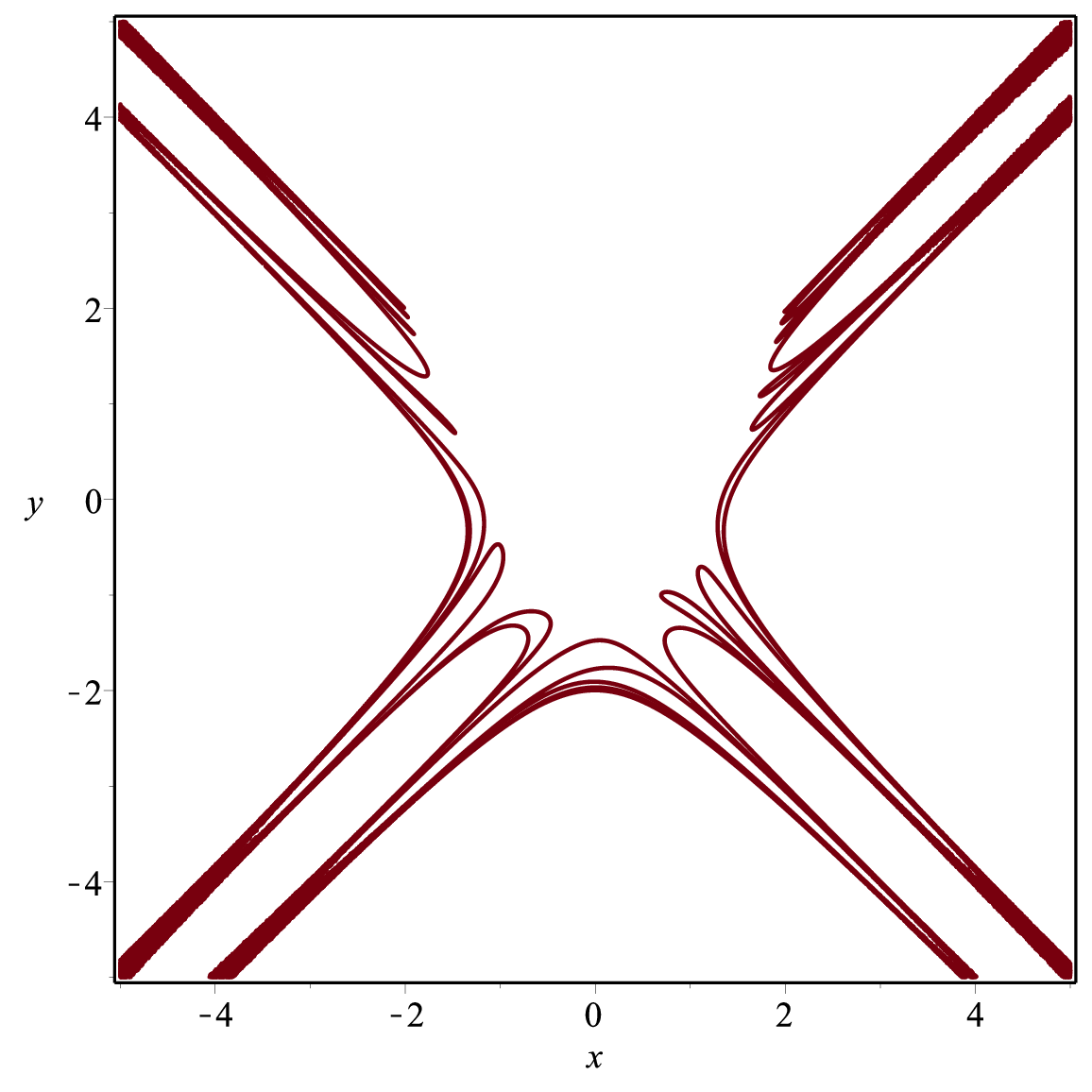}
\caption{Curves of the period six orbits on the $(x,y)$-plane.} \label{kuva1}
\end{center}
\end{figure}

\begin{figure}[h!]
\begin{center}
\includegraphics[width=0.57\textwidth,angle=360]{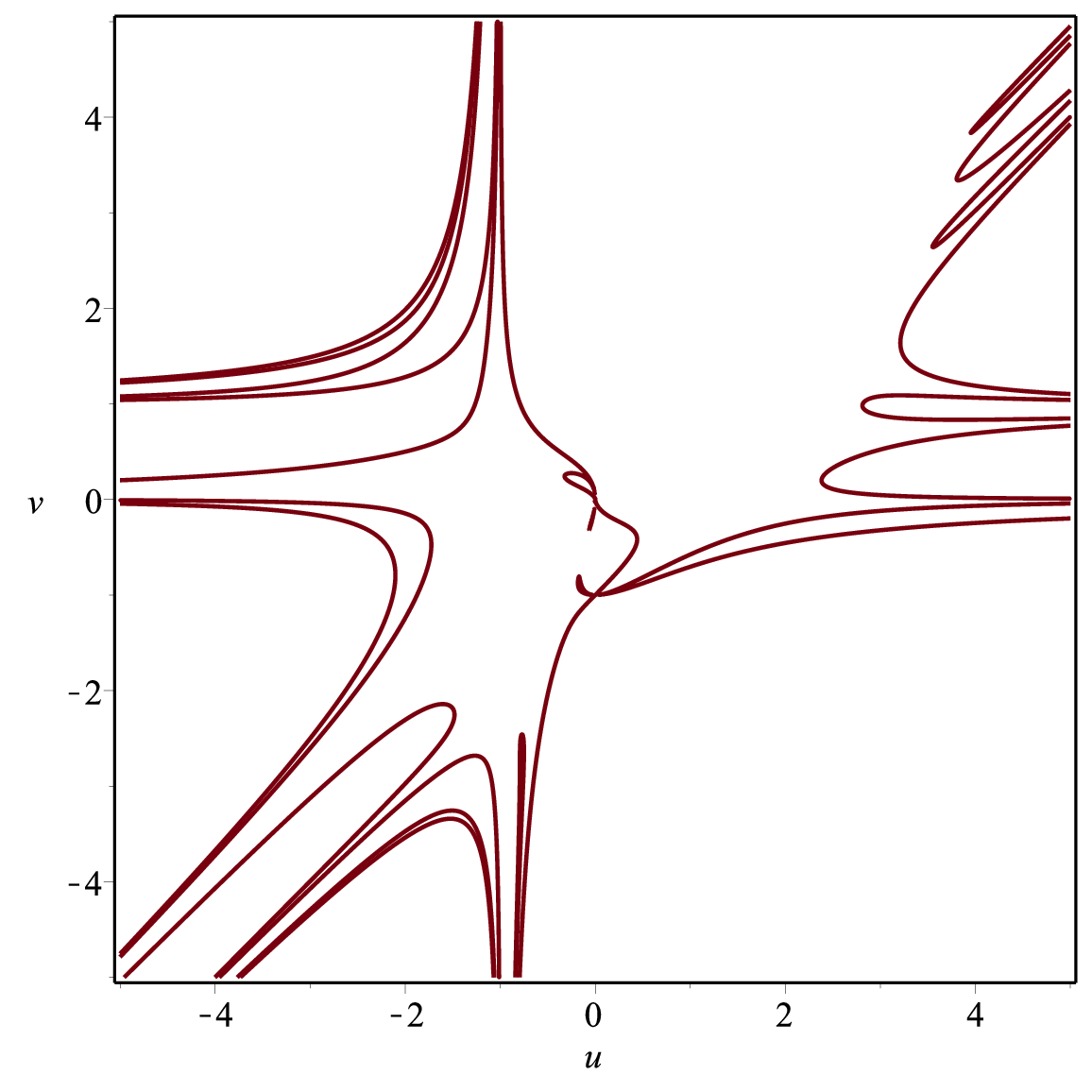}
\caption{Curves of the period six orbits on the $(u,v)$-plane.} \label{kuva2}
\end{center}
\end{figure}

\begin{figure}[h!]
\begin{center}
\includegraphics[width=0.6\textwidth,angle=360]{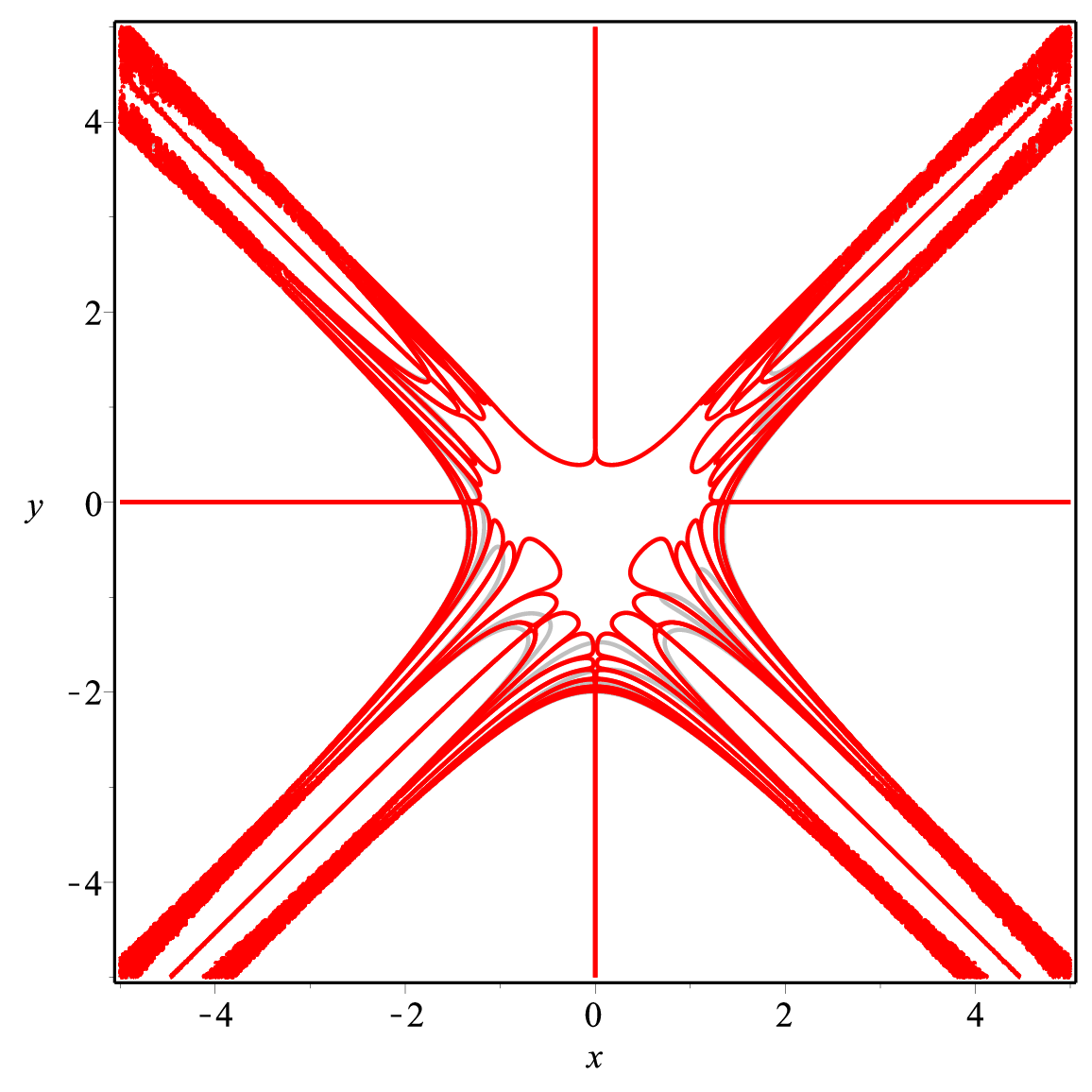}
\caption{Curves of the period six orbits and eigenvalues $|\lambda_{6}(x,y)|=1.$}
\label{kuva3}
\end{center}
\end{figure}

\begin{figure}[h!]
\begin{center}
\includegraphics[width=0.6\textwidth,angle=360]{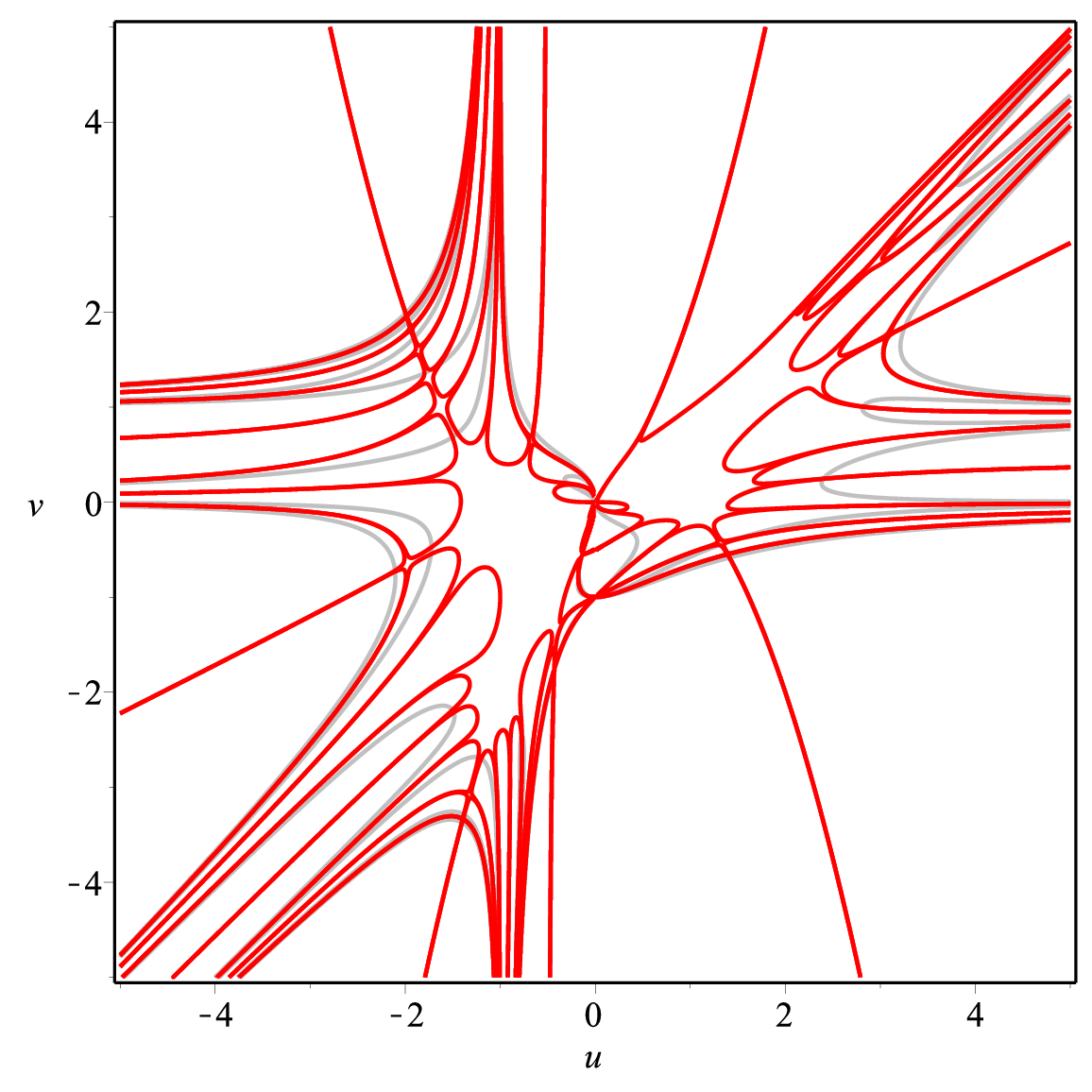}
\caption{Curves of the period six orbits and eigenvalues $|\lambda_{6}(u,v)|=1.$}
\label{kuva4}
\end{center}
\end{figure}

 Previously the equation of period four orbits has been presented in the literature before the article \cite{5} twice; T. Erkama \cite{1} presented nearly the same form as in \cite{5} at $2006$ and Morton \cite{2} presented a corresponding solution first, but in a more complicated form at $1998$. In the article \cite{5} we obtained equations of period $1-5$ orbits on the $(u,v)$-plane by iterating the function $G$. This equation of the period five orbits in the article \cite{5} is the first explicit presentation of the model with two variables.
In the articles \cite{4} and \cite{2} Brown and Morton have also studied period five and six cases and they have formed the so called \textit{trace} formulas in both cases using $c$ and the sum of period cycle points as parameters. Respectively E. V. Flynn, Bjorn Poonen and Edward F. Schaefer have studied the same parameter formula in the case of period five in the article \cite{13}. In Sections $3$ and $4$ we present first explicit equations of period six and seven orbits on the new $(u,v)$-plane. Moreover, we present all the central critical points and their real cycles of period six and seven cases. We see the curves of period six orbits on the $(u,v)$-plane in the Figure $\ref{kuva2}$ and on the $(x,y)$-plane in the Figure $\ref{kuva1}$. The curves of the period six orbits and curves of eigenvalues $|\lambda_{6}(x,y)|=1$ and $|\lambda_{6}(u,v)|=1$ are presented in the Figures $\ref{kuva3}$ and $\ref{kuva4}$.

\section{The $(u,v)$-model}
In the article \cite{5} we studied the dynamics of quadratic polynomials by iterating the function
\begin{equation}
 G(u,v)=\left(\frac{-u+v+uv}{u},\frac{u^{2}-u+v-u^{2}v-uv+uv^{2}+v^{2}}{u}\right), \label{1}
\end{equation}
and forming corresponding iterating system on the $(u,v)$-plane. 
We obtained a new model by using the change of variables
\begin{equation}
\left\{%
\begin{array}{l}
 u=x+y  \\
v=x+y^2+y-x^2   \label{3}
\end{array}%
\right.
\end{equation}
to the $(x,y)$-plane model (\cite{1}),
when
\begin{equation}
 x=\frac{u^{2}+u-v}{2u},\qquad y=\frac{u^{2}-u+v}{2u}. \label{4}
 \end{equation}
The new iteration system was defined as follows:
\begin{eqnarray*}
G(u,v)&=&(R(u,v),Q(u,v)) \\
&=&\left(\frac{-u+v+uv}{u},\frac{u^2-u+v-u^2v-uv+uv^2+v^2}{u}\right),
\end{eqnarray*}
when recursively
\begin{equation}
\left\{%
\begin{array}{l}
      (R_{0}(u,v),Q_{0}(u,v))=(u,v), \\ \\
      (R_{n+1}(u,v),Q_{n+1}(u,v))=G(R_{n}(u,v),Q_{n}(u,v)),\label{5}
\end{array}%
\right.
\end{equation}
where
\begin{equation}
\left\{%
\begin{array}{l}
      R_{n+1}(u,v)=Q_{n}(u,v)-1+Q_{n}(u,v)/R_{n}(u,v) \\ \\
      Q_{n+1}(u,v)=R_{n+1}(u,v)(1+Q_{n}(u,v)-R_{n}(u,v)),\label{6}
\end{array}%
\right.
\end{equation}
with $n\in \mathbb{N}\cup\{0\}$.

Now $(u,v)$ is fixed $G^{n}$, so $G^{n}(u,v)=(u,v)$, if and only if $(R_{n}(u,v),Q_{n}(u,v))$=$(u,v)$. The set of such points is the union of all orbits, whose period divides $n$, and the set of periodic points of period $n$ are the points with exact period dividing $n$. As mentioned in the introduction, the main motivation to develop the new $(u,v)$-model is that the equations of periodic orbits obtained through it are of significantly lower degree than before.

 When we solve the pair of equations ($\ref{6}$) to higher periods, the situation gets complicated very fast due to the growth of degree in the equations. The following result from the article \cite{5}, Theorem $3.1$, is very useful in sections $3-4$ below.
\begin{lause} \label{t1}
  Let $n \in \mathbb{N}\setminus{\{1\}}$. Any point of exact period $n+1$ satisfies
\begin{equation}
\left\{%
\begin{array}{ll}
      R_{0}R_{1}R_{2} \cdots R_{n}=&1\\
      Q_{0}Q_{1}Q_{2} \cdots Q_{n}=&1.\label{uviterkaava}
\end{array}%
\right.
\end{equation}
\end{lause}

 When we denote
$T_{n+1}(u,v)=1+Q_{n}(u,v)-R_{n}(u,v)$ the eigenvalue $\lambda_{n}(u,v)$ of the periodic point of the period $n$ by function $G^{n}(u,v)$ is
\begin{equation}
    \lambda_{n}(u,v)=u\Big(\frac{\partial R_{n}}{\partial u}\frac{\partial T_{n}}{\partial
  v}-\frac{\partial R_{n}}{\partial v}\frac{\partial T_{n}}{\partial  u}\Big),  \label{50}
\end{equation}
where $n=1,2,\ldots$.
On the $(u,v)$-plane fixed points and periodic points are classified in the following way:
\begin{maar} \label{luokitus}
Let us assume that $(u,v) \in \mathbb{C}^{2}$ is a periodic point of period $n$ of the function $G(u,v)$. In that case
$(u,v)$ is
\begin{eqnarray*}
1.&&  \textrm{attracting}, \ \textrm{if} \ 0<|\lambda_{n}(u,v)|<1, \\ \nonumber
2.&&  \textrm{super-attracting or critical}, \ \textrm{if} \ \lambda_{n}(u,v)=0, \\ \nonumber
3.&&  \textrm{repulsive}, \ \textrm{if} \ |\lambda_{n}(u,v)|>1, \\ \nonumber
4.&&  \textrm{indifferent (neutral)}, \ \textrm{if} \ |\lambda_{n}(u,v)|=1, \nonumber
\end{eqnarray*}
when $n \in \mathbb{N}$.
 \end{maar}

\section{Periodic orbits of period six}
By ($\ref{6}$) periodic orbits of the period six satisfy the pair of equations
\begin{equation}
\left\{%
\begin{array}{ll}
      R_{6}(u,v)=&u \\
      Q_{6}(u,v)=&v.
\end{array}%
\right.  \label{91}
\end{equation}
By Theorem $\ref{t1}$ we also know that the pair of equations
\begin{equation}
\left\{%
\begin{array}{ll}
      R_{0}R_{1}R_{2}R_{3}R_{4}R_{5}=&1\\
      Q_{0}Q_{1}Q_{2}Q_{3}Q_{4}Q_{5}=&1 \label{uviterkaava1}
\end{array}%
\right.
\end{equation}
satisfies the equation of the period six orbits.
Therefore, the equation
\begin{equation}
1=R_{0}R_{1}R_{2}R_{3}R_{4}R_{5} \label{a8}
\end{equation}
includes now orbits of the period six.
When we calculate formulas $R_{1}$, $R_{2}$, $R_{3}$, $R_{4}$ and $R_{5}$ by using ($\ref{5}$) and then substitute these formulas to the equation ($\ref{a8}$) and factorize the result, we obtain
\begin{equation}
\frac{(1+u)(uv+1+v)P_{6}L_{6}}{u^{15}}=0,\label{105a}
\end{equation}
where
\begin{eqnarray}
P_{6}&=& a_{13}{u}^{13}+ a_{12}{u}^{12}+ a_{11}{u}^{11}+ a_{10}{u}^{10}+a_{9}{u}^{9}+ a_{8}{u}^{8}+ a_{7}{u}^{7}  \nonumber  \\
&&+ a_{6}{u}^{6}+a_{5}{u}^{5}+ a_{4}{u}^{4}+a_{3}u^{3}+a_{2}{u}^{2}+a_{1}u+a_{0}, \label{122ab}
\end{eqnarray}

\begin{eqnarray*}
a_{13}&=& -4\,{v}^{4}+{v}^{3}-4\,{v}^{6}+{v}^{7}+6\,{v}^{5} \\
a_{12}&=&    {v}^{2}-7\,{v}^{8}-3\,{v}^{4}-4\,{v}^{3}+30\,{v}^{5}+30\,{v}^{7}-47\,{v}^{6}    \\
a_{11}&=&     -18\,{v}^{5}-15\,{v}^{3}+{v}^{2}+174\,{v}^{7}+46\,{v}^{4}-107\,{v}^{6}+21\,{v}^{9}-102\,{v}^{8}  \\
a_{10}&=&    v-35\,{v}^{10}+8\,{v}^{3}-7\,{v}^{2}+206\,{v}^{9}-402\,{v}^{8}+120\,{v}^{6}+251\,{v}^{7}-208\,{v}^{5}+66\,{v}^{4}  \\
a_{9}&=&    -270\,{v}^{10}+47\,{v}^{3}+546\,{v}^{6}+640\,{v}^{9}-310\,{v}^{7}-449\,{v}^{8}-160\,{v}^{5}-73\,{v}^{4}-v\\
&&-5\,{v}^{2}+35\,{v}^{11}   \\
a_{8}&=&    -726\,{v}^{10}-145\,{v}^{4}+10\,{v}^{2}+645\,{v}^{9}+278\,{v}^{6}+452\,{v}^{8}-21\,{v}^{12}+7\,{v}^{3}-964\,{v}^{7}\\
&&-v+234\,{v}^{11}+231\,{v}^{5} \\
a_{7}&=&    2\,{v}^{4}+{v}^{2}+278\,{v}^{5}+3\,v-1-33\,{v}^{3}-729\,{v}^{10}-417\,{v}^{6}-130\,{v}^{12}+582\,{v}^{11}\\
&&+7\,{v}^{13}-378\,{v}^{9}+1218\,{v}^{8}-401\,{v}^{7}  \\
a_{6}&=&    -314\,{v}^{12}-16\,{v}^{5}-4\,{v}^{3}-1116\,{v}^{9}+493\,{v}^{8}-4\,{v}^{2}+63\,{v}^{4}+118\,{v}^{10}+488\,{v}^{7}\\
&&-371\,{v}^{6}-{v}^{14}+42\,{v}^{13}+613\,{v}^{11}  \\
a_{5}&=&   -489\,{v}^{9}+8\,{v}^{4}-355\,{v}^{12}+720\,{v}^{10}+3\,{v}^{3}-77\,{v}^{5}+106\,{v}^{11}-6\,{v}^{14}+102\,{v}^{13}\\
&&-384\,{v}^{8}+20\,{v}^{6}+360\,{v}^{7} \\
a_{4}&=&   125\,{v}^{13}-296\,{v}^{11}-257\,{v}^{8}-15\,{v}^{14}+62\,{v}^{6}-13\,{v}^{7}+198\,{v}^{9}-{v}^{4}-148\,{v}^{12}\\
&&+365\,{v}^{10}-8\,{v}^{5} \\
a_{3}&=&     -191\,{v}^{11}+136\,{v}^{9}-20\,{v}^{14}+7\,{v}^{8}-33\,{v}^{7}+52\,{v}^{12}+4\,{v}^{6}+75\,{v}^{13}-56\,{v}^{10}  \\
a_{2}&=&    12\,{v}^{13}+62\,{v}^{12}+11\,{v}^{8}-54\,{v}^{10}-15\,{v}^{14}-2\,{v}^{11}-{v}^{7}-5\,{v}^{9}  \\
a_{1}&=&    16\,{v}^{11}+9\,{v}^{12}+3\,{v}^{10}-6\,{v}^{14}-8\,{v}^{13}-2\,{v}^{9}  \\
a_{0}&=&-{v}^{14}-3\,{v}^{12}-3\,{v}^{13}-{v}^{11},
\end{eqnarray*}
and
\begin{eqnarray}
L_{6}&=&b_{16}{u}^{16}+ b_{15}{u}^{15}+ b_{14}{u}^{14}+ b_{13}{u}^{13}+ b_{12}{u}^{12}+ b_{11}{u}^{11}+ b_{10}{u}^{10}  \nonumber \\
&&+b_{9}{u}^{9}+ b_{8}{u}^{8}+ b_{7}{u}^{7}+ b_{6}{u}^{6}+b_{5}{u}^{5}+ b_{4}{u}^{4}+b_{3}u^{3}+b_{2}{u}^{2}+b_{1}u+b_{0}, \label{122abc}
\end{eqnarray}
where
\begin{eqnarray*}
 b_{16}&=& {v}^{4}-4\,{v}^{5}+{v}^{8}+6\,{v}^{6}-4\,{v}^{7} \\
 b_{15}&=& 36\,{v}^{8}-8\,{v}^{9}+2\,{v}^{3}-6\,{v}^{5}-60\,{v}^{7}-6
\,{v}^{4}+42\,{v}^{6} \\
 b_{14}&=& {v}^{2}-31\,{v}^{4}+2\,{v}^{6}-148\,{v}^{9}-212\,{v}^{7}+28\,{v}^{10}+78\,{v}^{5}+282\,{v}^{8} \\
b_{13}&=& -440\,{v}^{6}+66\,{v}^{7}+694\,{v}^{8}-56\,
{v}^{11}+196\,{v}^{5}+364\,{v}^{10}+3\,{v}^{2}-24\,{v}^{3}-820\,{v}^{9
}\\&&+17\,{v}^{4} \\
b_{12}&=& 1472\,{v}^{7}-138\,{v}^{5}-201
\,{v}^{8}+1630\,{v}^{10}-748\,{v}^{6}-28\,{v}^{3}-588\,{v}^{11}-1648\,
{v}^{9}\\&&+70\,{v}^{12}+183\,{v}^{4}+v-5\,{v}^{2} \\
b_{11}&=& 2012\,{v}^{7}-2304\,{v}^{11}-3288\,{v}^{8}-15\,{v}^{2}+53\,{v}
^{3}+127\,{v}^{4}-753\,{v}^{5}+2954\,{v}^{10}\\&&+490\,{v}^{6}+136\,{v}^{9
}+644\,{v}^{12}-56\,{v}^{13} \\
b_{10}&=& 558\,{v}^{10}-
392\,{v}^{5}+5142\,{v}^{9}+5\,{v}^{2}-3\,v-4067\,{v}^{8}-238\,{v}^{4}-
476\,{v}^{13}+28\,{v}^{14}\\&&-958\,{v}^{7}+1997\,{v}^{6}-4008\,{v}^{11}+
84\,{v}^{3}+2330\,{v}^{12} \\
b_{9}&=& 19\,{v}^{2}-5640\,
{v}^{10}-23\,{v}^{3}-1806\,{v}^{11}+6284\,{v}^{9}-3696\,{v}^{7}+981\,{
v}^{8}+922\,{v}^{6}\\&&+228\,{v}^{14}+4042\,{v}^{12}+621\,{v}^{5}-8\,{v}^{
15}-1652\,{v}^{13}-276\,{v}^{4} \\
b_{8}&=& -2924\,{v}^{13}+604\,{v}^{5}+4104\,{v}^{11}+45\,{v}^{4}+{v}^{16}+3\,v-64\,{v}^{15}+
2720\,{v}^{12}\\&&-66\,{v}^{9}+4938\,{v}^{8}+780\,{v}^{14}-1-1684\,{v}^{7}
-7358\,{v}^{10}-60\,{v}^{3}-1050\,{v}^{6} \\
b_{7}&=& -1384\,{v}^{10}+117\,{v}^{4}-6\,{v}^{2}+8\,{v}^{16}-5\,{v}^{3}-42\,{v}^
{5}-2498\,{v}^{13}+2356\,{v}^{8}\\&&-220\,{v}^{15}-1562\,{v}^{12}-4783\,{v
}^{9}-927\,{v}^{6}+6352\,{v}^{11}+1428\,{v}^{14}+1198\,{v}^{7}\\
b_{6}&=& 28\,{v}^{16}+15\,{v}^{4}+3265\,{v}^{10}-232\,
{v}^{13}-152\,{v}^{5}-3822\,{v}^{12}+1024\,{v}^{7}-420\,{v}^{15}\\&&+4\,{v
}^{6}-2486\,{v}^{9}-919\,{v}^{8}+2262\,{v}^{11}+1434\,{v}^{14}+7\,{v}^
{3} \\
b_{5}&=& -1416\,{v}^{11}+442\,{v}^{9}+137\,{v}^{6}-
1980\,{v}^{12}-4\,{v}^{4}+1956\,{v}^{10}-19\,{v}^{5}+35\,{v}^{7}\\&&+624\,
{v}^{14}+1398\,{v}^{13}+56\,{v}^{16}-476\,{v}^{15}-825\,{v}^{8} \\
b_{4}&=& 1048\,{v}^{13}+490\,{v}^{9}-1128\,{v}^{11}-41
\,{v}^{8}+70\,{v}^{16}-88\,{v}^{10}+{v}^{5}+223\,{v}^{12}\\&&-162\,{v}^{14
}-88\,{v}^{7}+13\,{v}^{6}-308\,{v}^{15} \\
b_{3}&=& -5\,{v}^{7}+454\,{v}^{12}-300\,{v}^{14}-84\,{v}^{15}+21\,{v}^{9}+56\,{v}^{16
}-48\,{v}^{11}-220\,{v}^{10}\\&&+40\,{v}^{8}+146\,{v}^{13} \\
b_{2}&=& 54\,{v}^{12}-102\,{v}^{14}+{v}^{8}+78\,{v}^{11}-3\,{v}^{10}+
28\,{v}^{16}+20\,{v}^{15}-108\,{v}^{13}-12\,{v}^{9} \\
b_{1}&=& -26\,{v}^{13}-2\,{v}^{11}-22\,{v}^{12}+4\,{v}^{14}+20\,{v}^{15
}+2\,{v}^{10}+8\,{v}^{16} \\
b_{0}&=& u+{v}^{16}+4\,{v}^{15}+6\,{v}^{14}+4
\,{v}^{13}+{v}^{12}.
\end{eqnarray*}

We see that ($\ref{105a}$) includes also the periodic orbits of period two and period three.
 From these the factor $P_{6}=0$ satisfies the equation of period six orbits (Figure $\ref{kuva2}$), because it is satisfied also by period three orbits in the bifurcation points (see \cite{5}). This curve is of degree thirteen with respect to variable $u$ and of degree fourteen with respect to variable $v$. Respectively on the $(x,y)$-plane the value of degree of the period six equation is $26$ as for both factors. We obtain this equation of period six by using formulas ($\ref{3}$) to the earlier expression $P_{6}$. It is of the form
 \begin{eqnarray}
 &&c_{27}{x}^{27}+c_{26}{x}^{26}+ c_{25}{x}^{25}+ c_{24}{x}^{24}+ c_{23}{x}^{23}+ c_{22}{x}^{22}+ c_{21}{x}^{21}+ c_{20}{x}^{20}+ c_{19}{x}^{19} \nonumber \\
 &&+ c_{18}{x}^{18} +c_{17}{x}^{17}+ c_{16}{x}^{16}+ c_{15}{x}^{15}+ c_{14}{x}^{14}+ c_{13}{x}^{13}+ c_{12}{x}^{12}+ c_{11}{x}^{11}+ c_{10}{x}^{10}  \label{122} \\
&&+ c_{9}{x}^{9}+ c_{8}{x}^{8}+ c_{7}{x}^{7}+ c_{6}{x}^{6}+c_{5}{x}^{5}+ c_{4}{x}^{4}+c_{3}x^{3}+c_{2}{x}^{2}+c_{1}x+c_{0}=0,  \nonumber
\end{eqnarray}
 where
\begin{eqnarray*}
c_{27}&=&  1   \\
c_{26}&=&  -y-1   \\
c_{25}&=&  -6-13\,{y}^{2}-12\,y  \\
c_{24}&=&   25\,{y}^{2}+18\,y+7+13\,{y}^{3} \\
c_{23}&=&  14+78\,{y}^{4}+66\,y+144\,{y}^{3}+132\,{y}^{2}   \\
c_{22}&=&   -222\,{y}^{4}-276\,{y}^{3}-78\,{y}^{5}-92\,y-210\,{y}^{2}-20   \\
c_{21}&=&  -18-886\,{y}^{3}-1056\,{y}^{4}-454\,{y}^{2}-792\,{y}^{5}-286\,{y}^{6}-140\,y  \\
c_{20}&=&   286\,{y}^{7}+2008\,{y}^{4}+1848\,{y}^{5}+1078\,{y}^{6}+224\,y+
720\,{y}^{2}+1472\,{y}^{3}+33  \\
c_{19}&=&   2640\,{y}^{7}+
2120\,{y}^{3}+4010\,{y}^{4}+17+4620\,{y}^{6}+5230\,{y}^{5}+752\,{y}^{2}+715\,{y}^{8}\\
&&+162\,y  \\
c_{18}&=&   -40-7060\,{y}^{4}-3355\,{y}^{8}-7260\,{y}^{7}-10070\,{y}^{6}-1364\,{y}^{2}-715\,{y}^{9}\\
&&-329\,y-9900\,{y}^{5}-3692\,{y}^{3}  \\
c_{17}&=&   -12870\,{y}^{8}-11-5940\,{y}^{9}-17970\,{y}^{6}-1287\,{y}^{10}-2674\,{y}^{3}-12972\,{y}^{5}\\
&&-737\,{y}^{2}-136\,y-18090\,{y}^{7}-6918\,{y}^{4} \\
c_{16}&=&  7227\,{y}^{10}+31455\,{y}^{8}+34632\,{y}^{6}+1692\,{y}^{2}+355\,y+13207\,{y}^{4}+5481\,{y}^{3}\\
&&+24300\,{y}^{5}+38040\,{y}^{7}+18810\,{y}^{9}+1287\,{y}^{11}+38  \\
c_{15}&=&  2152\,{y}^{3}+40980\,{y}^{9}+76\,y+24552\,{y}^{10}+1716\,{y}^{12}+44256\,{y}^{7}+16232\,{y}^{5}\\
&&+9504\,{y}^{11}+30496\,{y}^{6}+2+6692\,{y}^{4}+49260\,{y}^{8}+512\,{y}^{2} \\
c_{14}&=&   -15528\,{y}^{4}-1556\,{y}^{2}-11220\,{y}^{12}-61248\,{y}^{6}-298\,y-34056\,{y}^{11}-66324\,{y}^{10}\\
&&-102696\,{y}^{8}-1716\,{y}^{13}-34276\,{y}^{5}-94200\,{y}^{9}-26-5624\,{y}^{3}-88432\,{y}^{7}  \\
c_{13}&=&  -12088\,{y}^{5}-64092\,{y}^{11}-52360\,{y}^{7}-33264\,{y}^{12}-1716\,{y}^{14}-79156\,{y}^{8}\\
&&-4284\,{y}^{4}-214\,{y}^{2}-11088\,{y}^{13}-1-1172\,{y}^{3}-89628\,{y}^{10}-95256\,{y}^{9}\\
&&-8\,y-27828\,{y}^{6}  \\
c_{12}&=&  1090\,{y}^{2}+98280\,{y}^{12}+183\,y+114292\,{y}^{7}+12964\,{y}^{4}+65768\,{y}^{6}+44352\,{y}^{13}\\
&&+31940\,{y}^{5}+166586\,{y}^{8}+4294\,{y}^{3}+17+12804\,{y}^{14}+200256\,{y}^{10}+159264\,{y}^{11}\\
&&+1716\,{y}^{15}+201712\,{y}^{9}   \\
c_{11}&=&  102620\,{y}^{9}+32472\,{y}^{14}+15760\,{y}^{6}+70812\,{y}^{13}+1+112308\,{y}^{12}+1706\,{y}^{4}\\
&&+9504\,{y}^{15}+65566\,{y}^{8}+136752\,{y}^{11}+131712\,{y}^{10}+10\,{y}^{2}+5884\,{y}^{5}+296\,{y}^{3}\\
&&+34952\,{y}^{7}+1287\,{y}^{16}  \\
c_{10}&=&   -21282\,{y}^{5}-93608\,{y}^{7}-7986\,{y}^{4}-231438\,{y}^{9}-48068\,{y}^{6}-289352\,{y}^{10}-12\\
&&-2430\,{y}^{3}-103\,y-188664\,{y}^{13}-1287\,{y}^{17}-10791\,{y}^{16}-266952\,{y}^{12}\\
&&-158232\,{y}^{8}-305256\,{y}^{11}-586\,{y}^{2}-104076\,{y}^{14}-41976\,{y}^{15} \\
c_{9}&=&  -97860\,{y}^{14}-31690\,{y}^{8}-5940\,{y}^{17}-59240\,{y}^{9}-145292\,{y}^{12}+2-129052\,{y}^{11}\\
&&+16\,{y}^{3}-22770\,{y}^{16}-1488\,{y}^{5}-95182\,{y}^{10}-14500\,{y}^{7}-5454\,{y}^{6}-133560\,{y}^{13}\\
&&-55380\,{y}^{15}+11\,{y}^{2}-179\,{y}^{4}-715\,{y}^{18} \\
c_{8}&=& 5+58\,y+6655\,{y}^{18}+28710\,{y}^{17}+3603\,{y}^{4}+280\,{y}^{2}+246000\,{y}^{14}+1099\,{y}^{3}\\
&&+332094\,{y}^{12}+52494\,{y}^{7}+97680\,{y}^{8}+24874\,{y}^{6}+300254\,{y}^{11}+312872\,{y}^{13}\\
&&+10221\,{y}^{5}+78645\,{y}^{16}+160930\,{y}^{9}+715\,{y}^{19}+157200\,{y}^{15}+234720\,{y}^{10} \\
c_{7}&=& 88032\,{y}^{15}-32\,{y}^{5}+61496\,{y}^{11}+636\,{y}^{6}+106208\,{y}^{14}+58470\,{y}^{16}-32\,{y}^{3}\\
&&+37536\,{y}^{10}+2640\,{y}^{19}+105224\,{y}^{13}+286\,{y}^{20}+19940\,{y}^{9}+11220\,{y}^{18}\\
&&-16\,{y}^{2}-10\,y+3152\,{y}^{7}+8966\,{y}^{8}+30090\,{y}^{17}+87252\,{y}^{12}-70\,{y}^{4}\\
c_{6}&=&  -114\,{y}^{2}-154692\,{y}^{16}-40526\,{y}^{8}-20\,y-13860\,{y}^{19}-8988\,{y}^{6}-3\\
&&-116988\,{y}^{10}-41530\,{y}^{18}-251872\,{y}^{14}-20244\,{y}^{7}-3606\,{y}^{5}-286\,{y}^{21}\\
&&-215600\,{y}^{15}-424\,{y}^{3}-252116\,{y}^{13}-2926\,{y}^{20}-169480\,{y}^{11}-1304\,{y}^{4}\\
&&-219832\,{y}^{12}-90540\,{y}^{17}-72606\,{y}^{9} \\
c_{5}&=& -919\,{y}^{8}+16\,{y}^{2}-3368\,{y}^{9}+52\,{y}^{4}-53992\,{y}^{15}-22910\,{y}^{18}-15676\,{y}^{11}\\
&&+38\,{y}^{3}+96\,{y}^{5}-49626\,{y}^{16}+1-3696\,{y}^{20}-49412\,{y}^{14}-26252\,{y}^{12}\\
&&-10830\,{y}^{19}-78\,{y}^{22}-38632\,{y}^{13}-8066\,{y}^{10}-792\,{y}^{21}-37596\,{y}^{17}+134\,{y}^{6} \\
c_{4}&=& 63536\,{y}^{18}+33\,{y}^{2}+4488\,{y}^{21}+78\,{y}^{23}+36390\,{y}^{10}+116\,{y}^{3}+20961\,{y}^{9}\\
&&+5254\,{y}^{7}+57522\,{y}^{11}+368\,{y}^{4}+34400\,{y}^{19}+14592\,{y}^{20}+972\,{y}^{5}\\
&&+127128\,{y}^{14}+121813\,{y}^{16}+95952\,{y}^{17}+10951\,{y}^{8}+870\,{y}^{22}+9\,y+1\\
&&+2352\,{y}^{6}+82716\,{y}^{12}+132996\,{y}^{15}+108004\,{y}^{13}  \\
c_{3}&=& -1-2\,y-2\,{y}^{2}+13\,{y}^{24}+144\,{y}^{23}+732\,{y}^{22}+2326\,{y}^{21}+5314\,{y}^{20}+9416\,{y}^{19}\\
&&+13456\,{y}^{18}+15890\,{y}^{17}-32\,{y}^{4}-8\,{y}^{3}-107\,{y}^{8}-96\,{y}^{7}-64\,{y}^{6}-46\,{y}^{5}\\
&&+602\,{y}^{10}+32\,{y}^{9}+3790\,{y}^{12}+1800\,{y}^{11}+15833\,{y}^{16}+13528\,{y}^{15}+10064\,{y}^{14}  \\
&& +6612\,{y}^{13} \\
c_{2}&=&  1-3\,y-8\,{y}^{2}-13\,{y}^{25}-157\,{y}^{24}-876\,{y}^{23}-3070\,{y}^{22}-7772\,{y}^{21}-15396\,{y}^{20}\\
&&-24972\,{y}^{19}-34116\,{y}^{18}-40177\,{y}^{17}-65\,{y}^{4}-26\,{y}^{3}-1852\,{y}^{8}-904\,{y}^{7}\\
&&-172\,{y}^{5}-6522\,{y}^{10}-3571\,{y}^{9}-17126\,{y}^{12}-11022\,{y}^{11}-41568\,{y}^{16}\\
&&-32044\,{y}^{14}-24470\,{y}^{13}-414\,{y}^{6}-38360\,{y}^{15} \\
c_{1}&=&  -1+2\,y+{y}^{2}-{y}^{26}-12\,{y}^{25}-66\,{y}^{24}-226\,{y}^{23}-554\,{y}^{22}-1052\,{y}^{21}+8\,{y}^{5}\\
&&-2050\,{y}^{19}-2201\,{y}^{18}-2032\,{y}^{17}+{y}^{4}+2\,{y}^{3}+26\,{y}^{8}+18\,{y}^{7}+16\,{y}^{6}\\
&&+31\,{y}^{10}+32\,{y}^{9}-149\,{y}^{12}-16\,{y}^{11}-1635\,{y}^{16}-1164\,{y}^{15}-730\,{y}^{14}-384\,{y}^{13} \\
&&-1614\,{y}^{20}\\
c_{0}&=&1-y+{y}^{2}+{y}^{27}+13\,{y}^{26}+78\,{y}^{25}+293\,{y}^{24}+792\,{y}^{23}+1672\,{y}^{22}+3\,{y}^{3}\\
&&+2892\,{y}^{21}+4219\,{y}^{20}+5313\,{y}^{19}+5892\,{y}^{18}+5843\,{y}^{17}+7\,{y}^{4}+927\,{y}^{11}\\
&&+155\,{y}^{8}+76\,{y}^{7}+35\,{y}^{6}+17\,{y}^{5}+536\,{y}^{10}+298\,{y}^{9}+1525\,{y}^{12}+2331\,{y}^{13}\\
&&+5258\,{y}^{16}+4346\,{y}^{15}+3310\,{y}^{14}.
\end{eqnarray*}
Curves of the period six ($\ref{122}$) are presented on the $(x,y)$-plane in Figure $\ref{kuva1}$.
The eigenvalue of the period six by ($\ref{50}$) on the $(u,v)$-plane is
\begin{equation}
    \lambda_{6}(u,v)=u\Big(\frac{\partial R_{6}}{\partial u}\frac{\partial T_{6}}{\partial
  v}-\frac{\partial R_{6}}{\partial v}\frac{\partial T_{6}}{\partial  u}\Big),  \label{500}
\end{equation}
where $T_{6}(u,v)=1+Q_{5}(u,v)-R_{5}(u,v)$. We obtain
\begin{equation}
\lambda_{6}(u,v)=\frac{a_{6}b_{6}c_{6}d_{6}e_{6}f_{6}}{u^{10}},\label{105}
\end{equation}
where
\begin{eqnarray*}
a_{6}&=&  {u}^{2}+u-v \\
b_{6}&=& {u}^{2}-u+v \\
c_{6}&=&  -2\,vu-v+u+{u}^{2} \\
d_{6}&=& -2\,{v}^{2}u-2\,{v}^{2}-v+2\,vu+u+2\,v{u}^{2}-{u}^{2} \\
e_{6}&=& -2\,v{u}^{4}+4\,v{u}^{2}+6\,{v}^{2}{u}^{3}+vu-4\,{v}^{2}u+4\,{v}^{4}u-8\,{v}^{3}{u}^{2}-2\,{v}^{3}u+2\,{v}^{4}\\
&&+2\,{u}^{4}{v}^{2}+2\,{u}^{2}{v}^{4}-4\,{u}^{3}{v}^{3}-{u}^{2}-{u}^{3}+2\,{v}^{3}\\
f_{6}&=& v{u}^{3}-4\,v{u}^{4}+8\,{v}^{2}{u}^{3}+62\,{v}^{5}{u}^{3}+8\,{v}^{8}u+2\,{v}^{4}u-6\,{v}^{3}{u}^{2}-30\,{v}^{6}{u}^{2}+4\,{v}^{7}u\\
&&-12\,{v}^{6}u+2\,{v}^{6}+20\,{u}^{2}{v}^{5}-22\,{u}^{3}{v}^{4}-6\,u{v}^{5}+2\,{u}^{8}{v}^{2}-2\,{u}^{7}{v}^{2}+12\,{u}^{6}{v}^{3}\\
&&+46\,{u}^{5}{v}^{3}-70\,{u}^{4}{v}^{4}+2\,{u}^{6}v-4\,{u}^{5}v+12\,{u}^{4}{v}^{2}+12\,{u}^{2}{v}^{4}+2\,{u}^{7}v-16\,{u}^{6}{v}^{2}\\
&&-8\,{u}^{5}{v}^{2}+16\,{u}^{4}{v}^{3}-16\,{u}^{3}{v}^{3}-8\,{u}^{5}{v}^{7}+2\,{u}^{4}{v}^{8}+8\,{u}^{3}{v}^{8}+2\,{u}^{8}{v}^{4}-20\,{u}^{2}{v}^{7}\\
&&-8\,{u}^{7}{v}^{5}+12\,{u}^{6}{v}^{6}-4\,{u}^{8}{v}^{3}+20\,{u}^{7}{v}^{4}-44\,{u}^{6}{v}^{5}+52\,{u}^{5}{v}^{6}-32\,{u}^{4}{v}^{7}\\
&&-12\,{u}^{7}{v}^{3}+34\,{u}^{6}{v}^{4}-60\,{u}^{5}{v}^{5}+68\,{u}^{4}{v}^{6}-44\,{u}^{3}{v}^{7}+12\,{u}^{2}{v}^{8}-22\,{u}^{5}{v}^{4}\\
&&+12\,{u}^{4}{v}^{5}+12\,{u}^{3}{v}^{6}-{u}^{4}+{u}^{5}+4\,{v}^{7}+2\,{v}^{8}.
\end{eqnarray*}
The bifurcation points are produced by the pair of equations
\begin{displaymath}
\left\{%
\begin{array}{ll}
      P_{6}(u,v)&=0\\
      |\lambda_{6}(u,v)|&=1,
\end{array}%
\right.
\end{displaymath}
for which the Singular program (\cite{10}) gives altogether $520$ solutions with multiplicity. In the Figure $\ref{kuva4}$ we see real bifurcation points by the intersection points of the period six and eigenvalue curves $|\lambda_{6}(u,v)|=1$.
The corresponding points on the $(x,y)$-plane are depicted in Figure $\ref{kuva3}$.
We obtain all critical points of the period six by solving
\begin{equation}
\left\{%
\begin{array}{ll}
      P_{6}(u,v)&=0\\
      |\lambda_{6}(u,v)|&=0, \label{ab}
\end{array}%
\right.
\end{equation}
which has $260$ solutions by Singular.
 Because $v=u^{2}+u$ is the critical curve on the $(u,v)$-plane (see \cite{5}), we obtain the central critical points from the pair of equations
\begin{displaymath}
\left\{%
\begin{array}{l}
      P_{6}(u,v)=0\\
      v=u^{2}+u,
\end{array}%
\right.
\end{displaymath}
which is equivalent with the equation
 \begin{eqnarray}
&& -{u}^{34}-13\,{u}^{33}-78\,{u}^{32}-293\,{u}^{31}-792\,{u}^{30}-1672\,{u}^{29}-2892\,{u}^{28}-4219\,{u}^{27} \nonumber \\
&&-5313\,{u}^{26}-5892\,{u}^{25}-5843\,{u}^{24}-5258\,{u}^{23}-4346\,{u}^{22}-3310\,{u}^{21}-2331\,{u}^{20} \nonumber \\
&&-1525\,{u}^{19}-927\,{u}^{18}-536\,{u}^{17}-298\,{u}^{16}-155\,{u}^{15}-76\,{u}^{14}-35\,{u}^{13}-17\,{u}^{12} \nonumber \\
&&-7\,{u}^{11}-3\,{u}^{10}-{u}^{9}+{u}^{8}-{u}^{7}=0. \nonumber
\end{eqnarray}
The central critical point is a critical point of the critical cycle, which is located on the critical curve. So by calculating central critical points first, we obtain other points of the cycle by iterating. A great advantage of this method is that the system of equations is of significantly lower degree than system ($\ref{ab}$).
Since $u=0$ is the point of singularity (see \cite{5}, Theorem $2.1$),
the real central critical points by a numerical approximation of six decimals are (see \cite{3} Table p.$106$ or \cite{4})
\begin{displaymath}
\begin{array}{l}
(u,v)\approx(-1.476015, 0.702605), \\ (u,v)\approx(-1.772893, 1.370256), \\  (u,v)\approx(-1.907280, 1.730437), \\  (u,v)\approx(-1.966773, 1.901424), \\  (u,v)\approx(-1.996376, 1.989142).
\end{array}
\end{displaymath}
Table $3.1$ includes these points with their cycles. These points has been ordered according to
the action of the map so the reader can easily see the orbits. The complex solutions are
\begin{displaymath}
\begin{array}{l}
 (u,v)\approx(0.389007+0.215851i, 0.493742+0.383785i),\\ (u,v)\approx(0.443326+0.372962i, 0.500762+0.703650i),\\ (u,v)\approx(0.396535+0.604182i, 0.188739+1.083340i),\\ (u,v)\approx(0.359893+0.684762i, 0.020516+1.177644i),\\ (u,v)\approx(-0.015570+1.020497i, -1.056743+.988718i),\\ (u,v)\approx(-0.113419+0.860570i, -0.841135+0.665360i),\\ (u,v)\approx(-0.163598+1.097781i, -1.341956+0.738591i),\\ (u,v)\approx(-0.217527+1.114454i, -1.412217+0.629607i),\\ (u,v)\approx(-0.596892+0.662981i, -0.680155-0.128475i),\\ (u,v)\approx(-1.284085+0.427269i, 0.182230-0.670030i),\\ (u,v)\approx(-1.138001+0.240332i, 0.099285-0.306664i),\\ (u,v)\approx(-1.138001-0.240332i, 0.099285+0.306664i) \\ (u,v)\approx(-1.284085-0.427269i, 0.182230+0.670030i),\\ (u,v)\approx(-0.596892-0.662981i, -0.680155+0.128475i),\\ (u,v)\approx(-0.217527-1.114455i, -1.412217-0.629607i),\\ (u,v)\approx(-0.163598-1.097781i, -1.341956-0.738591i),\\ (u,v)\approx(-0.113419-0.860569i, -0.841135-0.665360i),\\ (u,v)\approx(-0.015570-1.020497i, -1.056743-0.988718i),\\ (u,v)\approx(0.359893-0.684762i, 0.020516-1.177644i),\\ (u,v)\approx(0.396535-0.604182i, 0.188739-1.083340i),\\ (u,v)\approx(0.443326-0.372962i, 0.500762-0.703650i),\\ (u,v)\approx(0.389007-0.2158507i, 0.493742-0.383785i).
\end{array}
\end{displaymath}
There are altogether $162$ critical points.

\section{Periodic orbits of period seven}
By ($\ref{6}$) periodic points of the period seven satisfy the pair of equations
\begin{equation}
\left\{%
\begin{array}{ll}
      R_{7}(u,v)=&u \\
      Q_{7}(u,v)=&v.
\end{array}%
\right.  \label{91a}
\end{equation}
By Theorem $\ref{t1}$ we now know that the pair of equations
\begin{equation}
\left\{%
\begin{array}{ll}
      R_{0}R_{1}R_{2}R_{3}R_{4}R_{5}=&1\\
      Q_{0}Q_{1}Q_{2}Q_{3}Q_{4}Q_{5}=&1 \label{uviterkaava2}
\end{array}%
\right.
\end{equation}
must satisfies the equation of the period seven orbits and so the equation
\begin{equation}
1=R_{0}R_{1}R_{2}R_{3}R_{4}R_{5}R_{6}\label{aa8}
\end{equation}
includes the orbits of the period seven.
 After calculating $R_{1}$, $R_{2}$, $R_{3}$, $R_{4}$, $R_{5}$ and $R_{6}$ by using ($\ref{5}$) we substitute these formulas to the equation ($\ref{aa8}$) and factorize it. We obtain
\begin{equation}
\frac{P_{7}L_{7}}{u^{31}}=0,\label{105aa}
\end{equation}
where
\begin{eqnarray}
P_{7}&=& a_{31}{u}^{31}+ a_{30}{u}^{30}+ a_{29}{u}^{29}+ a_{28}{u}^{28}+a_{27}{u}^{27}+a_{26}{u}^{26}+ a_{25}{u}^{25}+ a_{24}{u}^{24}  \nonumber \\
 &&+ a_{23}{u}^{23}+ a_{22}{u}^{22}+ a_{21}{u}^{21}+ a_{20}{u}^{20}+ a_{19}{u}^{19}+a_{18}{u}^{18}+a_{17}{u}^{17}+ a_{16}{u}^{16}  \label{123a} \\
 &&+ a_{15}{u}^{15}+ a_{14}{u}^{14}+ a_{13}{u}^{13}+ a_{12}{u}^{12}+ a_{11}{u}^{11}+ a_{10}{u}^{10}+a_{9}{u}^{9}+ a_{8}{u}^{8}  \nonumber \\
&&+ a_{7}{u}^{7}+ a_{6}{u}^{6}+a_{5}{u}^{5}+ a_{4}{u}^{4}+a_{3}u^{3}+a_{2}{u}^{2}+a_{1}u+a_{0},  \nonumber
\end{eqnarray}
and
\begin{eqnarray}
L_{7}&=& b_{32}{u}^{32}+ b_{31}{u}^{31}+ b_{30}{u}^{30}+ b_{29}{u}^{29}+ b_{28}{u}^{28}+b_{27}{u}^{27}+b_{26}{u}^{26}+ b_{25}{u}^{25}  \nonumber \\
 &&+ b_{24}{u}^{24}+ b_{23}{u}^{23}+ b_{22}{u}^{22}+ b_{21}{u}^{21}+ b_{20}{u}^{20}+ b_{19}{u}^{19}+b_{18}{u}^{18}+b_{17}{u}^{17}  \label{123ab} \\
 &&+ b_{16}{u}^{16}+ b_{15}{u}^{15}+ b_{14}{u}^{14}+ b_{13}{u}^{13}+ b_{12}{u}^{12}+ b_{11}{u}^{11}+ b_{10}{u}^{10}+b_{9}{u}^{9}  \nonumber \\
&&+ b_{8}{u}^{8}+ b_{7}{u}^{7}+ b_{6}{u}^{6}+b_{5}{u}^{5}+ b_{4}{u}^{4}+b_{3}u^{3}+b_{2}{u}^{2}+b_{1}u+b_{0}.  \nonumber
\end{eqnarray}
The explicit forms of the coefficients $a_{0}$--$a_{31}$ and $b_{0}$--$b_{32}$ are very lengthy and they have been presented in the Appendix.

%
%
%
%

\begin{table}[h]
\begin{center}
Table $3.1$.
\end{center}
\begin{center}
 Real critical points of period six on the $(u,v)$-plane.
\end{center}
\begin{displaymath}
\begin{array}{|ll|}
       \hline
(u,v)\approx(-1.4760153,0.702605) &(u,v)\approx(-1.772893,1.370256)\\\hline (u,v)\approx(-0.773410,-2.458376)& (u,v)\approx(-0.402637,-1.668183)\\\hline (u,v)\approx(-0.279757,0.191624)& (u,v)\approx(1.474966,-0.391672)\\\hline (u,v)\approx(-1.493342,-2.197275)& (u,v)\approx(-1.657219,1.436210)\\\hline (u,v)\approx(-1.725894,-0.510981)&(u,v)\approx(-0.430429,-1.761929)\\\hline (u,v)\approx(-1.214913,-2.690928)&(u,v)\approx(1.331500,-0.441393)\\\hline

  (u,v)\approx(-1.907280,1.730437)&(u,v)\approx(-1.966773,1.901424)\\\hline (u,v)\approx(-0.176843,-0.820147)&(u,v)\approx(-0.065350,-0.318134)\\\hline (u,v)\approx(2.817570,1.005015)&(u,v)\approx(3.550062,2.652660)\\\hline (u,v)\approx(0.361711,-0.293910)&(u,v)\approx(2.399875,0.246222)\\\hline (u,v)\approx(-2.106465,-0.725422)&(u,v)\approx(-0.651180,0.751237)\\\hline (u,v)\approx(-1.381043,-3.288323)& (u,v)\approx(-1.402417,-3.369190)\\\hline

   (u,v)\approx-1.996376,1.989142)&\\\hline
    (u,v)\approx(-0.007235,-0.036068)& \\\hline (u,v)\approx(3.949450,3.835573)&\\\hline (u,v)\approx(3.806739,3.373240)&\\\hline (u,v)\approx(3.259363,1.846431)&\\\hline (u,v)\approx(1.412932,-0.583445)&\\\hline
\end{array}
\end{displaymath}
\end{table}

 From these factors $P_{7}=0$ satisfies the equation of period seven orbits, since, for example, the point $(u,v)=(-1.574889, 0.905387)$ is an approximation of a seven-periodic point under iteration of function $G(u,v)$, and it is also the central critical point of the period seven (see below) so it satisfies the equation $P_{7}(u,v)=0$. Instead, this point is not a solution of the equation $L_{7}(u,v)=0$ and hence $L_{7}(u,v)=0$ cannot be the equation of period seven orbits.

  The curve $P_{7}=0$ is of degree $31$ regarding variable $u$ and variable $v$, and the curve $L_{7}=0$ is of degree $32$ regarding both variables. Respectively on the $(x,y)$-plane the value of degree of the period seven equation is $63$ for both factors and it is of the form
\begin{eqnarray}
 &&c_{63}{x}^{63}+c_{62}{x}^{62}+ c_{61}{x}^{61}+ c_{60}{x}^{60}+ c_{59}{x}^{59}+ c_{58}{x}^{58}+ c_{57}{x}^{57}+ c_{56}{x}^{56}+ c_{55}{x}^{55}  \nonumber \\
 &&+c_{54}{x}^{54}+c_{53}{x}^{53}+ c_{52}{x}^{52}+ c_{51}{x}^{51}+ c_{50}{x}^{50}+ c_{49}{x}^{49}+ c_{48}{x}^{48}+ c_{47}{x}^{47}+ c_{46}{x}^{46}  \nonumber \\
 &&+c_{45}{x}^{45}+c_{44}{x}^{44}+ c_{43}{x}^{43}+ c_{42}{x}^{42}+ c_{41}{x}^{41}+ c_{40}{x}^{40}+ c_{39}{x}^{39}+ c_{38}{x}^{38}+ c_{37}{x}^{37}  \label{123bc} \\
 &&+c_{36}{x}^{36}+c_{35}{x}^{35}+ c_{34}{x}^{34}+ c_{33}{x}^{33}+ c_{32}{x}^{32}+ c_{31}{x}^{31}+ c_{30}{x}^{30}+ c_{29}{x}^{29}+ c_{28}{x}^{28}  \nonumber \\
 &&+c_{27}{x}^{27}+c_{26}{x}^{26}+ c_{25}{x}^{25}+ c_{24}{x}^{24}+ c_{23}{x}^{23}+ c_{22}{x}^{22}+ c_{21}{x}^{21}+ c_{20}{x}^{20}+ c_{19}{x}^{19}  \nonumber \\
 &&+c_{18}{x}^{18}+c_{17}{x}^{17}+ c_{16}{x}^{16}+ c_{15}{x}^{15}+ c_{14}{x}^{14}+ c_{13}{x}^{13}+ c_{12}{x}^{12}+ c_{11}{x}^{11}+ c_{10}{x}^{10}  \nonumber \\
&&+c_{9}{x}^{9}+ c_{8}{x}^{8}+ c_{7}{x}^{7}+ c_{6}{x}^{6}+c_{5}{x}^{5}+ c_{4}{x}^{4}+c_{3}x^{3}+c_{2}{x}^{2}+c_{1}x+c_{0}=0,  \nonumber
\end{eqnarray}
 where the coefficients $c_{0}$-$c_{63}$ have been presented in the Appendix.

The eigenvalue of the period seven by ($\ref{50}$) on the $(u,v)$-plane is
\begin{equation}
\lambda_{7}(u,v)=\frac{a_{7}b_{7}c_{7}d_{7}e_{7}f_{7}g_{7}}{u^{18}},\label{105a}
\end{equation}
where
\begin{eqnarray*}
a_{7}&=&  {u}^{2}+u-v \\
b_{7}&=& {u}^{2}-u+v \\
c_{7}&=& -2\,vu-v+u+{u}^{2} \\
d_{7}&=&   -2\,u{v}^{2}-2\,{v}^{2}-v+2
\,vu+u+2\,{u}^{2}v-{u}^{2} \\
e_{7}&=& 2\,{v}^{3}+2\,{u}^{2}{v}^{4}-4\,{u}^{3}{v}^{3}+2\,{v}^{4}+2\,{u}^{4}{v}^{2}+6\,{u}^{3}{v}^{2}-8\,{u}^{2}{v}^{3}-2\,u{v}^{3}-4\,u{v}^{2}\\
&&-2\,{u}^{4}v-{u}^{3}+4\,{u}^{2
}v-{u}^{2}+vu+4\,{v}^{4}u \\
f_{7}&=&   2\,{v}^{6}+12\,{u}^{2}{v}^{4}+2\,{u}^{7}v-16\,{u}^{6}{v}^{2}-8\,{u}^{5}{v}^{2}+16\,{u}^{4}{v}^{3}-16\,{u}^{3}{v}^{3}+{u}^{5}+2\,{u}^{6}v\\
&&-4\,{u}^{5}v+12\,{u}^{4}{v}^{2}+8\,{u}^{3}{v}^{2}-6\,{u}^{2}{v}^{3}-4\,{u}^{4}v+{u}^{3}v-{u}^{4}+8\,{u}^{3}{v}^{8}-8\,{u}^{5}{v}^{7}\\
&&+2\,{u}^{4}{v}^{8}+2\,{u}^{8}{v}^{4}-8\,{u}^{7}{v}^{5}+12\,{u}^{6}{v}^{6}+20\,{u}^{7}{v}^{4}-44\,{u}^{6}{v}^{5}+52\,{u}^{5}{v}^{6}-32\,{u}^{4}{v}^{7}\\
&&-4\,{u}^{8}{v}^{3}+20\,{u}^{2}{v}^{5}-22\,{u}^{3}{v}^{4}-70\,{u}^{4}{v}^{4}+2\,{u}^{8}{v}^{2}-2\,{u}^{7}{v}^{2}+12\,{u}^{6}{v}^{3}+46\,{u}^{5}{v}^{3}\\
&&-6\,u{v}^{5}-12\,{u}^{7}{v}^{3}+4\,{v}^{7}+2\,{v}^{8}+34\,{u}^{6}{v}^{4}-60\,{u}^{5}{v}^{5}-22\,{u}^{5}{v}^{4}+68\,{u}^{4}{v}^{6}\\
&&+12\,{u}^{4}{v}^{5}-44\,{u}^{3}{v}^{7}+12\,{u}^{3}{v}^{6}+62\,{u}^{3}{v}^{5}+4\,u{v}^{7}-12\,u{v}^{6}+12\,{v}^{8}{u}^{2}+8\,{v}^{8}u\\
&&+2\,{v}^{4}u-20\,{u}^{2}{v}^{7}-30\,{u}^{2}{v}^{6} \\
g_{7}&=&  2\,{v}^{12}+2\,{v}^{16}+8\,{v}^{15}+8\,{v}^{13}-{u}^{9}+6\,{u}^{8}v+{u}^{7}v+40\,u{v}^{15}+6\,{u}^{13}{v}^{2}-56\,{u}^{12}{v}^{3}\\
&&+254\,{u}^{11}{v}^{4}-784\,{u}^{10}{v}^{5}+1844\,{u}^{9}{v}^{6}-3368\,{u}^{8}{v}^{7}+4712\,{u}^{7}{v}^{8}-4972\,{u}^{6}{v}^{9}\\
&&+3912\,{u}^{5}{v}^{10}-2256\,{u}^{4}{v}^{11}+908\,{u}^{3}{v}^{12}-216\,{u}^{2}{v}^{13}+106\,{u}^{11}{v}^{3}-476\,{u}^{10}{v}^{4}\\
&&+1242\,{u}^{9}{v}^{5}-2100\,{u}^{8}{v}^{6}+2396\,{u}^{7}{v}^{7}-1838\,{u}^{6}{v}^{8}+884\,{u}^{5}{v}^{9}-176\,{u}^{4}{v}^{10}\\
&&-96\,{u}^{3}{v}^{11}+108\,{u}^{2}{v}^{12}-52\,u{v}^{13}+4524\,{u}^{6}{v}^{11}+2\,{u}^{14}{v}^{2}-48\,{u}^{13}{v}^{3}\\
&&+366\,{u}^{12}{v}^{4}-1506\,{u}^{11}{v}^{5}+3994\,{u}^{10}{v}^{6}-7392\,{u}^{9}{v}^{7}+9876\,{u}^{8}{v}^{8}-9566\,{u}^{7}{v}^{9}\\
&&+6530\,{u}^{6}{v}^{10}-2832\,{u}^{5}{v}^{11}+446\,{u}^{4}{v}^{12}-10\,{u}^{12}{v}^{2}+34\,{u}^{13}{v}^{4}-276\,{u}^{12}{v}^{5}\\
&&+980\,{u}^{11}{v}^{6}-1916\,{u}^{10}{v}^{7}+1962\,{u}^{9}{v}^{8}-132\,{u}^{8}{v}^{9}-2768\,{u}^{7}{v}^{10}-12\,{u}^{15}{v}^{4}\\
&&+156\,{u}^{14}{v}^{5}-880\,{u}^{13}{v}^{6}+2944\,{u}^{12}{v}^{7}-6576\,{u}^{11}{v}^{8}+10284\,{u}^{10}{v}^{9}+4\,{u}^{15}{v}^{3}\\
&&-62\,{u}^{14}{v}^{4}+392\,{u}^{13}{v}^{5}-1496\,{u}^{12}{v}^{6}+4024\,{u}^{11}{v}^{7}-8134\,{u}^{10}{v}^{8}+12568\,{u}^{9}{v}^{9}\\
&&-14716\,{u}^{8}{v}^{10}-{u}^{8}+12\,{v}^{14}-8\,{u}^{16}{v}^{5}+84\,{u}^{15}{v}^{6}-424\,{u}^{14}{v}^{7}+2\,{u}^{16}{v}^{4}\\
&&-12\,{u}^{15}{v}^{5}+4\,{u}^{14}{v}^{6}+132\,{u}^{13}{v}^{7}-402\,{u}^{12}{v}^{8}+12\,{u}^{16}{v}^{6}-16\,{u}^{9}{v}^{15}+2\,{u}^{8}{v}^{16}\\
&&+16\,{u}^{7}{v}^{16}-112\,{u}^{13}{v}^{11}+140\,{u}^{12}{v}^{12}-112\,{u}^{11}{v}^{13}+56\,{u}^{10}{v}^{14}+1288\,{u}^{11}{v}^{12}\\
&&-952\,{u}^{10}{v}^{13}+456\,{u}^{9}{v}^{14}-128\,{u}^{8}{v}^{15}-3304\,{u}^{9}{v}^{13}+1560\,{u}^{8}{v}^{14}-440\,{u}^{7}{v}^{15}\\
&&+56\,{u}^{6}{v}^{16}+2856\,{u}^{7}{v}^{14}-840\,{u}^{6}{v}^{15}+112\,{u}^{5}{v}^{16}-952\,{u}^{5}{v}^{15}+140\,{u}^{4}{v}^{16}\\
&&+112\,{u}^{3}{v}^{16}-600\,{u}^{3}{v}^{14}+40\,{u}^{2}{v}^{15}+2\,{u}^{16}{v}^{8}-16\,{u}^{15}{v}^{9}+56\,{u}^{14}{v}^{10}\\
&&+72\,{u}^{15}{v}^{8}-296\,{u}^{14}{v}^{9}+728\,{u}^{13}{v}^{10}-1176\,{u}^{12}{v}^{11}-1640\,{u}^{13}{v}^{9}+3260\,{u}^{12}{v}^{10}\\
&&-4608\,{u}^{11}{v}^{11}+4660\,{u}^{10}{v}^{12}+5908\,{u}^{11}{v}^{10}-8016\,{u}^{10}{v}^{11}+8084\,{u}^{9}{v}^{12}\\
&&-5848\,{u}^{8}{v}^{13}-3612\,{u}^{9}{v}^{11}+5440\,{u}^{8}{v}^{12}-4996\,{u}^{7}{v}^{13}+2868\,{u}^{6}{v}^{14}\\
&&-3124\,{u}^{7}{v}^{12}-464\,{u}^{6}{v}^{13}+1248\,{u}^{5}{v}^{14}-616\,{u}^{4}{v}^{15}+2796\,{u}^{5}{v}^{13}-324\,{u}^{4}{v}^{14}\\
&&-168\,{u}^{3}{v}^{15}+56\,{u}^{2}{v}^{16}-8\,{u}^{16}{v}^{7}-120\,{u}^{15}{v}^{7}+564\,{u}^{14}{v}^{8}+1388\,{u}^{13}{v}^{8}\\
&&-3296\,{u}^{12}{v}^{9}+272\,{u}^{11}{v}^{9}+1116\,{u}^{10}{v}^{10}-11280\,{u}^{9}{v}^{10}+8208\,{u}^{8}{v}^{11}\\
&&+12704\,{u}^{7}{v}^{11}-7644\,{u}^{6}{v}^{12}-3960\,{u}^{5}{v}^{12}+2096\,{u}^{4}{v}^{13}+292\,{u}^{3}{v}^{13}\\
&&-204\,{u}^{2}{v}^{14}+8\,u{v}^{14}-440\,{u}^{3}{v}^{10}+156\,{u}^{2}{v}^{11}-44\,u{v}^{12}+42\,{u}^{3}{v}^{9}-6\,{u}^{2}{v}^{10}\\
&&+80\,{u}^{3}{v}^{8}-24\,{u}^{2}{v}^{9}+4\,u{v}^{10}-1854\,{u}^{7}{v}^{6}+2048\,{u}^{6}{v}^{7}-1650\,{u}^{5}{v}^{8}+980\,{u}^{4}{v}^{9}\\
&&+70\,{u}^{5}{v}^{7}-82\,{u}^{4}{v}^{8}-30\,{u}^{11}{v}^{2}+168\,{u}^{10}{v}^{3}-552\,{u}^{9}{v}^{4}+1208\,{u}^{8}{v}^{5}-46\,{u}^{9}{v}^{3}\\
&&+90\,{u}^{8}{v}^{4}-84\,{u}^{7}{v}^{5}+8\,{u}^{6}{v}^{6}+234\,{u}^{7}{v}^{4}-304\,{u}^{6}{v}^{5}+274\,{u}^{5}{v}^{6}-176\,{u}^{4}{v}^{7}\\
&&+2\,{u}^{12}v+10\,{u}^{10}{v}^{2}+38\,{u}^{9}{v}^{2}-120\,{u}^{8}{v}^{3}-6\,{u}^{10}v-12\,{u}^{7}{v}^{2}+14\,{u}^{6}{v}^{3}\\
&&-4\,u{v}^{11}-10\,{u}^{7}{v}^{3}+30\,{u}^{6}{v}^{4}-38\,{u}^{5}{v}^{5}-8\,{u}^{5}{v}^{4}+26\,{u}^{4}{v}^{6}+2\,{u}^{4}{v}^{5}-10\,{u}^{3}{v}^{7}\\
&&+2\,{v}^{8}{u}^{2}+16\,{v}^{16}u.
\end{eqnarray*}
We obtain the central critical points from the pair of equations
\begin{displaymath}
\left\{%
\begin{array}{l}
      P_{7}(u,v)=0\\
      v=u^{2}+u,
\end{array}%
\right.
\end{displaymath}
which is equivalent with the equation
\begin{eqnarray}
&&-{u}^{78}-32\,{u}^{77}-496\,{u}^{76}-4976\,{u}^{75}-36440\,{u}^{74}-208336\,{u}^{73}-971272\,{u}^{72}  \nonumber \\
&&-3807704\,{u}^{71}-12843980\,{u}^{70}-37945904\,{u}^{69}-99582920\,{u}^{68}-234813592\,{u}^{67} \nonumber \\
&&-502196500\,{u}^{66}-981900168\,{u}^{65}-1766948340\,{u}^{64}-2943492972\,{u}^{63} \nonumber \\
&&-4562339774\,{u}^{62}-6609143792\,{u}^{61}-8984070856\,{u}^{60}-11500901864\,{u}^{59} \nonumber \\
&&-13910043524\,{u}^{58}-15941684776\,{u}^{57}-17357937708\,{u}^{56}-17999433372\,{u}^{55} \nonumber \\
&&-17813777994\,{u}^{54}-16859410792\,{u}^{53}-15286065700\,{u}^{52}-13299362332\,{u}^{51} \nonumber \\
&&-11120136162\,{u}^{50}-8948546308\,{u}^{49}-6939692682\,{u}^{48}-5193067630\,{u}^{47}\nonumber \\
&&-3754272037\,{u}^{46}-2625062128\,{u}^{45}-1777171560\,{u}^{44}-1166067016\,{u}^{43} \nonumber \\
&&-742179284\,{u}^{42}-458591432\,{u}^{41}-275276716\,{u}^{40}-160617860\,{u}^{39} \nonumber \\
&&-91143114\,{u}^{38}-50323496\,{u}^{37}-27049196\,{u}^{36}-14162220\,{u}^{35}-7228014\,{u}^{34}\nonumber\\
&&-3598964\,{u}^{33}-1749654\,{u}^{32}-831014\,{u}^{31}-385741\,{u}^{30}-175048\,{u}^{29}\nonumber\\
&&-77684\,{u}^{28}-33708\,{u}^{27}-14290\,{u}^{26}-5916\,{u}^{25}-2398\,{u}^{24}-950\,{u}^{23}-365\,{u}^{22}\nonumber \\
&&-132\,{u}^{21}-42\,{u}^{20}-14\,{u}^{19}-5\,{u}^{18}-2\,{u}^{17}-{u}^{16}-{u}^{15}=0.\nonumber
\end{eqnarray}
Since $u=0$ is the point of singularity,
the real central critical points using a nu\-me\-ri\-cal approximation of six decimals are
\begin{displaymath}
\begin{array}{l}
(u,v)\approx(-1.574889, .905387), \\ (u,v)\approx(-1.674066, 1.128431), \\  (u,v)\approx(-1.832315, 1.525064), \\  (u,v)\approx(-1.884804, 1.667681), \\  (u,v)\approx(-1.927148, 1.786751), \\ (u,v)\approx(-1.953706, 1.863261), \\  (u,v)\approx(-1.977180, 1.932060), \\ (u,v)\approx(-1.991814, 1.975510), \\  (u,v)\approx(-1.999096, 1.997288).
\end{array}
\end{displaymath}

Table $4.1$ includes these points with their cycles and they have again been ordered according to
the action of the map. The complex solutions are
\begin{displaymath}
\begin{array}{l}
 (u,v)\approx(0.376009+0.144749i,  0.496439 + 0.253603i),\\ (u,v)\approx(0.432376+0.226760i,  0.567905 + 0.422851i), \\ (u,v)\approx(0.456823+0.347759i, 0.544575 + 0.665487i),\\ (u,v)\approx(0.452774+0.396170i,  0.500828 + 0.754922i),\\ (u,v)\approx(0.386539+0.569325i,  0.211821 + 1.009457i),\\ (u,v)\approx(0.412916+0.614807i,  0.205428 + 1.122534i),\\ (u,v)\approx(0.376893+0.678569i, 0.058486 + 1.190065i),\\ (u,v)\approx(0.352483+0.698337i, -0.010948 + 1.190641i),\\ (u,v)\approx(0.121193+0.610612i, -0.236966 + 0.758615i),\\ (u,v)\approx(0.014895+0.848149i, -0.704239 + 0.873416i),\\ (u,v)\approx(-0.006984+1.003604i, -1.014156 + 0.989586i),\\ (u,v)\approx(-0.014233+1.032915i, -1.080944 + 1.003511i),\\ (u,v)\approx(-0.127500+0.987461i,  -1.086323 + 0.735658i),\\ (u,v)\approx(-0.157516+1.109007i, -1.362600 + 0.759634i),\\ (u,v)\approx(-0.174578+1.071428i, -1.292058 + 0.697332i),\\

 (u,v)\approx(-0.207284+1.117481i, -1.413081 + 0.654209i),\\ (u,v)\approx(-0.224916+1.116260i, -1.420366 + 0.614131i),\\ (u,v)\approx(-0.272102+.842365i, -0.907641 + 0.383946i),\\
 (u,v)\approx(-0.530828+.668289i, -0.695659 - 0.041204i),\\ (u,v)\approx(-0.623532+0.681064i,  -0.698588 - 0.168267i),\\ (u,v)\approx(-0.622436+0.424878i, -0.415531 - 0.104041i),\\

 (u,v)\approx(-1.028194+0.361377i, -0.101604 - 0.381754i),\\ (u,v)\approx(-1.292558+0.438199i, 0.186130 - 0.694596i),\\ (u,v)\approx(-1.262287+0.408104i, 0.164533 - 0.622185i),\\
 (u,v)\approx(-1.252736+0.342471i, 0.199325 - 0.515580i),\\ (u,v)\approx(-1.408446+0.1361720i,  0.556732 - 0.247410i),\\      (u,v)\approx(-1.769262+0.056920i, 1.357785 - 0.144491i),\\

 (u,v)\approx(-1.769262-0.056920i, 1.357785 + 0.144491i),\\ (u,v)\approx(-1.408446-0.136172i,  0.556732 + 0.247410i),\\ (u,v)\approx(-1.252736-0.342471i, 0.199325 + 0.515580i),\\

 (u,v)\approx(-1.262287-0.408104i, 0.164533 + 0.622185i),\\ (u,v)\approx(-1.292558-0.438199i, 0.186130 + 0.694596i),\\ (u,v)\approx(-1.028194-0.361377i, -0.101604 + 0.381754i),\\
  \end{array}
\end{displaymath}
 \begin{displaymath}
\begin{array}{l}
 (u,v)\approx(-0.622436-0.424878i, -0.415531 + 0.104041i),\\ (u,v)\approx(-0.623532-0.681064i, -0.698588 + 0.168267i),\\ (u,v)\approx(-0.530828-0.668289i, -0.695659 + 0.041204i),\\

 (u,v)\approx(-0.272102-0.842365i, -0.907641 - 0.383946i),\\ (u,v)\approx(-0.224916-1.116260i, -1.420366 - 0.614131i),\\ (u,v)\approx(-0.207284-1.117481i,  -1.413081 - 0.654209i),\\ (u,v)\approx(-0.174578-1.071428i, -1.292058 - 0.697332i),\\ (u,v)\approx(-0.157516-1.109007i, -1.362600 - 0.759634i),\\ (u,v)\approx(-0.127500-0.987461i, -1.086323 - 0.735658i),\\ (u,v)\approx(-0.014233-1.032915i, -1.080944 - 1.003511i),\\ (u,v)\approx(-0.006984-1.003604i, -1.014156 - 0.989586i),\\ (u,v)\approx(0.014895-0.848149i, -0.704239 - 0.873416i),\\
 (u,v)\approx(0.121193-0.610612i, -0.236966 - 0.758615i),\\ (u,v)\approx(0.352483-0.698337i, -0.010948 - 1.190641i),\\ (u,v)\approx(0.376893-0.678569i,  0.058486 - 1.190065i),\\ (u,v)\approx(0.412916-0.614807i,  0.205428 - 1.122534i),\\ (u,v)\approx(0.386539-0.569325i,  0.211821 - 1.009457i),\\ (u,v)\approx(0.452774-0.396170i, 0.500828 - 0.754922i),\\ (u,v)\approx(0.456823-0.347759i, 0.544575 - 0.665487i),\\ (u,v)\approx(0.432376-0.226760i,  0.567905 - 0.422851i),\\ (u,v)\approx(0.376009-0.144749i, 0.496439 - 0.253603i).
\end{array}
\end{displaymath}

 \begin{table}[h]
\begin{center}
Table $4.1$.
\end{center}
\begin{center}
 Real critical points of period seven on the $(u,v)$-plane.
\end{center}
\begin{displaymath}
\begin{array}{|ll|}\hline
  (u,v)\approx(-1.574889, .905387)&(u,v)\approx(-1.674066, 1.128431)\\\hline
                   (u,v)\approx(-0.669502, -2.330053)& (u,v)\approx(-0.545635, -2.074775)\\\hline

                  (u,v)\approx (0.150223, -0.099230)&(u,v)\approx (0.727722, -0.385067)\\\hline

            (u,v)\approx(-1.759780, -1.320799)&(u,v) \approx(-1.914207, 0.215902)\\\hline

                    (u,v)\approx(-1.570252, -2.259562)&(u,v) \approx(-0.896887, -2.807355)\\\hline

                  (u,v) \approx(-1.820581, -0.565635)&(u,v)\approx(-0.677246, 0.616611)\\\hline

                    (u,v)\approx(-1.254946, -2.829835)&(u,v)\approx (-1.293857, -2.967923)\\\hline

 (u,v)\approx(-1.832315, 1.525064)&(u,v)\approx(-1.884804, 1.667681)\\\hline
                   (u,v)\approx(-0.307251, -1.338811)&(u,v) \approx(-0.217123, -0.988447)\\\hline

                   (u,v)\approx(2.018568, -0.063705)&(u,v) \approx (2.564037, 0.586332)\\\hline

                    (u,v)\approx(-1.095264, 1.185375)&(u,v) \approx (-0.184993, 0.180869)\\\hline

                  (u,v) \approx(-0.896898, -2.942399)&(u,v)  \approx(-1.796837, -2.454231)\\\hline

                   (u,v)\approx(-0.661760, 0.691871)&(u,v) \approx (-2.088369, -0.715487)\\\hline

                    (u,v)\approx(-1.353630, -3.185946)& (u,v)   \approx(-1.372882, -3.257685)\\\hline

  (u,v)\approx(-1.927148, 1.786751)&(u,v) \approx(-1.953706, 1.863261)\\\hline
                  (u,v)\approx (-0.140397, -0.661818)&(u,v)  \approx(-0.090445, -0.435671)\\\hline

                    (u,v) \approx(3.052081, 1.460663)&(u,v)  \approx (3.381296, 2.213985)\\\hline

                (u,v) \approx (0.939242, -0.555484)&(u,v)     \approx(1.868759, -0.312663)\\\hline

                   (u,v) \approx(-2.146902, 1.062130)& (u,v)  \approx(-1.479974, 1.748475)\\\hline

                 (u,v) \approx (-0.432597, -1.820814)& (u,v)  \approx (-0.432948, -1.830698)\\\hline

                 (u,v)  \approx (1.388217, -0.538930)&(u,v) \approx (1.397750, -0.555956)\\\hline

  (u,v) \approx(-1.977180, 1.932060)&(u,v)\approx(-1.991814, 1.975510)\\\hline
                (u,v) \approx (-0.045120, -0.221505)& (u,v) \approx(-0.016305, -0.080990)\\\hline

                  (u,v) \approx ( 3.687734, 3.037273)&(u,v)  \approx(3.886333, 3.634943)\\\hline

               (u,v) \approx (2.860888, 0.999991)& (u,v) \approx (3.570257, 2.672728)\\\hline

                 (u,v)\approx( 0.349530, -0.300909)&(u,v)\approx (2.421337, 0.248117)\\\hline

                (u,v) \approx (-2.161806, -0.755684)&(u,v) \approx(-0.649412, 0.761904)\\\hline

                  (u,v) \approx (-1.406122, -3.383302)&(u,v) \approx (-1.411317, -3.403131)\\\hline

 (u,v)\approx(-1.999096, 1.997288)&\\\hline
                 (u,v)\approx(-0.001808, -0.009033)&\\\hline

                   (u,v) \approx (3.987352, 3.958544)&\\\hline

                    (u,v)\approx (3.951319, 3.837491)&\\\hline

                  (u,v) \approx(3.808684, 3.375150)&\\\hline

                   (u,v) \approx (3.261322, 1.847427)&\\\hline

                   (u,v)\approx (1.413893, -0.585203)&\\\hline

\end{array}
\end{displaymath}
\end{table}
There are altogether $441$ critical points.

\vspace*{1.1cm}
\begin{center}
''The world around us is very complicated. The tools at our
disposal to describe it are very weak.''
\end{center}
\vspace*{0.1cm}
\begin{center}
Benoit Mandelbrot
\end{center}

\clearpage
\section{Appendix}
Figures of periodic orbit curves and curves of eigenvalues on the $(u,v)$- and the $(x,y)$-plane.

\clearpage
\clearpage

\clearpage
The coefficients of equations ($\ref{123a}$), ($\ref{123ab}$) and ($\ref{123bc}$) in period seven case.

\begin{eqnarray*}
{a}_{31}&=& -8\,{v}^{15}-56\,{v}^{11}-56\,{v}^{13}+{v}^{8}-8\,{v}^{9}+70\,{v}^{12}+28\,{v}^{14}+{v}^{16}+28\,{v}^{10}\\
 {a}_{30}&=& 1004\,{v}^{14}-15\,{v}^{17}-1302\,{v}^{13}+128\,{v}^{16}-28\,{v}^{8}-448\,{v}^{11}-476\,{v}^{15}+4\,{v}^{7}\\&&
 +36\,{v}^{10}+1036\,{v}^{12}+61\,{v}^{9} \\
{a}_{29} &=& -34\,{v}^{8}+3864\,{v}^{16}-1637\,{v}^{10}+105\,{v}^{18}+11742\,{v}^{14}-8668\,{
v}^{15}+6\,{v}^{6}\\&&-968\,{v}^{17}+3170\,{v}^{12}+608\,{v}^{9}+1140\,{v}
^{11}-9296\,{v}^{13}-32\,{v}^{7} \\
 {a}_{28}&=& -17670\,{v}^{12}+54376\,{v}^{14}+47904\,{v}^{16}-455\,{v}^{19}-202\,{v}^{7}+18087
\,{v}^{11}\\&&-5778\,{v}^{10}-6\,{v}^{6}-19888\,{v}^{17}-68446\,{v}^{15}+
998\,{v}^{8}+4592\,{v}^{18}+4\,{v}^{5}\\&&-766\,{v}^{9}-12750\,{v}^{13}\\
  {a}_{27}&=& {v}^{4}+30066\,{v}^{14}+72660\,{v}^{18}-
12988\,{v}^{9}-123365\,{v}^{12}+14\,{v}^{5}+15492\,{v}^{10}\\&&-234920\,{v
}^{15}-15288\,{v}^{19}+131486\,{v}^{13}-189920\,{v}^{17}+2405\,{v}^{8}
+1365\,{v}^{20}\\&&+289688\,{v}^{16}-208\,{v}^{6}+466\,{v}^{7}+33046\,{v}^
{11} \\
{a}_{26}&=& 101142\,{v}^{10}-232\,{v}^{6}-33638\,{v}^{15}-134550\,{v}^{11}-199836\,{v}^{19}-646494\,{v}^{14}\\&&+803096\,{v}^{
16}-946624\,{v}^{17}+37856\,{v}^{20}-131968\,{v}^{12}-15309\,{v}^{9}-
79\,{v}^{5}\\&&+573388\,{v}^{18}+601165\,{v}^{13}-3003\,{v}^{21}-8502\,{v}
^{8}+3578\,{v}^{7}+10\,{v}^{4} \\
  {a}_{25}&=& 956\,{v}^{7}+743982\,{v}^{12}-549110\,{v}^{11}+2\,{v}^{3}-1365108\,{v}^{19}-72072\,
{v}^{21}\\&&+2338382\,{v}^{15}-7540\,{v}^{16}+1857\,{v}^{6}+5005\,{v}^{22}-318\,{v}^{5}+428064\,{v}^{20}\\&&+2475844\,{v}^{18}+80104\,{v}^{9}-
2262372\,{v}^{17}+61121\,{v}^{10}-{v}^{4}-33331\,{v}^{8}\\&&+413592\,{v}^{13}-2259057\,{v}^{14} \\
 {a}_{24}&=& 5103\,{v}^{8}+5370092\,{v}^{18}-6554460\,{v}^{16}+271\,{v}^{5}-728728\,{v}^{21}-35140\,{v}^{
17}\\&&+2260378\,{v}^{12}+3874\,{v}^{6}+2617416\,{v}^{20}-21145\,{v}^{7}+
107536\,{v}^{22}-169477\,{v}^{11}\\&&+205721\,{v}^{9}-488420\,{v}^{10}+6\,
{v}^{3}-5296032\,{v}^{19}-1152656\,{v}^{14}-2964152\,{v}^{13}\\&&+6846367
\,{v}^{15}-6435\,{v}^{23}-119\,{v}^{4} \\
  {a}_{23}&=& -937353\,{v}^{10}-4091208\,{v}^{21}+996996\,{v}^{22}-10845324\,{v}^{19}
+829388\,{v}^{18}\\&&-77263\,{v}^{9}+9030784\,{v}^{14}+149197\,{v}^{8}-
4989\,{v}^{6}+6435\,{v}^{24}-100\,{v}^{4}\\&&+14631112\,{v}^{17}-18\,{v}^{
3}+2011\,{v}^{5}-7423354\,{v}^{13}-27482\,{v}^{7}+9380882\,{v}^{20}\\&&+
387297\,{v}^{12}+2143092\,{v}^{11}+3157728\,{v}^{15}+{v}^{2}-17180848
\,{v}^{16}-126984\,{v}^{23} \\
  {a}_{22}&=& 43575\,{v}^{7}-3816896\,{v}^{19}-8521504\,{v}^{16}-13841806\,{v}^{21}+18666340\,{v}^{
20}\\&&-7179724\,{v}^{12}+3365459\,{v}^{11}-17755\,{v}^{6}-{v}^{2}+473\,{v
}^{5}+441131\,{v}^{10}\\&&-21712088\,{v}^{15}+136502\,{v}^{8}+499\,{v}^{4}
-49\,{v}^{3}+20064244\,{v}^{14}-26321700\,{v}^{18}\\&&-5005\,{v}^{25}-
1101100\,{v}^{23}+5243260\,{v}^{22}-739474\,{v}^{9}+118976\,{v}^{24}\\&&-
1028735\,{v}^{13}+36205378\,{v}^{17} \\
  {a}_{21}&=& 3549145\,{v}^{14}+69\,{v}^{3}+10559090\,{v}^{20}-45427268\,{v}^{15}+38113392
\,{v}^{19}\\&&-64411298\,{v}^{18}+21088408\,{v}^{17}-9909897\,{v}^{12}-
5445\,{v}^{5}-1577987\,{v}^{11}\\&&-5512628\,{v}^{23}+980408\,{v}^{24}+
19018372\,{v}^{13}-10\,{v}^{2}+2764411\,{v}^{10}\\&&-239348\,{v}^{8}-
27295632\,{v}^{21}+17031938\,{v}^{22}+41783680\,{v}^{16}+102773\,{v}^{
7}\\&&+158\,{v}^{6}-527838\,{v}^{9}+3003\,{v}^{26}-88088\,{v}^{25}+592\,{v}^{4} \\
 {a}_{20}&=& -v-40627714\,{v}^{14}+189\,{v}^{3}-
43601822\,{v}^{20}-11838681\,{v}^{15}+96634514\,{v}^{19}\\&&-44947970\,{v}
^{18}-64308896\,{v}^{17}+3907089\,{v}^{12}-3903\,{v}^{5}-8136439\,{v}^
{11}\\&&-17423322\,{v}^{23}+4734544\,{v}^{24}+24483979\,{v}^{13}+2\,{v}^{2
}+1693084\,{v}^{10}\\&&-435680\,{v}^{8}-21002066\,{v}^{21}+33714296\,{v}^{
22}+86855316\,{v}^{16}-10502\,{v}^{7}\\&&+34857\,{v}^{6}+929402\,{v}^{9}-
1365\,{v}^{27}+50960\,{v}^{26}-698880\,{v}^{25}-991\,{v}^{4} \\
 {a}_{19}&=& -v-51307375\,{v}^{14}-80\,{v}^{3}-121452630\,{v}^{20}+
70570902\,{v}^{15}+80029422\,{v}^{19}\\&&+77472658\,{v}^{18}-140393194\,{v
}^{17}+19355098\,{v}^{12}+6949\,{v}^{5}-4647032\,{v}^{11}\\&&-34898416\,{v
}^{23}+14714308\,{v}^{24}-6664985\,{v}^{13}+34\,{v}^{2}-2709254\,{v}^{
10}\\&&+49943\,{v}^{8}+37302526\,{v}^{21}+32060458\,{v}^{22}+32130528\,{v}
^{16}-153547\,{v}^{7}\\&&+18305\,{v}^{6}+1436462\,{v}^{9}+455\,{v}^{28}-
22568\,{v}^{27}+393484\,{v}^{26}-3290928\,{v}^{25}\\&&-1518\,{v}^{4}\\
  {a}_{18}&=& 3\,v+6538121\,{v}^{14}-324\,{v}^{3}-117627052\,{v}^{20}+91458619\,{v}^{15}-68162532\,{v}^{19}\\&&+190984568\,{
v}^{18}-69036680\,{v}^{17}+10997458\,{v}^{12}+7652\,{v}^{5}+6122512\,{
v}^{11}\\&&-38641518\,{v}^{23}+29968400\,{v}^{24}-37782254\,{v}^{13}+10\,{
v}^{2}-3810848\,{v}^{10}\\&&+503980\,{v}^{8}+126227076\,{v}^{21}-19696926
\,{v}^{22}-99352564\,{v}^{16}-68377\,{v}^{7}\\&&-30722\,{v}^{6}-116662\,{v
}^{9}-105\,{v}^{29}+7392\,{v}^{28}-171108\,{v}^{27}+1822964\,{v}^{26}\\&&-
10141796\,{v}^{25}+699\,{v}^{4} \\
 {u}_{17}&=& 3\,v+60951286\,{v}^{14}-36\,{v}^{3}+32832036\,{v}^{20}+1924473\,{v}^{15}-216375728
\,{v}^{19}\\&&+118469442\,{v}^{18}+111153274\,{v}^{17}-10878884\,{v}^{12}-
3101\,{v}^{5}+8294076\,{v}^{11}\\&&-1532018\,{v}^{23}+37170392\,{v}^{24}-
22338406\,{v}^{13}-39\,{v}^{2}+91121\,{v}^{10}\\&&+210743\,{v}^{8}+
141821436\,{v}^{21}-105994786\,{v}^{22}-138387138\,{v}^{16}+95916\,{v}
^{7}\\&&-28019\,{v}^{6}-1287560\,{v}^{9}+15\,{v}^{30}-1688\,{v}^{29}+55440
\,{v}^{28}-785740\,{v}^{27}\\&&+5608288\,{v}^{26}-21066500\,{v}^{25}+1662\,{v}^{4} \\
 {a}_{16}&=& -3\,v+38648946\,{v}^{14}+199\,{v}^{3}+200246592\,{v}^{20}-81251972\,{v}^{15}\\&&-164024172\,{v}^{19}-93395916
\,{v}^{18}+176716376\,{v}^{17}-14970644\,{v}^{12}-5715\,{v}^{5}\\&&+355049
\,{v}^{11}+68706502\,{v}^{23}+16524272\,{v}^{24}+15158994\,{v}^{13}-19
\,{v}^{2}\\&&+2625390\,{v}^{10}-224065\,{v}^{8}+15740502\,{v}^{21}-
139138104\,{v}^{22}-22457228\,{v}^{16}\\&&+79396\,{v}^{7}+8563\,{v}^{6}-
536892\,{v}^{9}+240\,{v}^{30}-12600\,{v}^{29}+253896\,{v}^{28}\\&&-2426284
\,{v}^{27}+11904060\,{v}^{26}-28525280\,{v}^{25}-{v}^{31}+111\,{v}^{4} \\
 {a}_{15}&=& -1-3\,v-16095562\,{v}^{14}+65\,{v}^{3}+183972136\,{v}^{20}-56575654\,{v}^{15}\\&&+48245628\,{v}^{19}-188584370\,{
v}^{18}+51804868\,{v}^{17}-1648193\,{v}^{12}-391\,{v}^{5}\\&&-4334750\,{v}
^{11}+109427632\,{v}^{23}-30628580\,{v}^{24}+22534164\,{v}^{13}+12\,{v
}^{2}\\&&+1125784\,{v}^{10}-179075\,{v}^{8}-145772816\,{v}^{21}-54807110\,
{v}^{22}+88750297\,{v}^{16}\\&&-15593\,{v}^{7}+14455\,{v}^{6}+403601\,{v}^
{9}+1792\,{v}^{30}-57840\,{v}^{29}+790380\,{v}^{28}\\&&-5269460\,{v}^{27}+
17293612\,{v}^{26}-20394332\,{v}^{25}-16\,{v}^{31}-616\,{v}^{4} \\
 {a}_{14}&=& -28336650\,{v}^{14}-26\,{v}^{3}+5316014\,{v}^{20}+11764102\,{v}^{15}+165450096\,{v}^{19}\\&&-78657862\,{v}^{18}-
77692282\,{v}^{17}+5839218\,{v}^{12}+1335\,{v}^{5}-1940602\,{v}^{11}\\&&+
67829118\,{v}^{23}-67007248\,{v}^{24}+3908723\,{v}^{13}+6\,{v}^{2}-
567876\,{v}^{10}+17740\,{v}^{8}\\&&-166440628\,{v}^{21}+76249250\,{v}^{22}
+69683704\,{v}^{16}-28161\,{v}^{7}+1266\,{v}^{6}\\&&+325840\,{v}^{9}+8280
\,{v}^{30}-182168\,{v}^{29}+1758120\,{v}^{28}-8116388\,{v}^{27}\\&&+
15696324\,{v}^{26}+5119990\,{v}^{25}-120\,{v}^{31}-157\,{v}^{4} \\
  {a}_{13}&=& -6563925\,{v}^{14}-7\,{v}^{3}-115795986\,{v}
^{20}+29723798\,{v}^{15}+90078860\,{v}^{19}\\&&+51622457\,{v}^{18}-
71808864\,{v}^{17}+2755485\,{v}^{12}+281\,{v}^{5}+623344\,{v}^{11}\\&&-
19727526\,{v}^{23}-55712752\,{v}^{24}-6430862\,{v}^{13}-482126\,{v}^{
10}+43318\,{v}^{8}\\&&-43482370\,{v}^{21}+119499168\,{v}^{22}-3184426\,{v}
^{16}-3138\,{v}^{7}-2137\,{v}^{6}\\&&-7274\,{v}^{9}+26460\,{v}^{30}-415352
\,{v}^{29}+2847880\,{v}^{28}-8548176\,{v}^{27}\\&&
+4946900\,{v}^{26}+
29998092\,{v}^{25}-560\,{v}^{31}+36\,{v}^{4} \\
  {a}_{12}&=&5769064\,{v}^{14}-80900268\,{v}^{20}+8496737\,{v}^{15}-21813289\,{v}^{
19}+61424066\,{v}^{18}\\&&-5854568\,{v}^{17}-523051\,{v}^{12}-39\,{v}^{5}+
583643\,{v}^{11}-65522026\,{v}^{23}\\&&-9179068\,{v}^{24}-3234021\,{v}^{13
}-15984\,{v}^{10}+5754\,{v}^{8}+60497808\,{v}^{21}\\&&+54055608\,{v}^{22}-
25868010\,{v}^{16}+2572\,{v}^{7}-385\,{v}^{6}-53383\,{v}^{9}+61880\,{v
}^{30}\\&&-703612\,{v}^{29}+3329872\,{v}^{28}-5220012\,{v}^{27}-8158108\,{
v}^{26}+32823380\,{v}^{25}\\&&-1820\,{v}^{31}+4\,{v}^{4} \\
  {a}_{11}&=& 3148867\,{v}^{14}-295688\,{v}^{20}-4164772\,{v}^{15}-43037982
\,{v}^{19}+11396578\,{v}^{18}\\&&+18466340\,{v}^{17}-581252\,{v}^{12}-{v}^
{5}+40916\,{v}^{11}-42476176\,{v}^{23}-4368\,{v}^{31}\\&&+24868670\,{v}^{24}+314411\,{v}^
{13}+53257\,{v}^{10}-2333\,{v}^{8}+57327534\,{v}^{21}\\&&-18991186\,{v}^{
22}-8785970\,{v}^{16}+401\,{v}^{7}+32\,{v}^{6}-7898\,{v}^{9}+109200\,{
v}^{30}\\&&-891072\,{v}^{29}+2656290\,{v}^{28}+3164\,{v}^{27}-13578104\,{v
}^{26}+14427864\,{v}^{25} \\
  {a}_{10}&=& -103599\,{v}^{14}+24091192\,{v}^{20}-2551683\,{v}^{15}-11998580\,{v}^{
19}-10567488\,{v}^{18}\\&&+7386354\,{v}^{17}-53468\,{v}^{12}-43445\,{v}^{
11}-1617460\,{v}^{23}+23553308\,{v}^{24}\\&&+478825\,{v}^{13}+8253\,{v}^{
10}-310\,{v}^{8}+9996760\,{v}^{21}-31477458\,{v}^{22}\\&&+2350018\,{v}^{16
}-18\,{v}^{7}+1581\,{v}^{9}+147576\,{v}^{30}-830258\,{v}^{29}+1145144
\,{v}^{28}\\&&+3453120\,{v}^{27}-9094932\,{v}^{26}-4268952\,{v}^{25}-8008\,{v}^{31} \\
  {a}_{9}&=& 6328152\,{v}^{24}-8991576\,{v}^{25
}+3439688\,{v}^{27}+4601992\,{v}^{19}-11440\,{v}^{31}\\&&+29332\,{v}^{12}-
959216\,{v}^{17}+6\,{v}^{8}-1678614\,{v}^{26}-328098\,{v}^{14}-6647\,{
v}^{11}\\&&+153010\,{v}^{30}-780\,{v}^{10}-536536\,{v}^{29}-5091762\,{v}^{
18}+9113936\,{v}^{20}\\&&-27735\,{v}^{15}+1725350\,{v}^{16}-10198304\,{v}^
{21}+48776\,{v}^{13}-211640\,{v}^{28}\\&&-10035264\,{v}^{22}+12644454\,{v}
^{23}+172\,{v}^{9} \\
  {a}_{8}&=& -16626\,{v}^{13}-12870\,{v}^{31}+2766896\,{v}^{22}-975308\,{v}^{17}+1524362\,{v}^{27}-{v}^{9}\\&&+
2885166\,{v}^{19}+68316\,{v}^{16}-716976\,{v}^{28}-1298628\,{v}^{20}+
4148\,{v}^{12}\\&&+188122\,{v}^{15}-5317002\,{v}^{21}+1960588\,{v}^{26}-
192764\,{v}^{29}-33468\,{v}^{14}\\&&-3074416\,{v}^{24}+6162960\,{v}^{23}+
204116\,{v}^{18}+259\,{v}^{11}+120120\,{v}^{30}\\&&-4245668\,{v}^{25}-66\,{v}^{10} \\
{a}_{7}&=& -43\,{v}^{12}+31556\,{v}^{27}+27632\,{v}^{21}-55936\,{v}^{23}+21840\,{v}^{29}+461494\,{v}^{18}\\&&-55754\,{v}^{
17}+61660\,{v}^{19}-90884\,{v}^{16}-11440\,{v}^{31}-1340044\,{v}^{20}-
39500\,{v}^{25}\\&&+1612380\,{v}^{26}+17632\,{v}^{15}+2403422\,{v}^{22}-
2002\,{v}^{13}+8006\,{v}^{14}\\&&-2610751\,{v}^{24}-516516\,{v}^{28}+68640
\,{v}^{30}+16\,{v}^{11} \\
 {a}_{6}&=& 74256\,{v}^{29}+738\,{v}^{14}+217698\,{v}^{22}-168168\,{v}^{28}-328580\,{v}^{27}-88466\,{v}^
{20}\\&&-3274\,{v}^{15}+505624\,{v}^{21}-7042\,{v}^{16}+386764\,{v}^{26}+
25480\,{v}^{30}+724351\,{v}^{25}\\&&-8008\,{v}^{31}-396536\,{v}^{24}+30598
\,{v}^{18}-182902\,{v}^{19}+37301\,{v}^{17}-817698\,{v}^{23}\\&&-6\,{v}^{
13}-2\,{v}^{12} \\
  {a}_{5}&=& -12536\,{v}^{19}+1110\,{v}^{16}-4368\,{v}^{31}+226372\,{v}^{25}+3276\,{v}^{30}+60820\,{v}^{20}\\&&-
158480\,{v}^{27}-151912\,{v}^{22}+45304\,{v}^{29}-146870\,{v}^{23}+
188024\,{v}^{24}-202\,{v}^{15}\\&&-90615\,{v}^{26}-13092\,{v}^{18}+51602\,
{v}^{21}+2002\,{v}^{17}+5\,{v}^{14}+5208\,{v}^{28} \\
 {a}_{4}&=& -17162\,{v}^{21}+59508\,{v}^{24}+13132\,{v}^{29}-17367\,{v}^{
27}-324\,{v}^{18}+26320\,{v}^{28}\\&&-67204\,{v}^{26}-2520\,{v}^{30}+3900
\,{v}^{19}-{v}^{15}-1820\,{v}^{31}-21112\,{v}^{22}-16244\,{v}^{25}\\&&+38
\,{v}^{16}-296\,{v}^{17}+34410\,{v}^{23}+3846\,{v}^{20} \\
  {a}^{3}&=& 9180\,{v}^{27}+6872\,{v}^{23}+608\,{v}^{29}-4\,{v}^{17}-6748
\,{v}^{26}+8911\,{v}^{28}-1840\,{v}^{30}\\&&+4276\,{v}^{22}-788\,{v}^{21}+
56\,{v}^{18}-16256\,{v}^{25}-4674\,{v}^{24}-941\,{v}^{20}-12\,{v}^{19}\\&&
-560\,{v}^{31} \\
 {a}_{2}&=& -300\,{v}^{25}+552\,{v}^{28}-861\,{v}^{29}-1892\,{v}^{24}+48\,{v}^{22}-600\,{v}^{30}-998\,{v}^{23}\\&&+
2576\,{v}^{26}-6\,{v}^{19}+20\,{v}^{20}-120\,{v}^{31}+161\,{v}^{21}+
2828\,{v}^{27} \\
  a_{1}&=& 208\,{v}^{24}+28\,{v}^{23}-28\,{v}^{27}-13\,{v}^{22}-105\,{v}^{30}+420\,{v}^{25}-16\,{v}^{31}+350\,{v
}^{26}\\&&-4\,{v}^{21}-312\,{v}^{28}-272\,{v}^{29} \\
a_{0}&=&-8\,{v}^{30}-28\,{v}^{29}-56\,{v}^{28}-70\,{v}^{27}-56\,{v}^{26}-28\,{v}^{25}-{v}^{
31}-{v}^{23}-8\,{v}^{24},
\end{eqnarray*}

\begin{eqnarray*}
b_{32}&=&     28\,{v}^{14}-56\,{v}^{13}+70\,{v}^{12}+28\,{v}^{10}-8\,{v}^{9}+{v}^{8}+{v}^{16}-56\,{v}^{11}-8\,{v}^{15}  \\
b_{31}&=&     -504\,{v}^{15}+44\,{v}^{10}-16\,{v}^{17}+1092\,{v}^{12}+4\,{v}^{7}-1372\,{v}^{13}-28\,{v}^{8}-476\,{v}^{11}\\
&&+1060\,{v}^{14}+136\,{v}^{16}+60\,{v}^{9}   \\
b_{30}&=&    -32\,{v}^{7}+6\,{v}^{6}+120\,{v}^{18}-1096\,{v}^{17}+1104\,{v}^{11}-38\,{v}^{8}-9672\,{v}^{15}+3618\,{v}^{12}\\
&&+13044\,{v}^{14}-1698\,{v}^{10}-10332\,{v}^{13}+636\,{v}^{9}+4340\,{v}^{16} \\
b_{29}&=&    -560\,{v}^{19}+19724\,{v}^{11}+63672\,{v}^{14}+1030\,{v}^{8}+5560\,{v}^{18}-6386\,{v}^{10}-6\,{v}^{6} \\
&&-18810\,{v}^{12}-80188\,{v}^{15}-23752\,{v}^{17}-732\,{v}^{9}-15920\,{v}^{13}-208\,{v}^{7}+4\,{v}^{5}\\
&&+56572\,{v}^{16}  \\
b_{28}&=&    42816\,{v}^{14}+472\,{v}^{7}+1820\,{v}^{20}+2607\,{v}^{8}-212\,{v}^{6}+16258\,{v}^{10}-19880\,{v}^{19}\\
&&+14\,{v}^{5}+38824\,{v}^{11}-141452\,{v}^{12}+149156\,{v}^{13}-13986\,{v}^{9}+{v}^{4}-289296\,{v}^{15}\\
&&-237824\,{v}^{17}+92548\,{v}^{18}+358134\,{v}^{16}   \\
b_{27}&=&    -8968\,{v}^{8}+53144\,{v}^{20}+724530\,{v}^{13}-63704\,{v}^{15}-246\,{v}^{6}-777980\,{v}^{14}\\
&&+114130\,{v}^{10}-4368\,{v}^{21}+10\,{v}^{4}+763308\,{v}^{18}-165014\,{v}^{12}-272496\,{v}^{19}\\
&&-17714\,{v}^{9}-1236312\,{v}^{17}+1038016\,{v}^{16}-150042\,{v}^{11}-80\,{v}^{5}+3786\,{v}^{7}   \\
b_{26}&=&    1188\,{v}^{7}+26098\,{v}^{16}+878532\,{v}^{12}-328\,{v}^{5}+88606\,{v}^{9}-109928\,{v}^{21}\\
&&+545560\,{v}^{13}+2\,{v}^{3}+1936\,{v}^{6}-1938496\,{v}^{19}-{v}^{4}+3422468\,{v}^{18}\\
&&-650252\,{v}^{11}-2860222\,{v}^{14}+627900\,{v}^{20}+8008\,{v}^{22}+2984876\,{v}^{15}   \\
&&+76430\,{v}^{10}-3065468\,{v}^{17}-36909\,{v}^{8}\\
b_{25}&=&   -3708134\,{v}^{13}-27600\,{v}^{17}-1566248\,{v}^{14}+6\,{v}^{3}+3982524\,{v}^{20}+4147\,{v}^{8}\\
&&+179608\,{v}^{22}+9105424\,{v}^{15}-23002\,{v}^{7}+4192\,{v}^{6}-121\,{v}^{4}-1156792\,{v}^{21}\\
&&+272\,{v}^{5}-11440\,{v}^{23}-568524\,{v}^{10}-7771876\,{v}^{19}+239052\,{v}^{9}-8892842\,{v}^{16}\\
&&-230598\,{v}^{11}+2809488\,{v}^{12}+7632464\,{v}^{18}  \\
b_{24}&=&   864528\,{v}^{18}+4310384\,{v}^{15}-82366\,{v}^{9}-24027215\,{v}^{16}-18\,{v}^{3}+12870\,{v}^{24}\\
&&+2631512\,{v}^{11}-234520\,{v}^{23}+{v}^{2}-5260\,{v}^{6}+11994936\,{v}^{14}+21185572\,{v}^{17}\\
&&-106\,{v}^{4}+2130\,{v}^{5}-6708624\,{v}^{21}-16215416\,{v}^{19}+1725724\,{v}^{22}\\
&&+14676914\,{v}^{20}-1143074\,{v}^{10}-31356\,{v}^{7}+170342\,{v}^{8}+556774\,{v}^{12}\\
&&-9683732\,{v}^{13}   \\
b_{23}&=&     -40952812\,{v}^{18}+53386226\,{v}^{17}-11679232\,{v}^{16}-9322816\,{v}^{12}+517\,{v}^{4}\\
&&-4646284\,{v}^{19}-50\,{v}^{3}-2098096\,{v}^{23}-11440\,{v}^{25}+245960\,{v}^{24}-888671\,{v}^{9}\\
&&+9334468\,{v}^{22}+27487598\,{v}^{14}-30742872\,{v}^{15}-1416032\,{v}^{13}-{v}^{2}-19766\,{v}^{6}\\
&&+48564\,{v}^{7}-23222688\,{v}^{21}+518394\,{v}^{10}+163984\,{v}^{8}+4302812\,{v}^{11}\\
&&+29511664\,{v}^{20}+573\,{v}^{5}   \\
b_{22}&=&      -2019118\,{v}^{11}-664340\,{v}^{9}-10755888\,{v}^{23}+2081508\,{v}^{24}+8008\,{v}^{26}\\
&&-207064\,{v}^{25}-13275356\,{v}^{12}+70\,{v}^{3}+29609912\,{v}^{17}+26198096\,{v}^{13}\\
&&+4577880\,{v}^{14}-65491512\,{v}^{15}+63495768\,{v}^{16}-100616676\,{v}^{18}\\
&&+14375986\,{v}^{20}-45961972\,{v}^{21}+30873744\,{v}^{22}-282923\,{v}^{8}+120528\,{v}^{7}\\
&&-315\,{v}^{6}-10\,{v}^{2}-5944\,{v}^{5}+641\,{v}^{4}+3503885\,{v}^{10}+64435092\,{v}^{19}  \\
b_{21}&=&    -v-10900850\,{v}^{11}+1168750\,{v}^{9}-34455260\,{v}^{23}+10247172\,{v}^{24}-4368\,{v}^{27}\\
&&+139048\,{v}^{26}-1679288\,{v}^{25}+5485076\,{v}^{12}+199\,{v}^{3}-106092576\,{v}^{17}\\
&&+34393876\,{v}^{13}-59646086\,{v}^{14}-15387826\,{v}^{15}+132282584\,{v}^{16}-4495\,{v}^{5}\\
&&-66036378\,{v}^{18}+161045812\,{v}^{19}-81715214\,{v}^{20}-31561156\,{v}^{21}-1060\,{v}^{4}\\
&&+61009928\,{v}^{22}-538453\,{v}^{8}-10660\,{v}^{7}+40302\,{v}^{6}+2\,{v}^{2}+2220922\,{v}^{10}\\
b_{20}&=&   -v-6340116\,{v}^{11}+1872142\,{v}^{9}-68612712\,{v}^{23}+32137630\,{v}^{24}+1820\,{v}^{28}\\
&&-73528\,{v}^{27}+1092364\,{v}^{26}-8025472\,{v}^{25}+27491537\,{v}^{12}-82\,{v}^{3}\\
&&-227248510\,{v}^{17}-10572074\,{v}^{13}-75791354\,{v}^{14}+111198616\,{v}^{15}\\
&&+43969209\,{v}^{16}+141781554\,{v}^{18}+124977392\,{v}^{19}-218087144\,{v}^{20}\\
&&+80904348\,{v}^{21}+53062524\,{v}^{22}+60445\,{v}^{8}-188404\,{v}^{7}+22208\,{v}^{6}\\
&&+35\,{v}^{2}+7940\,{v}^{5}-1707\,{v}^{4}-3638656\,{v}^{10}  \\
b_{19}&=&   3\,v+8831766\,{v}^{11}-166605\,{v}^{9}-70701976\,{v}^{23}+64866816\,{v}^{24}-560\,{v}^{29}\\
&&+29960\,{v}^{28}-564592\,{v}^{27}+5113892\,{v}^{26}-24856104\,{v}^{25}+15644490\,{v}^{12}\\
&&-358\,{v}^{3}-101167208\,{v}^{17}-57137352\,{v}^{13}+13203106\,{v}^{14}+142765994\,{v}^{15}\\
&&-169923466\,{v}^{16}+331377762\,{v}^{18}-145635190\,{v}^{19}-197656474\,{v}^{20}\\
&&+247679706\,{v}^{21}-56999452\,{v}^{22}+657527\,{v}^{8}-86682\,{v}^{7}-37671\,{v}^{6}\\
&&+11\,{v}^{2}+9170\,{v}^{5}+779\,{v}^{4}-5247310\,{v}^{10}   \\
b_{18}&=&    3\,v+12104924\,{v}^{11}-1791540\,{v}^{9}+18164908\,{v}^{23}+75811910\,{v}^{24}+120\,{v}^{30}\\
&&-9080\,{v}^{29}+226548\,{v}^{28}-2608704\,{v}^{27}+15750084\,{v}^{26}-51034900\,{v}^{25}\\
&&-17001396\,{v}^{12}-46\,{v}^{3}+210505838\,{v}^{17}-33335864\,{v}^{13}+98733540\,{v}^{14}\\
&&-4613648\,{v}^{15}-229845757\,{v}^{16}+187506122\,{v}^{18}-407360296\,{v}^{19}\\
&&+259448488\,{v}^{21}-232221862\,{v}^{22}+279120\,{v}^{8}+126638\,{v}^{7}-35671\,{v}^{6}\\
&&-3800\,{v}^{5}+1986\,{v}^{4}+207783\,{v}^{10}+100994568\,{v}^{20}-42\,{v}^{2}   \\
b_{17}&=&   1-3\,v+263928\,{v}^{11}-747635\,{v}^{9}+174701288\,{v}^{23}+18056290\,{v}^{24}-16\,{v}^{31}\\
&&+1928\,{v}^{30}-68040\,{v}^{29}+1039636\,{v}^{28}-8034572\,{v}^{27}+32970560\,{v}^{26}\\
&&-65695672\,{v}^{25}-23264720\,{v}^{12}+238\,{v}^{3}+315103514\,{v}^{17}+26037878\,{v}^{13}\\
&&+60987352\,{v}^{14}-142203258\,{v}^{15}-24381701\,{v}^{16}-204549190\,{v}^{18}\\
&&-282493614\,{v}^{19}+416622320\,{v}^{20}-17091534\,{v}^{21}-280959540\,{v}^{22}+147\,{v}^{4}\\
&&-319981\,{v}^{8}+107415\,{v}^{7}+11664\,{v}^{6}-22\,{v}^{2}-7377\,{v}^{5}+3912950\,{v}^{10}\\
b_{16}&=&    -1-3\,v-6960140\,{v}^{11}+{v}^{32}+627666\,{v}^{9}+248565736\,{v}^{23}-99335082\,{v}^{24}\\
&&-256\,{v}^{31}+14392\,{v}^{30}-311736\,{v}^{29}+3216664\,{v}^{28}-17173520\,{v}^{27}\\
&&+45818892\,{v}^{26}-36918604\,{v}^{25}-2003242\,{v}^{12}+84\,{v}^{3}+74262096\,{v}^{17}\\
&&+37504808\,{v}^{13}-31254556\,{v}^{14}-95224600\,{v}^{15}+170002269\,{v}^{16}\\
&&+141641544\,{v}^{19}+347996308\,{v}^{20}-346019408\,{v}^{21}-70547612\,{v}^{22}\\
&&-24156\,{v}^{7}+20170\,{v}^{6}+15\,{v}^{2}-502\,{v}^{5}-815\,{v}^{4}+1662676\,{v}^{10} \\
&&-365300746\,{v}^{18}-258471\,{v}^{8}\\
b_{15}&=&  -3066386\,{v}^{11}+16\,{v}^{32}+504915\,{v}^{9}+122636228\,{v}^{23}-176434880\,{v}^{24}\\
&&-1912\,{v}^{31}+66120\,{v}^{30}-972548\,{v}^{29}+7027580\,{v}^{28}-25410000\,{v}^{27}-222\,{v}^{4}\\
&&+36090656\,{v}^{26}+35748570\,{v}^{25}+10173968\,{v}^{12}-38\,{v}^{3}-166442579\,{v}^{17}\\
&&+5556916\,{v}^{13}-50870814\,{v}^{14}+27859664\,{v}^{15}+126259358\,{v}^{16}+1951\,{v}^{5}\\
&&-130462730\,{v}^{18}+354034466\,{v}^{19}-42929614\,{v}^{20}-350412764\,{v}^{21}\\
&&+222022066\,{v}^{22}+33333\,{v}^{8}-42616\,{v}^{7}+1657\,{v}^{6}+9\,{v}^{2}-971477\,{v}^{10}\\
b_{14}&=&    1191220\,{v}^{11}+120\,{v}^{32}-25014\,{v}^{9}-95976776\,{v}^{23}-123541870\,{v}^{24}-8840\,{v}^{31}\\
&&+208628\,{v}^{30}-2173472\,{v}^{29}+10964268\,{v}^{28}-24244500\,{v}^{27}-173090\,{v}^{26}\\
&&+97005340\,{v}^{25}+4696087\,{v}^{12}-13\,{v}^{3}-141492568\,{v}^{17}-12270080\,{v}^{13}\\
&&-10472648\,{v}^{14}+58060448\,{v}^{15}-14948528\,{v}^{16}+129314739\,{v}^{18}\\
&&-281246082\,{v}^{20}-48798384\,{v}^{21}+285939796\,{v}^{22}+71479\,{v}^{8}-4404\,{v}^{7}\\
&&-3472\,{v}^{6}+438\,{v}^{5}+62\,{v}^{4}-807966\,{v}^{10}+168736722\,{v}^{19}    \\
b_{13}&=&  1065769\,{v}^{11}+560\,{v}^{32}-96701\,{v}^{9}-185021194\,{v}^{23}+10548458\,{v}^{24}\\
&&-28280\,{v}^{31}+477232\,{v}^{30}-3551492\,{v}^{29}+11878048\,{v}^{28}-10166912\,{v}^{27}\\
&&-38156200\,{v}^{26}+88536132\,{v}^{25}-1146395\,{v}^{12}-2670142\,{v}^{17}-5989506\,{v}^{13}\\
&&+12199926\,{v}^{14}+15060662\,{v}^{15}-55591808\,{v}^{16}+133232930\,{v}^{18}\\
&&-73435746\,{v}^{19}-170979128\,{v}^{20}+176293794\,{v}^{21}+97537978\,{v}^{22}+8892\,{v}^{8}\\
&&+4709\,{v}^{7}-666\,{v}^{6}-75\,{v}^{5}+11\,{v}^{4}-8710\,{v}^{10}   \\
b_{12}&=&    56900\,{v}^{11}+1820\,{v}^{32}-13652\,{v}^{9}-96531784\,{v}^{23}+90390696\,{v}^{24}\\
&&-66248\,{v}^{31}+812812\,{v}^{30}-4220944\,{v}^{29}+7876302\,{v}^{28}+8161272\,{v}^{27}\\
&&-46401484\,{v}^{26}+23606932\,{v}^{25}-1164895\,{v}^{12}+44334350\,{v}^{17}+837462\,{v}^{13}\\
&&+6382888\,{v}^{14}-9933836\,{v}^{15}-17282707\,{v}^{16}+17251146\,{v}^{18}-104462048\,{v}^{19}\\
&&+21517601\,{v}^{20}+138227802\,{v}^{21}-79488994\,{v}^{22}-4905\,{v}^{8}+786\,{v}^{7}+71\,{v}^{6}\\
&&-5\,{v}^{5}+106640\,{v}^{10}   \\
b_{11}&=&    -96702\,{v}^{11}+4368\,{v}^{32}+3914\,{v}^{9}+17373726\,{v}^{23}+66029484\,{v}^{24}+16151\,{v}^{10}\\
&&-117208\,{v}^{31}+1038648\,{v}^{30}-3486548\,{v}^{29}+1141980\,{v}^{28}+17031224\,{v}^{27}\\
&&-23522796\,{v}^{26}-29137622\,{v}^{25}-94384\,{v}^{12}+16172324\,{v}^{17}+1060077\,{v}^{13}\\
&&-418010\,{v}^{14}-5700550\,{v}^{15}+6514790\,{v}^{16}-29033828\,{v}^{18}-23395158\,{v}^{19}\\
&&+67129174\,{v}^{20}+10292448\,{v}^{21}-88804992\,{v}^{22}-711\,{v}^{8}-50\,{v}^{7}+{v}^{6}\\
b_{10}&=&    -14900\,{v}^{11}+8008\,{v}^{32}+482\,{v}^{9}+44121912\,{v}^{23}+7945612\,{v}^{24}-159016\,{v}^{31}\\
&&+983268\,{v}^{30}-1681680\,{v}^{29}-3664760\,{v}^{28}+12534620\,{v}^{27}+2590338\,{v}^{26}\\
&&-32544884\,{v}^{25}+72777\,{v}^{12}-3309234\,{v}^{17}+102244\,{v}^{13}-806923\,{v}^{14}\\
&&+75864\,{v}^{15}+4277033\,{v}^{16}-12478116\,{v}^{18}+15169480\,{v}^{19}+21112516\,{v}^{20}\\
&&-34289496\,{v}^{21}-20032024\,{v}^{22}+24\,{v}^{8}-2361\,{v}^{10}   \\
b_{9}&=&   -14430938\,{v}^{21}-165880\,{v}^{31}-15718870\,{v}^{24}-82244\,{v}^{14}-10573820\,{v}^{25}\\
&&+656656\,{v}^{30}-45958\,{v}^{13}-4156664\,{v}^{28}+7976928\,{v}^{19}+1039\,{v}^{11}\\
&&-2700658\,{v}^{17}+10952164\,{v}^{26}+3202976\,{v}^{27}+96051\,{v}^{16}+18876\,{v}^{29}\\
&&+11440\,{v}^{32}+516220\,{v}^{15}-5900620\,{v}^{20}-7\,{v}^{9}+12965200\,{v}^{22}-238\,{v}^{10}\\
&&+10795\,{v}^{12}+16198224\,{v}^{23}+1163332\,{v}^{18} \\
b_{8}&=&  51100\,{v}^{15}+5858048\,{v}^{26}-124070\,{v}^{17}+7720424\,{v}^{22}+1326260\,{v}^{21}-302\,{v}^{12}\\
&&-131560\,{v}^{31}+24632\,{v}^{14}+3034916\,{v}^{25}+261404\,{v}^{30}-2040878\,{v}^{28}\\
&&+1436802\,{v}^{18}+{v}^{10}-4225210\,{v}^{20}-1929032\,{v}^{27}+738816\,{v}^{29}-2822832\,{v}^{23}\\
&&-142456\,{v}^{19}+82\,{v}^{11}-279006\,{v}^{16}+12870\,{v}^{32}-8773711\,{v}^{24}-6150\,{v}^{13} \\
b_{7}&=&   590772\,{v}^{29}+1845668\,{v}^{21}-3221120\,{v}^{23}+3640\,{v}^{30}+86352\,{v}^{18}+2740\,{v}^{14}\\
&&+37\,{v}^{13}-340600\,{v}^{24}-24674\,{v}^{16}-11280\,{v}^{15}+128185\,{v}^{17}+190066\,{v}^{22}\\
&&-1940960\,{v}^{27}-76648\,{v}^{31}-18\,{v}^{12}+426264\,{v}^{26}-150126\,{v}^{20}+11440\,{v}^{32}\\
&&+3335102\,{v}^{25}-199724\,{v}^{28}-644396\,{v}^{19}\\
b_{6}&=&   243722\,{v}^{20}-657536\,{v}^{22}+140068\,{v}^{21}+4384\,{v}^{16}-70980\,{v}^{30}-940\,{v}^{15}\\
&&+213472\,{v}^{29}-814966\,{v}^{26}+8008\,{v}^{32}-43134\,{v}^{19}+1005722\,{v}^{24}-50393\,{v}^{18}\\
&&+9044\,{v}^{17}-545244\,{v}^{27}-364568\,{v}^{23}+11\,{v}^{14}+333788\,{v}^{28}-29848\,{v}^{31}\\
&&+622908\,{v}^{25}+2\,{v}^{13}  \\
b_{5}&=&  184800\,{v}^{28}-5096\,{v}^{31}+240\,{v}^{16}-47824\,{v}^{30}+73248\,{v}^{27}+186322\,{v}^{23}\\
&&+7924\,{v}^{29}+16382\,{v}^{20}-1406\,{v}^{17}-72714\,{v}^{22}-204268\,{v}^{25}+16992\,{v}^{19}\\
&&+4368\,{v}^{32}-2326\,{v}^{18}+206378\,{v}^{24}-77982\,{v}^{21}-6\,{v}^{15}-293576\,{v}^{26}  \\
b_{4}&=&  26278\,{v}^{28}-42\,{v}^{17}+1960\,{v}^{31}-4634\,{v}^{21}+312\,{v}^{19}-4841\,{v}^{20}-25712\,{v}^{29}\\
&&+27984\,{v}^{23}+21438\,{v}^{22}-75764\,{v}^{25}+352\,{v}^{18}+76384\,{v}^{27}-14972\,{v}^{30}\\
&&+9496\,{v}^{26}+1820\,{v}^{32}+{v}^{16}-39084\,{v}^{24}  \\
b_{3}&=&    -62\,{v}^{19}+4374\,{v}^{25}+4\,{v}^{18}-5274\,{v}^{23}+836\,{v}^{22}+1720\,{v}^{31}-9772\,{v}^{29}\\
&&+32\,{v}^{20}-1208\,{v}^{30}-8764\,{v}^{24}+9576\,{v}^{27}+560\,{v}^{32}-8628\,{v}^{28}+1102\,{v}^{21}\\
&&+18832\,{v}^{26} \\
b_{2}&=&  -2604\,{v}^{27}-174\,{v}^{22}+1206\,{v}^{24}+6\,{v}^{20}-824\,{v}^{29}+650\,{v}^{26}+584\,{v}^{31}\\
&&-3140\,{v}^{28}-24\,{v}^{21}+2312\,{v}^{25}+120\,{v}^{32}-20\,{v}^{23}+756\,{v}^{30}  \\
b_{1}&=&  4\,{v}^{22}-420\,{v}^{27}+264\,{v}^{30}-476\,{v}^{26}-28\,{v}^{28}+284\,{v}^{29}+104\,{v}^{31}-36\,{v}^{24}\\
&&-236\,{v}^{25}+12\,{v}^{23}+16\,{v}^{32}  \\
b_{0}&=&{v}^{32}+{v}^{24}+8\,{v}^{31}+28\,{v}^{30}+56\,{v}^{29}+70\,{v}^{28}+56\,{v}^{27}+28\,{v}^{26}+8\,{v}^{25},
\end{eqnarray*}

\begin{eqnarray*}
c_{63}&=&1     \\
c_{62}&=&y    \\
c_{61}&=&  -16-32\,y-31\,{y}^{2}    \\
c_{60}&=&  -16\,y-32\,{y}^{2}-31\,{y}^{3}    \\
c_{59}&=&    960\,{y}^{2}+120+960\,{y}^{3}+496\,y+465\,{y}^{4}    \\
c_{58}&=&     496\,{y}^{2}+960\,{y}^{4}+465\,{y}^{5}+120\,y+960\,{y}^{3}  \\
c_{57}&=&     -20880\,{y}^{4}-568-13920\,{y}^{5}-3600\,y-10680\,{y}^{2}-4495\,{y}^{6}-18864\,{y}^{3}   \\
c_{56}&=&   -3600\,{y}^{2}-568\,y-10680\,{y}^{3}-18864\,{y}^{4}-13920\,{y}^{6}-20880\,{y}^{5}-4495\,{y}^{7} \\
c_{55}&=&     259840\,{y}^{6}+1932+66416\,{y}^{2}+16472\,y+279440\,{y}^{4}+129920\,{y}^{7}+31465\,{y}^{8}\\
&&+165760\,{y}^{3}+326816\,{y}^{5}  \\
c_{54}&=&    66416\,{y}^{3}+165760\,{y}^{4}+31465\,{y}^{9}+279440\,{y}^{5}+326816\,{y}^{6}+259840\,{y}^{7}\\
&&+1932\,y+129920\,{y}^{8}+16472\,{y}^{2}  \\
c_{53}&=&    -275268\,{y}^{2}-2192400\,{y}^{8}-169911\,{y}^{10}-3946320\,{y}^{6}-884680\,{y}^{3}\\
&&-876960\,{y}^{9}-54096\,y-1986208\,{y}^{4}-5096-3254496\,{y}^{5}-3505824\,{y}^{7}\\
c_{52}&=&    -169911\,{y}^{11}-884680\,{y}^{4}-3946320\,{y}^{7}-2192400\,{y}^{9}-1986208\,{y}^{5}\\
&&-3505824\,{y}^{8}-3254496\,{y}^{6}-5096\,y-275268\,{y}^{3}-54096\,{y}^{2}-876960\,{y}^{10} \\
c_{51}&=&    736281\,{y}^{12}+19106984\,{y}^{5}+30702672\,{y}^{6}+26456976\,{y}^{9}+13680576\,{y}^{10}\\
&&+3279360\,{y}^{3}+838768\,{y}^{2}+9140404\,{y}^{4}+137592\,y+36658440\,{y}^{8}+10948\\
&&+4560192\,{y}^{11}+38228736\,{y}^{7}  \\
c_{50}&=&   30702672\,{y}^{7}+38228736\,{y}^{8}+36658440\,{y}^{9}+19106984\,{y}^{6}+736281\,{y}^{13}\\
&&+3279360\,{y}^{4}+838768\,{y}^{3}+137592\,{y}^{2}+26456976\,{y}^{10}+9140404\,{y}^{5}\\
&&+4560192\,{y}^{12}+10948\,y+13680576\,{y}^{11} \\
c_{49}&=& -239798520\,{y}^{7}-304686200\,{y}^{8}-30215528\,{y}^{4}-76260464\,{y}^{5}-19788\\
&&-66502800\,{y}^{12}-151361756\,{y}^{6}-150630480\,{y}^{11}-2002644\,{y}^{2}-2629575\,{y}^{14}\\
&&-309145200\,{y}^
{9}-246846600\,{y}^{10}-9144408\,{y}^{3}-284648\,y-19000800\,{y}^{13}\\
c_{48}&=&    -309145200\,{y}^{10}-151361756\,{y}^{7}-2002644\,{y}^{3}-304686200\,{y}^{9}-284648\,{y}^{2}\\
&&-66502800\,{y}^{13}-246846600\,{y}^{11}-76260464\,{y}^{6}-19000800\,{y}^{14}-19788\,y\\
&&-150630480\,{y}^{12}-2629575\,{y}^{15}-239798520\,{y}^{8}-9144408\,{y}^{4}-30215528\,{y}^{5}  \\
c_{47}&=&   7888725\,{y}^{16}+1005396736\,{y}^{7}+527024064\,{y}^{6}+2034153680\,{y}^{9}\\
&&+224838432\,{y}^{5}+20243552\,{y}^{3}+3911880\,{y}^{2}+2156030240\,{y}^{10}\\
&&+1575645864\,{y}^{8}+675904320\,{y}^{13}+1857273600\,{y}^{11}+76666576\,{y}^{4}\\
&&+65145600\,{y}^{15}+494700\,y+30782+1275674400\,{y}^{12}+260582400\,{y}^{14}   \\
c_{46}&=&    30782\,y+1275674400\,{y}^{13}+20243552\,{y}^{4}+675904320\,{y}^{14}+2034153680\,{y}^{10}\\
&&+1005396736\,{y}^{8}+260582400\,{y}^{15}+1857273600\,{y}^{12}+527024064\,{y}^{7}\\
&&+494700\,{y}^{2}+76666576\,{y}^{5}+2156030240\,{y}^{11}+3911880\,{y}^{3}+1575645864\,{y}^{9}\\
&&+65145600\,{y}^{16}+224838432\,{y}^{6}+7888725\,{y}^{17}  \\
c_{45}&=&   -41944-12651024880\,{y}^{11}-11565415680\,{y}^{12}-517756288\,{y}^{5}\\
&&-3069905696\,{y}^{7}-5674585584\,{y}^{8}-2461942080\,{y}^{15}-36935732\,{y}^{3}\\
&&-738768\,y-8690136000\,{y}^{13}-842821200\,{y}^{16}-1387485456\,{y}^{6}-6437842\,{y}^{2}\\
&&-8815886944\,{y}^{9}-20160075\,{y}^{18}-187293600\,{y}^{17}-156342280\,{y}^{4}\\
&&-5250679200\,{y}^{14}-11527653176\,{y}^{10}  \\
c_{44}&=&  -5250679200\,{y}^{15}-12651024880\,{y}^{12}-517756288\,{y}^{6}-2461942080\,{y}^{16}\\
&&-842821200\,{y}^{17}-11527653176\,{y}^{11}-3069905696\,{y}^{8}-1387485456\,{y}^{7}\\
&&-8690136000\,{y}^{14}-11565415680\,{y}^{13}-156342280\,{y}^{5}-8815886944\,{y}^{10}\\
&&-5674585584\,{y}^{9}-187293600\,{y}^{18}-36935732\,{y}^{4}-20160075\,{y}^{19}\\
&&-6437842\,{y}^{3}-738768\,{y}^{2}-41944\,y \\
c_{43}&=&    264696582\,{y}^{4}+15262104264\,{y}^{8}+42393023200\,{y}^{10}+9115712\,{y}^{2}+964712\,y\\
&&+2896080616\,{y}^{6}+968788028\,{y}^{5}+7227407968\,{y}^{7}+60648896880\,{y}^{13}\\
&&+48925059040\,{y}^{14}+457828800\,{y}^{19}+56009848576\,{y}^{11}+27488903728\,{y}^{9}\\
&&+7433403120\,{y}^{17}+32685818880\,{y}^{15}+57014784\,{y}^{3}+63235287544\,{y}^{12}\\
&&+17659956600\,{y}^{16}+44352165\,{y}^{20}+2289144000\,{y}^{18}+50788  \\
c_{42}&=&  7227407968\,{y}^{8}+50788\,y+63235287544\,{y}^{13}+964712\,{y}^{2}+48925059040\,{y}^{15}\\
&&+60648896880\,{y}^{14}+56009848576\,{y}^{12}+7433403120\,{y}^{18}+15262104264\,{y}^{9}\\
&&+44352165\,{y}^{21}+17659956600\,{y}^{17}+2896080616\,{y}^{7}+9115712\,{y}^{3}\\
&&+2289144000\,{y}^{19}+42393023200\,{y}^{11}+27488903728\,{y}^{10}+457828800\,{y}^{20}\\
&&+264696582\,{y}^{5}+57014784\,{y}^{4}+968788028\,{y}^{6}+32685818880\,{y}^{16} \\
c_{41}&=&    -961440480\,{y}^{21}-235104125872\,{y}^{12}-13663739684\,{y}^{7}-4958853354\,{y}^{6}\\
&&-18891752880\,{y}^{19}-231728994640\,{y}^{15}-1117336\,y-176458354800\,{y}^{11}\\
&&-115163329512\,{y}^{10}-32174130384\,{y}^{8}-271728739744\,{y}^{13}\\
&&-101056235280\,{y}^{17}-11289684\,{y}^{2}-167711503960\,{y}^{16}-75898184\,{y}^{3}\\
&&-49439590200\,{y}^{18}-84672315\,{y}^{22}-1515101520\,{y}^{5}-55308-380782360\,{y}^{4}\\
&&-65357234096\,{y}^{9}-5287922640\,{y}^{20}-271122347496\,{y}^{14}  \\
c_{40}&=&   -49439590200\,{y}^{19}-5287922640\,{y}^{21}-167711503960\,{y}^{17}-101056235280\,{y}^{18}\\
&&-1515101520\,{y}^{6}-75898184\,{y}^{4}-380782360\,{y}^{5}-18891752880\,{y}^{20}\\
&&-115163329512\,{y}^{11}-65357234096\,{y}^{10}-55308\,y-231728994640\,{y}^{16}\\
&&-271122347496\,{y}^{15}-961440480\,{y}^{22}-13663739684\,{y}^{8}-176458354800\,{y}^{12}\\
&&-32174130384\,{y}^{9}-4958853354\,{y}^{7}-84672315\,{y}^{23}-11289684\,{y}^{3}\\
&&-1117336\,{y}^{2}-235104125872\,{y}^{13}-271728739744\,{y}^{14} \\
c_{39}&=&   54746+1161468\,y+12395472\,{y}^{2}+141120525\,{y}^{24}+1748073600\,{y}^{23}\\
&&+10488441600\,{y}^{22}+40877796960\,{y}^{21}+116781865200\,{y}^{20}\\
&&+261025564800\,{y}^{19}+475001326800\,{y}^{18}+722381883080\,{y}^{17}\\
&&+473596296\,{y}^{4}+88466928\,{y}^{3}+55058724882\,{y}^{8}+21364460416\,{y}^{7}\\
&&+7138798144\,{y}^{6}+2021478160\,{y}^{5}+242729271928\,{y}^{10}+857514144032\,{y}^{13}\\
&&+123496243892\,{y}^{9}+639258711056\,{y}^{12}+419669647264\,{y}^{11}\\
&&+934774965124\,{y}^{16}+1042368845056\,{y}^{15}+1010698496192\,{y}^{14}\\
c_{38}&=&   54746\,y+1161468\,{y}^{2}+141120525\,{y}^{25}+1748073600\,{y}^{24}+10488441600\,{y}^{23}\\
&&+40877796960\,{y}^{22}+116781865200\,{y}^{21}+261025564800\,{y}^{20}\\
&&+475001326800\,{y}^{19}+722381883080\,{y}^{18}+934774965124\,{y}^{17}\\
&&+88466928\,{y}^{4}+12395472\,{y}^{3}+21364460416\,{y}^{8}+7138798144\,{y}^{7}\\
&&+2021478160\,{y}^{6}+473596296\,{y}^{5}+123496243892\,{y}^{10}+55058724882\,{y}^{9}\\
&&+419669647264\,{y}^{12}+242729271928\,{y}^{11}+1042368845056\,{y}^{16}\\
&&+1010698496192\,{y}^{15}+857514144032\,{y}^{14}+639258711056\,{y}^{13}  \\
c_{37}&=&  -49700-1094920\,y-12201950\,{y}^{2}-206253075\,{y}^{26}-2767783200\,{y}^{25}\\
&&-17990590800\,{y}^{24}-75951788640\,{y}^{23}-235118410800\,{y}^{22}\\
&&-570067898400\,{y}^{21}-1127467621280\,{y}^{20}-1868663373720\,{y}^{19}\\
&&-2644651365356\,{y}^{18}-3239549996144\,{y}^{17}-516409928\,{y}^{4}\\
&&-91427092\,{y}^{3}-78464517736\,{y}^{8}-28228295952\,{y}^{7}-8799133768\,{y}^{6}\\
&&-2337402240\,{y}^{5}-411333491870\,{y}^{10}-191153981104\,{y}^{9}\\
&&-1335956846424\,{y}^{12}-785717717036\,{y}^{11}-3468792768056\,{y}^{16}\\
&&-3270491250208\,{y}^{15}-2729493910864\,{y}^{14}-2023781557056\,{y}^{13} \\
c_{36}&=&    -49700\,y-1094920\,{y}^{2}-206253075\,{y}^{27}-2767783200\,{y}^{26}\\
&&-17990590800\,{y}^{25}-75951788640\,{y}^{24}-235118410800\,{y}^{23}\\
&&-570067898400\,{y}^{22}-1127467621280\,{y}^{21}-1868663373720\,{y}^{20}\\
&&-2644651365356\,{y}^{19}-3239549996144\,{y}^{18}-3468792768056\,{y}^{17}\\
&&-91427092\,{y}^{4}-12201950\,{y}^{3}-28228295952\,{y}^{8}-8799133768\,{y}^{7}\\
&&-2337402240\,{y}^{6}-516409928\,{y}^{5}-191153981104\,{y}^{10}-78464517736\,{y}^{9}\\
&&-785717717036\,{y}^{12}-411333491870\,{y}^{11}-3270491250208\,{y}^{16}\\
&&-2729493910864\,{y}^{15}-2023781557056\,{y}^{14}-1335956846424\,{y}^{13}    \\
c_{35}&=&    41658+944300\,y+10875816\,{y}^{2}+265182525\,{y}^{28}+3832315200\,{y}^{27}\\
&&+26826206400\,{y}^{26}+121944710160\,{y}^{25}+406506227400\,{y}^{24}\\
&&+1062065222400\,{y}^{23}+2266467672560\,{y}^{22}+4061277591800\,{y}^{21}\\
&&+6230647907484\,{y}^{20}+8300632978176\,{y}^{19}+9704772508816\,{y}^{18}\\
&&+10037591056584\,{y}^{17}+499392190\,{y}^{4}+84663024\,{y}^{3}+94890173796\,{y}^{8}\\
&&+32030345872\,{y}^{7}+9418433808\,{y}^{6}+2372198484\,{y}^{5}+575308646016\,{y}^{10}\\
&&+247809863848\,{y}^{9}+2224609904090\,{y}^{12}+1194203144832\,{y}^{11}\\
&&+9240215595324\,{y}^{16}+7605354581664\,{y}^{15}+5615059610136\,{y}^{14}\\
&&+3726298760324\,{y}^{13}   \\
c_{34}&=&    41658\,y+944300\,{y}^{2}+265182525\,{y}^{29}+3832315200\,{y}^{28}+26826206400\,{y}^{27}\\
&&+121944710160\,{y}^{26}+406506227400\,{y}^{25}+1062065222400\,{y}^{24}\\
&&+2266467672560\,{y}^{23}+4061277591800\,{y}^{22}+6230647907484\,{y}^{21}\\
&&+8300632978176\,{y}^{20}+9704772508816\,{y}^{19}+10037591056584\,{y}^{18}\\
&&+9240215595324\,{y}^{17}+84663024\,{y}^{4}+10875816\,{y}^{3}+32030345872\,{y}^{8}\\
&&+9418433808\,{y}^{7}+2372198484\,{y}^{6}+499392190\,{y}^{5}+247809863848\,{y}^{10}\\
&&+94890173796\,{y}^{9}+1194203144832\,{y}^{12}+575308646016\,{y}^{11}\\
&&+5615059610136\,{y}^{15}+3726298760324\,{y}^{14}+2224609904090\,{y}^{13}   \\
&&+7605354581664\,{y}^{16}\\
c_{33}&=&     -32398-749844\,y-8849490\,{y}^{2}-300540195\,{y}^{30}-4653525600\,{y}^{29}\\
&&-34901442000\,{y}^{28}-169949964240\,{y}^{27}-606843808200\,{y}^{26}\\
&&-1698913048560\,{y}^{25}-3888387037080\,{y}^{24}-7483593343720\,{y}^{23}\\
&&-12355689365396\,{y}^{22}-17758908907344\,{y}^{21}-22468474332920\,{y}^{20}\\
&&-25237557100392\,{y}^{19}-25335391302956\,{y}^{18}-22849211615160\,{y}^{17}\\
&&-432498180\,{y}^{4}-70894108\,{y}^{3}-98786488380\,{y}^{8}-31612018812\,{y}^{7}\\
&&-8856044570\,{y}^{6}-2135472168\,{y}^{5}-677034433908\,{y}^{10}-273541049032\,{y}^{9}\\
&&-3031077040872\,{y}^{12}-1507215783752\,{y}^{11}-18587267839548\,{y}^{16}\\
&&-13678571283612\,{y}^{15}-9124132183414\,{y}^{14}-5521648486832\,{y}^{13}  \\
c_{32}&=&   -32398\,y-749844\,{y}^{2}-300540195\,{y}^{31}-4653525600\,{y}^{30}-34901442000\,{y}^{29}\\
&&-169949964240\,{y}^{28}-606843808200\,{y}^{27}-1698913048560\,{y}^{26}\\
&&-3888387037080\,{y}^{25}-7483593343720\,{y}^{24}-12355689365396\,{y}^{23}\\
&&-17758908907344\,{y}^{22}-22468474332920\,{y}^{21}-25237557100392\,{y}^{20}\\
&&-25335391302956\,{y}^{19}-22849211615160\,{y}^{18}-18587267839548\,{y}^{17}\\
&&-70894108\,{y}^{4}-8849490\,{y}^{3}-31612018812\,{y}^{8}-8856044570\,{y}^{7}\\
&&-2135472168\,{y}^{6}-432498180\,{y}^{5}-273541049032\,{y}^{10}-98786488380\,{y}^{9}\\
&&-1507215783752\,{y}^{12}-677034433908\,{y}^{11}-13678571283612\,{y}^{16}\\
&&-9124132183414\,{y}^{15}-5521648486832\,{y}^{14}-3031077040872\,{y}^{13}    \\
c_{31}&=&   23461+550766\,y+6612316\,{y}^{2}+300540195\,{y}^{32}+4963760640\,{y}^{31}\\
&&+39710085120\,{y}^{30}+206210021760\,{y}^{29}+785108606400\,{y}^{28}\\
&&+2343979061760\,{y}^{27}+5724440317120\,{y}^{26}+11768038794720\,{y}^{25}\\
&&+20784246546544\,{y}^{24}+32017471806976\,{y}^{23}+43517224083328\,{y}^{22}\\
&&+52656189618240\,{y}^{21}+57125944786848\,{y}^{20}+55882561771456\,{y}^{19}\\
&&+49516167229392\,{y}^{18}+39883640561464\,{y}^{17}+338249088\,{y}^{4}\\
&&+54081072\,{y}^{3}+89545258664\,{y}^{8}+27418986624\,{y}^{7}+7385214272\,{y}^{6}\\
&&+1720386848\,{y}^{5}+680135364160\,{y}^{10}+260380673904\,{y}^{9}\\
&&+3445948191456\,{y}^{12}+1606265391040\,{y}^{11}+29281089322476\,{y}^{16}\\
&&+19630028733952\,{y}^{15}+12028540910848\,{y}^{14}+6737207143360\,{y}^{13}   \\
c_{30}&=&    23461\,y+550766\,{y}^{2}+4963760640\,{y}^{32}+39710085120\,{y}^{31}\\
&&+206210021760\,{y}^{30}+785108606400\,{y}^{29}+2343979061760\,{y}^{28}\\
&&+5724440317120\,{y}^{27}+11768038794720\,{y}^{26}+20784246546544\,{y}^{25}\\
&&+32017471806976\,{y}^{24}+43517224083328\,{y}^{23}+52656189618240\,{y}^{22}\\
&&+57125944786848\,{y}^{21}+55882561771456\,{y}^{20}+49516167229392\,{y}^{19}\\
&&+39883640561464\,{y}^{18}+29281089322476\,{y}^{17}+300540195\,{y}^{33}\\
&&+54081072\,{y}^{4}+6612316\,{y}^{3}+27418986624\,{y}^{8}+7385214272\,{y}^{7}\\
&&+1720386848\,{y}^{6}+338249088\,{y}^{5}+260380673904\,{y}^{10}+89545258664\,{y}^{9}\\
&&+1606265391040\,{y}^{12}+680135364160\,{y}^{11}+19630028733952\,{y}^{16}\\
&&+12028540910848\,{y}^{15}+6737207143360\,{y}^{14}+3445948191456\,{y}^{13}   \\
c_{29}&=&   -15864-375376\,y-4555563\,{y}^{2}-39554967600\,{y}^{32}-218191512960\,{y}^{31}\\
&&-882244267200\,{y}^{30}-2797321603200\,{y}^{29}-7257693203200\,{y}^{28}\\
&&-15862083433120\,{y}^{27}-29816683632336\,{y}^{26}-48957289784384\,{y}^{25}\\
&&-71053495586208\,{y}^{24}-92003796675648\,{y}^{23}-107080614147104\,{y}^{22}\\
&&-112699479483904\,{y}^{21}-107790255861456\,{y}^{20}-94064767870216\,{y}^{19}\\
&&-75137868974516\,{y}^{18}-55073622939168\,{y}^{17}-265182525\,{y}^{34}\\\
&&-4653525600\,{y}^{33}-240447060\,{y}^{4}-37778802\,{y}^{3}-71327065296\,{y}^{8}\\
&&-21079547040\,{y}^{7}-5505704480\,{y}^{6}-1249565808\,{y}^{5}-589797314744\,{y}^{10}\\
&&-215898531488\,{y}^{9}-3316370807392\,{y}^{12}-1463762798544\,{y}^{11}\\
&&-37104707404400\,{y}^{16}-23000462824896\,{y}^{15}-13120524730080\,{y}^{14}\\
&&-6883455615616\,{y}^{13} \\
c_{28}&=&    -15864\,y-375376\,{y}^{2}-218191512960\,{y}^{32}-882244267200\,{y}^{31}\\
&&-2797321603200\,{y}^{30}-7257693203200\,{y}^{29}-15862083433120\,{y}^{28}\\
&&-29816683632336\,{y}^{27}-48957289784384\,{y}^{26}-71053495586208\,{y}^{25}\\
&&-92003796675648\,{y}^{24}-107080614147104\,{y}^{23}-112699479483904\,{y}^{22}\\
&&-107790255861456\,{y}^{21}-94064767870216\,{y}^{20}-75137868974516\,{y}^{19}\\
&&-55073622939168\,{y}^{18}-37104707404400\,{y}^{17}-265182525\,{y}^{35}\\
&&-4653525600\,{y}^{34}-39554967600\,{y}^{33}-37778802\,{y}^{4}-4555563\,{y}^{3}\\
&&-21079547040\,{y}^{8}-5505704480\,{y}^{7}-1249565808\,{y}^{6}-240447060\,{y}^{5}\\
&&-215898531488\,{y}^{10}-71327065296\,{y}^{9}-1463762798544\,{y}^{12}\\
&&-589797314744\,{y}^{11}-23000462824896\,{y}^{16}-13120524730080\,{y}^{15}\\
&&-6883455615616\,{y}^{14}-3316370807392\,{y}^{13}    \\
c_{27}&=&   10068+237960\,y+2901504\,{y}^{2}+861824917800\,{y}^{32}+2891533209600\,{y}^{31}\\
&&+7939784758080\,{y}^{30}+18373993554080\,{y}^{29}+36600446098384\,{y}^{28}\\
&&+63755034646016\,{y}^{27}+98303363102784\,{y}^{26}+135461150729632\,{y}^{25}\\
&&+168115942730736\,{y}^{24}+189102273515456\,{y}^{23}+193798659736592\,{y}^{22}\\
&&+181743697614808\,{y}^{21}+156524571863548\,{y}^{20}+124160351514112\,{y}^{19}\\
&&+90915946208640\,{y}^{18}+61553704697040\,{y}^{17}+206253075\,{y}^{36}\\
&&+3832315200\,{y}^{35}+34490836800\,{y}^{34}+201401378640\,{y}^{33}+156030381\,{y}^{4}\\
&&+24248704\,{y}^{3}+50301165160\,{y}^{8}+14464607280\,{y}^{7}+3692732284\,{y}^{6}\\
&&+822736158\,{y}^{5}+445393105888\,{y}^{10}+157195789520\,{y}^{9}\\
&&+2731352951576\,{y}^{12}+1151711923264\,{y}^{11}+38569662181960\,{y}^{16}\\
&&+22374507261760\,{y}^{15}+12012381259328\,{y}^{14}+5962478863376\,{y}^{13}   \\
c_{26}&=&    10068\,y+237960\,{y}^{2}+2891533209600\,{y}^{32}+7939784758080\,{y}^{31}\\
&&+18373993554080\,{y}^{30}+36600446098384\,{y}^{29}+63755034646016\,{y}^{28}\\
&&+98303363102784\,{y}^{27}+135461150729632\,{y}^{26}+168115942730736\,{y}^{25}\\
&&+189102273515456\,{y}^{24}+193798659736592\,{y}^{23}+181743697614808\,{y}^{22}\\
&&+156524571863548\,{y}^{21}+124160351514112\,{y}^{20}+90915946208640\,{y}^{19}\\
&&+61553704697040\,{y}^{18}+38569662181960\,{y}^{17}+206253075\,{y}^{37}\\
&&+3832315200\,{y}^{36}+34490836800\,{y}^{35}+201401378640\,{y}^{34}\\
&&+861824917800\,{y}^{33}+24248704\,{y}^{4}+2901504\,{y}^{3}+14464607280\,{y}^{8}\\
&&+3692732284\,{y}^{7}+822736158\,{y}^{6}+156030381\,{y}^{5}+157195789520\,{y}^{10}\\
&&+50301165160\,{y}^{9}+1151711923264\,{y}^{12}+445393105888\,{y}^{11}\\
&&+22374507261760\,{y}^{16}+12012381259328\,{y}^{15}+5962478863376\,{y}^{14}\\
&&+2731352951576\,{y}^{13}  \\
c_{25}&=&   -6036-140952\,y-1711748\,{y}^{2}-7497840382760\,{y}^{32}-18307707855840\,{y}^{31}\\
&&-38500788815216\,{y}^{30}-70862893244864\,{y}^{29}-115577732666144\,{y}^{28}\\
&&-168697385650976\,{y}^{27}-222115579375536\,{y}^{26}-265542242629664\,{y}^{25}\\
&&-289833757036928\,{y}^{24}-290151785862120\,{y}^{23}-267452722477668\,{y}^{22}\\
&&-227729992194848\,{y}^{21}-179594286555632\,{y}^{20}-131450824320400\,{y}^{19}\\
&&-89431228294824\,{y}^{18}-56609342891504\,{y}^{17}-141120525\,{y}^{38}\\
&&-2767783200\,{y}^{37}-26293940400\,{y}^{36}-162029846160\,{y}^{35}\\
&&-731482151400\,{y}^{34}-2588688538800\,{y}^{33}-92661976\,{y}^{4}\\
&&-14328872\,{y}^{3}-31589613164\,{y}^{8}-8905006466\,{y}^{7}-2237847235\,{y}^{6}\\
&&-492686736\,{y}^{5}-294908965224\,{y}^{10}-101142237120\,{y}^{9}\\
&&-1940993334032\,{y}^{12}-788267379728\,{y}^{11}-33353648804696\,{y}^{16}\\
&&-18289630975280\,{y}^{15}-9327903364424\,{y}^{14}-4419115551904\,{y}^{13}  \\
c_{24}&=&     -6036\,y-140952\,{y}^{2}-18307707855840\,{y}^{32}-38500788815216\,{y}^{31}\\
&&-70862893244864\,{y}^{30}-115577732666144\,{y}^{29}-168697385650976\,{y}^{28}\\
&&-222115579375536\,{y}^{27}-265542242629664\,{y}^{26}-289833757036928\,{y}^{25}\\
&&-290151785862120\,{y}^{24}-267452722477668\,{y}^{23}-227729992194848\,{y}^{22}\\
&&-179594286555632\,{y}^{21}-131450824320400\,{y}^{20}-89431228294824\,{y}^{19}\\
&&-56609342891504\,{y}^{18}-33353648804696\,{y}^{17}-141120525\,{y}^{39}\\
&&-2767783200\,{y}^{38}-26293940400\,{y}^{37}-162029846160\,{y}^{36}\\
&&-731482151400\,{y}^{35}-2588688538800\,{y}^{34}-7497840382760\,{y}^{33}\\
&&-14328872\,{y}^{4}-1711748\,{y}^{3}-8905006466\,{y}^{8}-2237847235\,{y}^{7}\\
&&-492686736\,{y}^{6}-92661976\,{y}^{5}-101142237120\,{y}^{10}-31589613164\,{y}^{9}\\
&&-788267379728\,{y}^{12}-294908965224\,{y}^{11}-18289630975280\,{y}^{16}\\
&&-9327903364424\,{y}^{15}-4419115551904\,{y}^{14}-1940993334032\,{y}^{13}  \\
c_{23}&=&   3434+78468\,y+939080\,{y}^{2}+34698266843604\,{y}^{32}+67244843649536\,{y}^{31}\\
&&+115581600821120\,{y}^{30}+177976221503808\,{y}^{29}+247527653786784\,{y}^{28}\\
&&+313044440817472\,{y}^{27}+362051571096048\,{y}^{26}+384774038456808\,{y}^{25}\\
&&+377305281107300\,{y}^{24}+342567709116672\,{y}^{23}+288828205184896\,{y}^{22}\\
&&+226682829211488\,{y}^{21}+165925005540528\,{y}^{20}+113432524543520\,{y}^{19}\\
&&+72495738192736\,{y}^{18}+43337018585352\,{y}^{17}+84672315\,{y}^{40}\\
&&+1748073600\,{y}^{39}+17480736000\,{y}^{38}+113364582240\,{y}^{37}\\
&&+538429273200\,{y}^{36}+2004180447360\,{y}^{35}+6105039384080\,{y}^{34}\\
&&+15679936536680\,{y}^{33}+50417880\,{y}^{4}+7806416\,{y}^{3}+17734396617\,{y}^{8}\\
&&+4933332096\,{y}^{7}+1227992256\,{y}^{6}+268745712\,{y}^{5}+172079552524\,{y}^{10}\\
&&+57768036230\,{y}^{9}+1197375768160\,{y}^{12}+471909608688\,{y}^{11}\\
&&+24233161211084\,{y}^{16}+12670853369728\,{y}^{15}+6189870335424\,{y}^{14}\\
&&+2821260753888\,{y}^{13}  \\
c_{22}&=& 3434\,y+78468\,{y}^{2}+67244843649536\,{y}^{32}+115581600821120\,{y}^{31}\\
&&+177976221503808\,{y}^{30}+247527653786784\,{y}^{29}+313044440817472\,{y}^{28}\\
&&+362051571096048\,{y}^{27}+384774038456808\,{y}^{26}+377305281107300\,{y}^{25}\\
&&+342567709116672\,{y}^{24}+288828205184896\,{y}^{23}+226682829211488\,{y}^{22}\\
&&+165925005540528\,{y}^{21}+113432524543520\,{y}^{20}+72495738192736\,{y}^{19}\\
&&+43337018585352\,{y}^{18}+24233161211084\,{y}^{17}+84672315\,{y}^{41}\\
&&+1748073600\,{y}^{40}+17480736000\,{y}^{39}+113364582240\,{y}^{38}\\
&&+538429273200\,{y}^{37}+2004180447360\,{y}^{36}+6105039384080\,{y}^{35}\\
&&+15679936536680\,{y}^{34}+34698266843604\,{y}^{33}+7806416\,{y}^{4}\\
&&+939080\,{y}^{3}+4933332096\,{y}^{8}+1227992256\,{y}^{7}+268745712\,{y}^{6}\\
&&+50417880\,{y}^{5}+57768036230\,{y}^{10}+17734396617\,{y}^{9}+471909608688\,{y}^{12}\\
&&+172079552524\,{y}^{11}+12670853369728\,{y}^{16}+6189870335424\,{y}^{15}\\
&&+2821260753888\,{y}^{14}+1197375768160\,{y}^{13} \\
c_{21}&=&    -1860-41208\,y-481390\,{y}^{2}-98205182723032\,{y}^{32}-158975487711040\,{y}^{31}\\
&&-232679061489952\,{y}^{30}-310041634001792\,{y}^{29}-378312376903472\,{y}^{28}\\
&&-424824424882648\,{y}^{27}-440914214605564\,{y}^{26}-424505680464736\,{y}^{25}\\
&&-380335644625872\,{y}^{24}-317952218824032\,{y}^{23}-248556355004464\,{y}^{22}\\
&&-182021476061568\,{y}^{21}-125037946100144\,{y}^{20}-80649214256600\,{y}^{19}\\
&&-48871493922372\,{y}^{18}-27829911367280\,{y}^{17}-44352165\,{y}^{42}\\
&&-961440480\,{y}^{41}-10095125040\,{y}^{40}-68726417760\,{y}^{39}-342554612400\,{y}^{38}\\
&&-1337717340960\,{y}^{37}-4274354404320\,{y}^{36}-11515763991160\,{y}^{35}\\
&&-26737946431196\,{y}^{34}-54393771958064\,{y}^{33}-25155248\,{y}^{4}-3933420\,{y}^{3}\\
&&-8913416664\,{y}^{8}-2461329648\,{y}^{7}-610337816\,{y}^{6}-133547584\,{y}^{5}\\
&&-88749319431\,{y}^{10}-29355926672\,{y}^{9}-643793038084\,{y}^{12}\\
&&-247995308890\,{y}^{11}-14890744366584\,{y}^{16}-7482590828128\,{y}^{15}\\
&&-3527920071680\,{y}^{14}-1558496855536\,{y}^{13}  \\
c_{20}&=&   -1860\,y-41208\,{y}^{2}-158975487711040\,{y}^{32}-232679061489952\,{y}^{31}\\
&&-310041634001792\,{y}^{30}-378312376903472\,{y}^{29}-424824424882648\,{y}^{28}\\
&&-440914214605564\,{y}^{27}-424505680464736\,{y}^{26}-380335644625872\,{y}^{25}\\
&&-317952218824032\,{y}^{24}-248556355004464\,{y}^{23}-182021476061568\,{y}^{22}\\
&&-125037946100144\,{y}^{21}-80649214256600\,{y}^{20}-48871493922372\,{y}^{19}\\
&&-27829911367280\,{y}^{18}-14890744366584\,{y}^{17}-44352165\,{y}^{43}\\
&&-961440480\,{y}^{42}-10095125040\,{y}^{41}-68726417760\,{y}^{40}-643793038084\,{y}^{13}\\
&&-342554612400\,{y}^{39}-1337717340960\,{y}^{38}-4274354404320\,{y}^{37}\\
&&-11515763991160\,{y}^{36}-26737946431196\,{y}^{35}-54393771958064\,{y}^{34}\\
&&-98205182723032\,{y}^{33}-3933420\,{y}^{4}-481390\,{y}^{3}-2461329648\,{y}^{8}\\
&&-610337816\,{y}^{7}-133547584\,{y}^{6}-25155248\,{y}^{5}-29355926672\,{y}^{10}\\
&&-8913416664\,{y}^{9}-247995308890\,{y}^{12}-88749319431\,{y}^{11}\\
&&-7482590828128\,{y}^{16}-3527920071680\,{y}^{15}-1558496855536\,{y}^{14}\\
c_{19}&=&   958+20460\,y+231152\,{y}^{2}+184072524213116\,{y}^{32}+257524494604608\,{y}^{31}\\
&&+330288464159280\,{y}^{30}+390331660078856\,{y}^{29}+426932172806836\,{y}^{28}\\
&&+433845313567488\,{y}^{27}+410966437011968\,{y}^{26}+363928613524112\,{y}^{25}\\
&&+302008912761192\,{y}^{24}+235340174568928\,{y}^{23}+172488228305888\,{y}^{22}\\
&&+119061638600344\,{y}^{21}+77473784139620\,{y}^{20}+47555285107616\,{y}^{19}\\
&&+27546629792112\,{y}^{18}+15059190513160\,{y}^{17}+20160075\,{y}^{44}\\
&&+457828800\,{y}^{43}+5036116800\,{y}^{42}+35910354480\,{y}^{41}+187411824600\,{y}^{40}\\
&&+766058092800\,{y}^{39}+2561529770800\,{y}^{38}+7221452461080\,{y}^{37}\\
&&+17547533041324\,{y}^{36}+37370557075456\,{y}^{35}+70667860599952\,{y}^{34}\\
&&+119898127453832\,{y}^{33}+11532806\,{y}^{4}+1839264\,{y}^{3}+4007479388\,{y}^{8}\\
&&+1104644976\,{y}^{7}+274478488\,{y}^{6}+60458908\,{y}^{5}+40475460160\,{y}^{10}\\
&&+13269333080\,{y}^{9}+302227968993\,{y}^{12}+114532596480\,{y}^{11}\\
&&+7767883706028\,{y}^{16}+3778587358192\,{y}^{15}+1731687401452\,{y}^{14}\\
&&+746622879894\,{y}^{13}  \\
c_{18}&=&   958\,y+20460\,{y}^{2}+257524494604608\,{y}^{32}+330288464159280\,{y}^{31}\\
&&+390331660078856\,{y}^{30}+426932172806836\,{y}^{29}+433845313567488\,{y}^{28}\\
&&+410966437011968\,{y}^{27}+363928613524112\,{y}^{26}+302008912761192\,{y}^{25}\\
&&+235340174568928\,{y}^{24}+172488228305888\,{y}^{23}+119061638600344\,{y}^{22}\\
&&+77473784139620\,{y}^{21}+47555285107616\,{y}^{20}+27546629792112\,{y}^{19}\\
&&+15059190513160\,{y}^{18}+7767883706028\,{y}^{17}+20160075\,{y}^{45}\\
&&+457828800\,{y}^{44}+5036116800\,{y}^{43}+35910354480\,{y}^{42}+302227968993\,{y}^{13}\\
&&+187411824600\,{y}^{41}+766058092800\,{y}^{40}+2561529770800\,{y}^{39}\\
&&+7221452461080\,{y}^{38}+17547533041324\,{y}^{37}+37370557075456\,{y}^{36}\\
&&+70667860599952\,{y}^{35}+119898127453832\,{y}^{34}+184072524213116\,{y}^{33}\\
&&+1839264\,{y}^{4}+231152\,{y}^{3}+1104644976\,{y}^{8}+274478488\,{y}^{7}\\
&&+60458908\,{y}^{6}+11532806\,{y}^{5}+13269333080\,{y}^{10}+4007479388\,{y}^{9}\\
&&+114532596480\,{y}^{12}+40475460160\,{y}^{11}+3778587358192\,{y}^{16}\\
&&+1731687401452\,{y}^{15}+746622879894\,{y}^{14} \\
c_{17}&=&  -470-9580\,y-104030\,{y}^{2}-239995902220932\,{y}^{32}-297452787946808\,{y}^{31}\\
&&-341592891568364\,{y}^{30}-364919215746912\,{y}^{29}-363903464593872\,{y}^{28}\\
&&-339772161667024\,{y}^{27}-297806217436296\,{y}^{26}-245578506693520\,{y}^{25}\\
&&-190884200458392\,{y}^{24}-140069671803336\,{y}^{23}-97153052600812\,{y}^{22}\\
&&-63758016163248\,{y}^{21}-39618645403896\,{y}^{20}-23322138259880\,{y}^{19}\\
&&-13009335747612\,{y}^{18}-6876470437704\,{y}^{17}-7888725\,{y}^{46}-187293600\,{y}^{45}\\
&&-2153876400\,{y}^{44}-16053214320\,{y}^{43}-87541997400\,{y}^{42}-373773800400\,{y}^{41}\\
&&-1305134270440\,{y}^{40}-3841741051720\,{y}^{39}-9747162009156\,{y}^{38}\\
&&-21678813590544\,{y}^{37}-42828036093656\,{y}^{36}-75952348723368\,{y}^{35}\\
&&-121961553432684\,{y}^{34}-178605863828952\,{y}^{33}-4865820\,{y}^{4}-798700\,{y}^{3}\\
&&-1607047332\,{y}^{8}-444849428\,{y}^{7}-111553074\,{y}^{6}-24949832\,{y}^{5}\\
&&-16284818508\,{y}^{10}-5320278232\,{y}^{9}-123743046968\,{y}^{12}-46403526648\,{y}^{11}\\
&&-3443385293384\,{y}^{16}-1632526042378\,{y}^{15}-732079699215\,{y}^{14}\\
&&-310045351120\,{y}^{13}\\
c_{16}&=&    -470\,y-9580\,{y}^{2}-297452787946808\,{y}^{32}-341592891568364\,{y}^{31}\\
&&-364919215746912\,{y}^{30}-363903464593872\,{y}^{29}-339772161667024\,{y}^{28}\\
&&-297806217436296\,{y}^{27}-245578506693520\,{y}^{26}-190884200458392\,{y}^{25}\\
&&-140069671803336\,{y}^{24}-97153052600812\,{y}^{23}-63758016163248\,{y}^{22}\\
&&-39618645403896\,{y}^{21}-23322138259880\,{y}^{20}-13009335747612\,{y}^{19}\\
&&-6876470437704\,{y}^{18}-3443385293384\,{y}^{17}-7888725\,{y}^{47}-187293600\,{y}^{46}\\
&&-2153876400\,{y}^{45}-16053214320\,{y}^{44}-87541997400\,{y}^{43}-373773800400\,{y}^{42}\\
&&-1305134270440\,{y}^{41}-3841741051720\,{y}^{40}-9747162009156\,{y}^{39}\\
&&-21678813590544\,{y}^{38}-42828036093656\,{y}^{37}-75952348723368\,{y}^{36}\\
&&-121961553432684\,{y}^{35}-178605863828952\,{y}^{34}-239995902220932\,{y}^{33}\\
&&-798700\,{y}^{4}-104030\,{y}^{3}-444849428\,{y}^{8}-111553074\,{y}^{7}-24949832\,{y}^{6}\\
&&-4865820\,{y}^{5}-5320278232\,{y}^{10}-1607047332\,{y}^{9}-46403526648\,{y}^{12}\\
&&-16284818508\,{y}^{11}-1632526042378\,{y}^{16}-732079699215\,{y}^{15}\\
&&-310045351120\,{y}^{14}-123743046968\,{y}^{13}  \\
c_{15}&=&   221+4230\,y+43696\,{y}^{2}+224453875410934\,{y}^{32}+251204504233472\,{y}^{31}\\
&&+262741215297792\,{y}^{30}+257628708088256\,{y}^{29}+237475576589280\,{y}^{28}\\
&&+206270655047872\,{y}^{27}+169176128332608\,{y}^{26}+131244284158128\,{y}^{25}\\
&&+96449234110088\,{y}^{24}+67224175268736\,{y}^{23}+44483237082816\,{y}^{22}\\
&&+27967868049376\,{y}^{21}+16717751610432\,{y}^{20}+9504532480656\,{y}^{19}\\
&&+5140527848084\,{y}^{18}+2644832316938\,{y}^{17}+2629575\,{y}^{48}+65145600\,{y}^{47}\\
&&+781747200\,{y}^{46}+6078646080\,{y}^{45}+34571815200\,{y}^{44}+153893625600\,{y}^{43}\\
&&+560065128480\,{y}^{42}+1717893903440\,{y}^{41}+4541584786024\,{y}^{40}\\
&&+10526132702976\,{y}^{39}+21675724202176\,{y}^{38}+40083624525984\,{y}^{37}\\
&&+67150920705424\,{y}^{36}+102659859883488\,{y}^{35}+144113043662568\,{y}^{34}\\
&&+186755849437596\,{y}^{33}+1886008\,{y}^{4}+321512\,{y}^{3}+572231076\,{y}^{8}\\
&&+160253440\,{y}^{7}+40907456\,{y}^{6}+9375440\,{y}^{5}+5749591888\,{y}^{10}\\
&&+1882331432\,{y}^{9}+43987583888\,{y}^{12}+16410546464\,{y}^{11}+1294108820623\,{y}^{16}\\
&&+601787533568\,{y}^{15}+265670307712\,{y}^{14}+111164450976\,{y}^{13}\\
c_{14}&=&    221\,y+4230\,{y}^{2}+251204504233472\,{y}^{32}+262741215297792\,{y}^{31}\\
&&+257628708088256\,{y}^{30}+237475576589280\,{y}^{29}+206270655047872\,{y}^{28}\\
&&+169176128332608\,{y}^{27}+131244284158128\,{y}^{26}+96449234110088\,{y}^{25}\\
&&+67224175268736\,{y}^{24}+44483237082816\,{y}^{23}+27967868049376\,{y}^{22}\\
&&+16717751610432\,{y}^{21}+9504532480656\,{y}^{20}+5140527848084\,{y}^{19}\\
&&+2644832316938\,{y}^{18}+1294108820623\,{y}^{17}+2629575\,{y}^{49}+65145600\,{y}^{48}\\
&&+781747200\,{y}^{47}+6078646080\,{y}^{46}+34571815200\,{y}^{45}+153893625600\,{y}^{44}\\
&&+560065128480\,{y}^{43}+1717893903440\,{y}^{42}+4541584786024\,{y}^{41}\\
&&+10526132702976\,{y}^{40}+21675724202176\,{y}^{39}+40083624525984\,{y}^{38}\\
&&+67150920705424\,{y}^{37}+102659859883488\,{y}^{36}+144113043662568\,{y}^{35}\\
&&+186755849437596\,{y}^{34}+224453875410934\,{y}^{33}+321512\,{y}^{4}+43696\,{y}^{3}\\
&&+160253440\,{y}^{8}+40907456\,{y}^{7}+9375440\,{y}^{6}+1886008\,{y}^{5}\\
&&+1882331432\,{y}^{10}+572231076\,{y}^{9}+16410546464\,{y}^{12}+5749591888\,{y}^{11}\\
&&+601787533568\,{y}^{16}+265670307712\,{y}^{15}+111164450976\,{y}^{14}\\
&&+43987583888\,{y}^{13}  \\
c_{13}&=&   -100-1768\,y-17083\,{y}^{2}-153233296391224\,{y}^{32}-157250659653568\,{y}^{31}\\
&&-151878629674464\,{y}^{30}-138406407201664\,{y}^{29}-119267409857632\,{y}^{28}\\
&&-97368253050768\,{y}^{27}-75431670307928\,{y}^{26}-55531823760672\,{y}^{25}\\
&&-38896049239568\,{y}^{24}-25947141564512\,{y}^{23}-16499352006432\,{y}^{22}\\
&&-10007856635088\,{y}^{21}-5793562425276\,{y}^{20}-3202132928406\,{y}^{19}\\
&&-1690028073569\,{y}^{18}-851676813040\,{y}^{17}-736281\,{y}^{50}-19000800\,{y}^{49}\\
&&-237510000\,{y}^{48}-1923405120\,{y}^{47}-11389341600\,{y}^{46}-52765702080\,{y}^{45}\\
&&-199792087040\,{y}^{44}-637437406320\,{y}^{43}-1752659764856\,{y}^{42}\\
&&-4224896000864\,{y}^{41}-9049932306800\,{y}^{40}-17413503368352\,{y}^{39}\\
&&-30366290882256\,{y}^{38}-48348103454336\,{y}^{37}-70726238869416\,{y}^{36}\\
&&-95575763081092\,{y}^{35}-119875889776986\,{y}^{34}-140131546291984\,{y}^{33}\\
&&-668072\,{y}^{4}-119178\,{y}^{3}-179857704\,{y}^{8}-51395408\,{y}^{7}-13472664\,{y}^{6}\\
&&-3191560\,{y}^{5}-1766676236\,{y}^{10}-583385904\,{y}^{9}-13459829376\,{y}^{12}\\
&&-5022155192\,{y}^{11}-409652090600\,{y}^{16}-187926274976\,{y}^{15}\\
&&-82125226896\,{y}^{14}-34131247424\,{y}^{13}\\
c_{12}&=&   -100\,y-1768\,{y}^{2}-157250659653568\,{y}^{32}-151878629674464\,{y}^{31}\\
&&-138406407201664\,{y}^{30}-119267409857632\,{y}^{29}-97368253050768\,{y}^{28}\\
&&-75431670307928\,{y}^{27}-55531823760672\,{y}^{26}-38896049239568\,{y}^{25}\\
&&-25947141564512\,{y}^{24}-16499352006432\,{y}^{23}-10007856635088\,{y}^{22}\\
&&-5793562425276\,{y}^{21}-3202132928406\,{y}^{20}-1690028073569\,{y}^{19}\\
&&-851676813040\,{y}^{18}-409652090600\,{y}^{17}-736281\,{y}^{51}-19000800\,{y}^{50}\\
&&-237510000\,{y}^{49}-1923405120\,{y}^{48}-11389341600\,{y}^{47}-13459829376\,{y}^{13}\\
&&-52765702080\,{y}^{46}-199792087040\,{y}^{45}-637437406320\,{y}^{44}\\
&&-1752659764856\,{y}^{43}-4224896000864\,{y}^{42}-9049932306800\,{y}^{41}\\
&&-17413503368352\,{y}^{40}-30366290882256\,{y}^{39}-48348103454336\,{y}^{38}\\
&&-70726238869416\,{y}^{37}-95575763081092\,{y}^{36}-119875889776986\,{y}^{35}\\
&&-140131546291984\,{y}^{34}-153233296391224\,{y}^{33}-119178\,{y}^{4}-17083\,{y}^{3}\\
&&-51395408\,{y}^{8}-13472664\,{y}^{7}-3191560\,{y}^{6}-668072\,{y}^{5}\\
&&-583385904\,{y}^{10}-179857704\,{y}^{9}-5022155192\,{y}^{12}-1766676236\,{y}^{11}\\
&&-187926274976\,{y}^{16}-82125226896\,{y}^{15}-34131247424\,{y}^{14}  \\
c_{11}&=&  42+700\,y+6216\,{y}^{2}+76808425621716\,{y}^{32}+73185606168640\,{y}^{31}\\
&&+66018061921600\,{y}^{30}+56493484492368\,{y}^{29}+45941140332152\,{y}^{28}\\
&&+35558795229248\,{y}^{27}+26231648793312\,{y}^{26}+18465315874512\,{y}^{25}\\
&&+12416392346504\,{y}^{24}+7982466313872\,{y}^{23}+4910432839540\,{y}^{22}\\
&&+2892162420314\,{y}^{21}+1631773224327\,{y}^{20}+882200588160\,{y}^{19}\\
&&+457073117824\,{y}^{18}+226902468024\,{y}^{17}+169911\,{y}^{52}+4560192\,{y}^{51}\\
&&+59282496\,{y}^{50}+499196880\,{y}^{49}+3072724200\,{y}^{48}+14792494080\,{y}^{47}\\
&&+58180609760\,{y}^{46}+192763266800\,{y}^{45}+550294423224\,{y}^{44}\\
&&+1377218757376\,{y}^{43}+3063057109728\,{y}^{42}+6120734583600\,{y}^{41}\\
&&+11087868148168\,{y}^{40}+18346259909600\,{y}^{39}+27904269575304\,{y}^{38}\\
&&+39228726807340\,{y}^{37}+51219029799966\,{y}^{36}+62371671684096\,{y}^{35}\\
&&+71104825778240\,{y}^{34}+76139818059720\,{y}^{33}+214017\,{y}^{4}+40400\,{y}^{3}\\
&&+49420548\,{y}^{8}+14531784\,{y}^{7}+3945168\,{y}^{6}+974334\,{y}^{5}+8853161400\,{y}^{13}\\
&&+466960416\,{y}^{10}+156747336\,{y}^{9}+3498053324\,{y}^{12}+1313106496\,{y}^{11}\\
&&+107874349260\,{y}^{16}+49075383776\,{y}^{15}+21338054064\,{y}^{14}  \\
c_{10}&=&   42\,y+700\,{y}^{2}+73185606168640\,{y}^{32}+66018061921600\,{y}^{31}\\
&&+56493484492368\,{y}^{30}+45941140332152\,{y}^{29}+35558795229248\,{y}^{28}\\
&&+26231648793312\,{y}^{27}+18465315874512\,{y}^{26}+12416392346504\,{y}^{25}\\
&&+7982466313872\,{y}^{24}+4910432839540\,{y}^{23}+2892162420314\,{y}^{22}\\
&&+1631773224327\,{y}^{21}+882200588160\,{y}^{20}+457073117824\,{y}^{19}\\
&&+226902468024\,{y}^{18}+107874349260\,{y}^{17}+169911\,{y}^{53}+3498053324\,{y}^{13}\\
&&+4560192\,{y}^{52}+59282496\,{y}^{51}+499196880\,{y}^{50}+3072724200\,{y}^{49}\\
&&+14792494080\,{y}^{48}+58180609760\,{y}^{47}+192763266800\,{y}^{46}\\
&&+550294423224\,{y}^{45}+1377218757376\,{y}^{44}+3063057109728\,{y}^{43}\\
&&+6120734583600\,{y}^{42}+11087868148168\,{y}^{41}+18346259909600\,{y}^{40}\\
&&+27904269575304\,{y}^{39}+39228726807340\,{y}^{38}+51219029799966\,{y}^{37}\\
&&+62371671684096\,{y}^{36}+71104825778240\,{y}^{35}+76139818059720\,{y}^{34}\\
&&+76808425621716\,{y}^{33}+40400\,{y}^{4}+6216\,{y}^{3}+14531784\,{y}^{8}\\
&&+3945168\,{y}^{7}+974334\,{y}^{6}+214017\,{y}^{5}+156747336\,{y}^{10}+8853161400\,{y}^{14}\\
&&+49420548\,{y}^{9}+1313106496\,{y}^{12}+466960416\,{y}^{11}+49075383776\,{y}^{16}\\
&&+21338054064\,{y}^{15} \\
c_{9}&=&    -14-252\,y-2090\,{y}^{2}-28077372499868\,{y}^{32}-25098955113968\,{y}^{31}\\
&&-21347223917672\,{y}^{30}-17303673015328\,{y}^{29}-13387253668816\,{y}^{28}\\
&&-9898697135120\,{y}^{27}-7003509266248\,{y}^{26}-4746453967584\,{y}^{25}\\
&&-3084283289652\,{y}^{24}-1923259624582\,{y}^{23}-1151691937865\,{y}^{22}\\
&&-662684126384\,{y}^{21}-366550654856\,{y}^{20}-194947311000\,{y}^{19}\\
&&-99690728668\,{y}^{18}-49004477896\,{y}^{17}-31465\,{y}^{54}-876960\,{y}^{53}\\
&&-11838960\,{y}^{52}-103508496\,{y}^{51}-661326120\,{y}^{50}-3303409200\,{y}^{49}\\
&&-13476208200\,{y}^{48}-46296237520\,{y}^{47}-137010274856\,{y}^{46}\\
&&-355427925664\,{y}^{45}-819405519408\,{y}^{44}-1697441917808\,{y}^{43}\\
&&-3188473544424\,{y}^{42}-5472176918320\,{y}^{41}-8636341292464\,{y}^{40}\\
&&-12604032196532\,{y}^{39}-17092716269202\,{y}^{38}-21632053806736\,{y}^{37}\\
&&-25646306254712\,{y}^{36}-28580395790376\,{y}^{35}-30029268395908\,{y}^{34}\\
&&-29828198869240\,{y}^{33}-61164\,{y}^{4}-12444\,{y}^{3}-11665092\,{y}^{8}\\
&&-3559218\,{y}^{7}-1009351\,{y}^{6}-262264\,{y}^{5}-104312228\,{y}^{10}-1900252432\,{y}^{13} \\
&&-35884960\,{y}^{9}-756947512\,{y}^{12}-287897544\,{y}^{11}-23143828508\,{y}^{16}\\
&&-10493317928\,{y}^{15}-4562670436\,{y}^{14}  \\
c_{8}&=&    -14\,y-252\,{y}^{2}-25098955113968\,{y}^{32}-21347223917672\,{y}^{31}\\
&&-17303673015328\,{y}^{30}-13387253668816\,{y}^{29}-9898697135120\,{y}^{28}\\
&&-7003509266248\,{y}^{27}-4746453967584\,{y}^{26}-3084283289652\,{y}^{25}\\
&&-1923259624582\,{y}^{24}-1151691937865\,{y}^{23}-662684126384\,{y}^{22}\\
&&-366550654856\,{y}^{21}-194947311000\,{y}^{20}-99690728668\,{y}^{19}\\
&&-49004477896\,{y}^{18}-23143828508\,{y}^{17}-31465\,{y}^{55}-876960\,{y}^{54}\\
&&-11838960\,{y}^{53}-103508496\,{y}^{52}-661326120\,{y}^{51}-3303409200\,{y}^{50}\\
&&-13476208200\,{y}^{49}-46296237520\,{y}^{48}-137010274856\,{y}^{47}\\
&&-355427925664\,{y}^{46}-819405519408\,{y}^{45}-1697441917808\,{y}^{44}\\
&&-3188473544424\,{y}^{43}-5472176918320\,{y}^{42}-8636341292464\,{y}^{41}\\
&&-12604032196532\,{y}^{40}-17092716269202\,{y}^{39}-21632053806736\,{y}^{38}\\
&&-25646306254712\,{y}^{37}-28580395790376\,{y}^{36}-30029268395908\,{y}^{35}\\
&&-29828198869240\,{y}^{34}-28077372499868\,{y}^{33}-12444\,{y}^{4}-2090\,{y}^{3}\\
&&-3559218\,{y}^{8}-1009351\,{y}^{7}-262264\,{y}^{6}-61164\,{y}^{5}-35884960\,{y}^{10}\\
&&-11665092\,{y}^{9}-287897544\,{y}^{12}-104312228\,{y}^{11}-10493317928\,{y}^{16}\\
&&-4562670436\,{y}^{15}-1900252432\,{y}^{14}-756947512\,{y}^{13}   \\
c_{7}&=&   5+70\,y+616\,{y}^{2}+7313081337386\,{y}^{32}+6187019672704\,{y}^{31}\\
&&+5001900819008\,{y}^{30}+3869759209760\,{y}^{29}+2868764374752\,{y}^{28}\\
&&+2040261727824\,{y}^{27}+1393577428740\,{y}^{26}+915089088130\,{y}^{25}\\
&&+578194161219\,{y}^{24}+351808762496\,{y}^{23}+206277316544\,{y}^{22}\\
&&+116608603248\,{y}^{21}+63575384600\,{y}^{20}+33433959952\,{y}^{19}+1776429824\,{y}^{15}\\
&&+16959111912\,{y}^{18}+8295244308\,{y}^{17}+4495\,{y}^{56}+129920\,{y}^{55}\\
&&+1818880\,{y}^{54}+16488864\,{y}^{53}+109201680\,{y}^{52}+565219200\,{y}^{51}\\
&&+2388353968\,{y}^{50}+8495842680\,{y}^{49}+26027492764\,{y}^{48}+326893136\,{y}^{13} \\
&&+69884743936\,{y}^{47}+166748284864\,{y}^{46}+357531140512\,{y}^{45}+51357240\,{y}^{11}\\
&&+695222438736\,{y}^{44}+1235449012000\,{y}^{43}+2019536911384\,{y}^{42}\\
&&+3053870732836\,{y}^{41}+4292981378154\,{y}^{40}+5634619819904\,{y}^{39}\\
&&+6931822075968\,{y}^{38}+8020643537936\,{y}^{37}+8755782372104\,{y}^{36}\\
&&+9042912134000\,{y}^{35}+8857775392592\,{y}^{34}+8247140069852\,{y}^{33}\\
&&+15316\,{y}^{4}+3408\,{y}^{3}+2295159\,{y}^{8}+733024\,{y}^{7}+219248\,{y}^{6}\\
&&+60696\,{y}^{5}+19109072\,{y}^{10}+6790994\,{y}^{9}+132235736\,{y}^{12}+776871360\,{y}^{14}\\
&&+3911002194\,{y}^{16}\\
c_{6}&=&   5\,y+70\,{y}^{2}+6187019672704\,{y}^{32}+5001900819008\,{y}^{31}+3869759209760\,{y}^{30}\\
&&+2868764374752\,{y}^{29}+2040261727824\,{y}^{28}+1393577428740\,{y}^{27}\\
&&+578194161219\,{y}^{25}+351808762496\,{y}^{24}+206277316544\,{y}^{23}\\
&&+63575384600\,{y}^{21}+33433959952\,{y}^{20}+16959111912\,{y}^{19}+8295244308\,{y}^{18}\\
&&+3911002194\,{y}^{17}+4495\,{y}^{57}+129920\,{y}^{56}+1818880\,{y}^{55}+16488864\,{y}^{54}\\
&&+109201680\,{y}^{53}+565219200\,{y}^{52}+2388353968\,{y}^{51}+8495842680\,{y}^{50}\\
&&+26027492764\,{y}^{49}+69884743936\,{y}^{48}+166748284864\,{y}^{47}\\
&&+695222438736\,{y}^{45}+1235449012000\,{y}^{44}+2019536911384\,{y}^{43}\\
&&+4292981378154\,{y}^{41}+5634619819904\,{y}^{40}+6931822075968\,{y}^{39}\\
&&+8755782372104\,{y}^{37}+9042912134000\,{y}^{36}+8857775392592\,{y}^{35}\\
&&+7313081337386\,{y}^{33}+3408\,{y}^{4}+616\,{y}^{3}+733024\,{y}^{8}+219248\,{y}^{7}\\
&&+60696\,{y}^{6}+15316\,{y}^{5}+6790994\,{y}^{10}+2295159\,{y}^{9}+51357240\,{y}^{12}\\
&&+19109072\,{y}^{11}+1776429824\,{y}^{16}+776871360\,{y}^{15}+326893136\,{y}^{14}\\
&&+132235736\,{y}^{13}+915089088130\,{y}^{26}+116608603248\,{y}^{22}\\
&&+3053870732836\,{y}^{42}+8020643537936\,{y}^{38}+8247140069852\,{y}^{34}   \\
&&+357531140512\,{y}^{46}\\
c_{5}&=&   -2-20\,y-135\,{y}^{2}-1291240176900\,{y}^{32}-1041748018848\,{y}^{31}\\
&&-806281172416\,{y}^{30}-599429462736\,{y}^{29}-428584659180\,{y}^{28}\\
&&-295028984862\,{y}^{27}-195735333933\,{y}^{26}-125273542192\,{y}^{25}\\
&&-77410559816\,{y}^{24}-46217249904\,{y}^{23}-26676063960\,{y}^{22}-491718996\,{y}^{16}\\
&&-14891452032\,{y}^{21}-8042059472\,{y}^{20}-4202120092\,{y}^{19}-99677176\,{y}^{14}\\
&&-2124447526\,{y}^{18}-1039137368\,{y}^{17}-465\,{y}^{58}-13920\,{y}^{57}-7062174\,{y}^{11}\\
&&-201840\,{y}^{56}-1894816\,{y}^{55}-12991440\,{y}^{54}-69588960\,{y}^{53}-17671720\,{y}^{12}\\
&&-304194016\,{y}^{52}-1119005992\,{y}^{51}-3544114548\,{y}^{50}-42689912\,{y}^{13}\\
&&-9836099344\,{y}^{49}-24256384264\,{y}^{48}-53753140896\,{y}^{47}-225093584\,{y}^{15}\\
&&-108038522256\,{y}^{46}-198480699200\,{y}^{45}-335496053496\,{y}^{44}-1007064\,{y}^{9}\\
&&-524756939068\,{y}^{43}-763293723014\,{y}^{42}-1037044904688\,{y}^{41}-2719985\,{y}^{10}\\
&&-1321226300232\,{y}^{40}-1583998322704\,{y}^{39}-1792664023752\,{y}^{38}\\
&&-1920611818240\,{y}^{37}-1952903458120\,{y}^{36}-1888915844532\,{y}^{35}\\
&&-1741500968654\,{y}^{34}-1533253370056\,{y}^{33}-3204\,{y}^{4}-770\,{y}^{3}\\
&&-357052\,{y}^{8}-120568\,{y}^{7}-38564\,{y}^{6}-11592\,{y}^{5}\\
c_{4}&=&  -2\,y-20\,{y}^{2}-1041748018848\,{y}^{32}-806281172416\,{y}^{31}-599429462736\,{y}^{30}\\
&&-428584659180\,{y}^{29}-295028984862\,{y}^{28}-195735333933\,{y}^{27}\\
&&-77410559816\,{y}^{25}-46217249904\,{y}^{24}-26676063960\,{y}^{23}-14891452032\,{y}^{22}\\
&&-8042059472\,{y}^{21}-4202120092\,{y}^{20}-2124447526\,{y}^{19}-1039137368\,{y}^{18}\\
&&-491718996\,{y}^{17}-465\,{y}^{59}-13920\,{y}^{58}-201840\,{y}^{57}-1894816\,{y}^{56}\\
&&-12991440\,{y}^{55}-69588960\,{y}^{54}-304194016\,{y}^{53}-1119005992\,{y}^{52}\\
&&-9836099344\,{y}^{50}-24256384264\,{y}^{49}-53753140896\,{y}^{48}-108038522256\,{y}^{47}\\
&&-198480699200\,{y}^{46}-335496053496\,{y}^{45}-524756939068\,{y}^{44}\\
&&-1037044904688\,{y}^{42}-1321226300232\,{y}^{41}-1583998322704\,{y}^{40}\\
&&-1920611818240\,{y}^{38}-1952903458120\,{y}^{37}-1888915844532\,{y}^{36}\\
&&-1533253370056\,{y}^{34}-1291240176900\,{y}^{33}-770\,{y}^{4}-135\,{y}^{3}-120568\,{y}^{8}\\
&&-38564\,{y}^{7}-11592\,{y}^{6}-3204\,{y}^{5}-1007064\,{y}^{10}-357052\,{y}^{9}-7062174\,{y}^{12}\\
&&-2719985\,{y}^{11}-225093584\,{y}^{16}-99677176\,{y}^{15}-42689912\,{y}^{14}\\
&&-125273542192\,{y}^{26}-3544114548\,{y}^{51}-763293723014\,{y}^{43}-17671720\,{y}^{13}\\
&&-1741500968654\,{y}^{35}-1792664023752\,{y}^{39}   \\
c_{3}&=&   1+6\,y+28\,{y}^{2}+139134657962\,{y}^{32}+107772280144\,{y}^{31}+80374075940\,{y}^{30}\\
&&+57780713586\,{y}^{29}+40086413227\,{y}^{28}+26866822144\,{y}^{27}+17412522304\,{y}^{26}\\
&&+10922249064\,{y}^{25}+6635841796\,{y}^{24}+3907393968\,{y}^{23}+2231027256\,{y}^{22}\\
&&+1235719404\,{y}^{21}+664148942\,{y}^{20}+346461216\,{y}^{19}+175470704\,{y}^{18}\\
&&+86311084\,{y}^{17}+31\,{y}^{60}+960\,{y}^{59}+14400\,{y}^{58}+139824\,{y}^{57}+1650\,{y}^{5}\\
&&+991320\,{y}^{56}+5488896\,{y}^{55}+24792208\,{y}^{54}+94201160\,{y}^{53}+1622763\,{y}^{12}\\
&&+308074372\,{y}^{52}+882664192\,{y}^{51}+2246800496\,{y}^{50}+5139121208\,{y}^{49}\\
&&+10661724164\,{y}^{48}+20219969056\,{y}^{47}+35289119800\,{y}^{46}+3802026\,{y}^{13}\\
&&+57003970868\,{y}^{45}+85655656018\,{y}^{44}+120260511360\,{y}^{43}+41250862\,{y}^{16}\\
&&+158389551488\,{y}^{42}+196385714120\,{y}^{41}+229965381556\,{y}^{40}+672720\,{y}^{11}\\
&&+255056476176\,{y}^{39}+268633279376\,{y}^{38}+269307919028\,{y}^{37}+14352\,{y}^{7}\\
&&+257525633614\,{y}^{36}+235343516912\,{y}^{35}+205895850024\,{y}^{34}+5000\,{y}^{6}\\
&&+172721468332\,{y}^{33}+495\,{y}^{4}+120\,{y}^{3}+39586\,{y}^{8}+8657072\,{y}^{14} \\
&&+270536\,{y}^{10}+105324\,{y}^{9}+19164440\,{y}^{15}\\
c_{2}&=&  y+6\,{y}^{2}+107772280144\,{y}^{32}+80374075940\,{y}^{31}+57780713586\,{y}^{30}\\
&&+40086413227\,{y}^{29}+26866822144\,{y}^{28}+17412522304\,{y}^{27}+10922249064\,{y}^{26}\\
&&+6635841796\,{y}^{25}+3907393968\,{y}^{24}+2231027256\,{y}^{23}+1235719404\,{y}^{22}\\
&&+664148942\,{y}^{21}+346461216\,{y}^{20}+175470704\,{y}^{19}+86311084\,{y}^{18}\\
&&+41250862\,{y}^{17}+31\,{y}^{61}+960\,{y}^{60}+14400\,{y}^{59}+139824\,{y}^{58}\\
&&+991320\,{y}^{57}+5488896\,{y}^{56}+24792208\,{y}^{55}+94201160\,{y}^{54}\\
&&+308074372\,{y}^{53}+882664192\,{y}^{52}+2246800496\,{y}^{51}+5139121208\,{y}^{50}\\
&&+10661724164\,{y}^{49}+20219969056\,{y}^{48}+35289119800\,{y}^{47}+57003970868\,{y}^{46}\\
&&+85655656018\,{y}^{45}+120260511360\,{y}^{44}+158389551488\,{y}^{43}\\
&&+229965381556\,{y}^{41}+255056476176\,{y}^{40}+268633279376\,{y}^{39}\\
&&+257525633614\,{y}^{37}+235343516912\,{y}^{36}+205895850024\,{y}^{35}\\
&&+139134657962\,{y}^{33}+120\,{y}^{4}+28\,{y}^{3}+14352\,{y}^{8}+5000\,{y}^{7}+1650\,{y}^{6}\\
&&+495\,{y}^{5}+105324\,{y}^{10}+39586\,{y}^{9}+672720\,{y}^{12}+270536\,{y}^{11}\\
&&+19164440\,{y}^{16}+8657072\,{y}^{15}+3802026\,{y}^{14}+1622763\,{y}^{13}  \\
&&+196385714120\,{y}^{42}+269307919028\,{y}^{38}+172721468332\,{y}^{34}\\
c_{1}&=&   -1-2\,y-5\,{y}^{2}-6939692682\,{y}^{32}-5193067630\,{y}^{31}-3754272037\,{y}^{30}\\
&&-2625062128\,{y}^{29}-1777171560\,{y}^{28}-1166067016\,{y}^{27}-742179284\,{y}^{26}\\
&&-458591432\,{y}^{25}-275276716\,{y}^{24}-160617860\,{y}^{23}-91143114\,{y}^{22}-132\,{y}^{5}\\
&&-50323496\,{y}^{21}-27049196\,{y}^{20}-14162220\,{y}^{19}-7228014\,{y}^{18}-5916\,{y}^{9}\\
&&-3598964\,{y}^{17}-{y}^{62}-32\,{y}^{61}-496\,{y}^{60}-4976\,{y}^{59}-175048\,{y}^{13}-14290\,{y}^{10}\\
&&-36440\,{y}^{58}-208336\,{y}^{57}-971272\,{y}^{56}-3807704\,{y}^{55}-385741\,{y}^{14}\\
&&-12843980\,{y}^{54}-37945904\,{y}^{53}-99582920\,{y}^{52}-234813592\,{y}^{51}-365\,{y}^{6}\\
&&-502196500\,{y}^{50}-981900168\,{y}^{49}-1766948340\,{y}^{48}-2943492972\,{y}^{47}\\
&&-4562339774\,{y}^{46}-6609143792\,{y}^{45}-8984070856\,{y}^{44}-831014\,{y}^{15} \\
&&-11500901864\,{y}^{43}-13910043524\,{y}^{42}-15941684776\,{y}^{41}-77684\,{y}^{12}\\
&&-17357937708\,{y}^{40}-17999433372\,{y}^{39}-17813777994\,{y}^{38}-33708\,{y}^{11}\\
&&-16859410792\,{y}^{37}-15286065700\,{y}^{36}-13299362332\,{y}^{35}-1749654\,{y}^{16}\\
&&-11120136162\,{y}^{34}-8948546308\,{y}^{33}-42\,{y}^{4}-14\,{y}^{3}-2398\,{y}^{8}-950\,{y}^{7}\\
c_{0}&=&-1-y-2\,{y}^{2}-5193067630\,{y}^{32}-3754272037\,{y}^{31}-2625062128\,{y}^{30}\\
&&-1777171560\,{y}^{29}-1166067016\,{y}^{28}-742179284\,{y}^{27}-458591432\,{y}^{26}\\
&&-275276716\,{y}^{25}-160617860\,{y}^{24}-91143114\,{y}^{23}-50323496\,{y}^{22}\\
&&-27049196\,{y}^{21}-14162220\,{y}^{20}-7228014\,{y}^{19}-3598964\,{y}^{18}-1749654\,{y}^{17}\\
&&-{y}^{63}-32\,{y}^{62}-496\,{y}^{61}-4976\,{y}^{60}-36440\,{y}^{59}-208336\,{y}^{58}-971272\,{y}^{57}\\
&&-3807704\,{y}^{56}-12843980\,{y}^{55}-37945904\,{y}^{54}-99582920\,{y}^{53}-175048\,{y}^{14}\\
&&-502196500\,{y}^{51}-981900168\,{y}^{50}-1766948340\,{y}^{49}-2943492972\,{y}^{48}\\
&&-6609143792\,{y}^{46}-8984070856\,{y}^{45}-11500901864\,{y}^{44}-13910043524\,{y}^{43}\\
&&-15941684776\,{y}^{42}-17357937708\,{y}^{41}-17999433372\,{y}^{40}-17813777994\,{y}^{39}\\
&&-16859410792\,{y}^{38}-15286065700\,{y}^{37}-13299362332\,{y}^{36}-11120136162\,{y}^{35}\\
&&-8948546308\,{y}^{34}-6939692682\,{y}^{33}-14\,{y}^{4}-5\,{y}^{3}-950\,{y}^{8}-365\,{y}^{7}\\
&&-132\,{y}^{6}-42\,{y}^{5}-5916\,{y}^{10}-2398\,{y}^{9}-33708\,{y}^{12}-14290\,{y}^{11}-77684\,{y}^{13}\\
&&-831014\,{y}^{16}-385741\,{y}^{15}-234813592\,{y}^{52}-4562339774\,{y}^{47}.
\end{eqnarray*}

\end{document}